\documentclass[a4paper,11pt]{article}
\usepackage[english]{babel}
\usepackage{amsmath,amsfonts,amssymb,pifont}
\usepackage{textcomp}
\usepackage{latexsym}
\usepackage{theorem,srcltx}
\usepackage{graphicx,float}
\usepackage{mathpazo}

\usepackage{rotate,graphics,epsfig,float}
\theoremstyle{break}
\usepackage{amsmath,verbatim,a4wide}

\usepackage{color}

\definecolor{blue}{rgb}{0,0,1}
\definecolor{red}{rgb}{1,0,0}
\definecolor{purple}{rgb}{1,0,1}

\long\def\red#1\endred{\textcolor{red}{#1}}
\long\def\blue#1\endblue{\textcolor{blue}{#1}}
\long\def\purple#1\endpurple{\textcolor{purple}{#1}}

\usepackage{srcltx}
\newcommand\ld{\lambda}
\newcommand\Ld{\Lambda}
\newcommand\prf{\par\noindent%
\emph{Proof. }\ignorespaces}
\newcommand\qed{\hspace{\fill}$\square$\smallskip\par}
\newcommand\RR{{\mathbb R}} {\theorembodyfont{\rmfamily}
\newtheorem{prop}{Proposition}[section] } {\theorembodyfont{\rmfamily}
\newtheorem{lemma}[prop]{Lemma} } { \theorembodyfont{\rmfamily}
\newtheorem{corol}[prop]{Corollary} }

\title{The inhomogeneous Fermi-Pasta-Ulam chain, \\a case study of the
$1:2:3$ resonance}
\author{Roelof Bruggeman and Ferdinand Verhulst\\
Mathematisch Instituut, PO Box 80.010\\
3508TA Utrecht, Netherlands}

\begin{document}


\maketitle

\noindent
Running title: The inhomogeneous Fermi-Pasta-Ulam chain
\bigskip

\noindent
\begin{abstract}
The inhomogeneous Fermi-Pasta-Ulam chain is studied by identifying the mass
ratios that produce prominent resonances. This is a technically complicated
problem as we have to solve an inverse problem for the spectrum of the
corresponding linearized equations of motion. In the case of the
 inhomogeneous periodic Fermi-Pasta-Ulam chain with four particles each mass
 ratio determines a frequency ratio for the quadratic part of the
 Hamiltonian. Most prominent frequency ratios occur but not all. In general
 we find a one-dimensional variety of mass ratios for a given frequency
 ratio.

For the resonance $1:2:3$ a small cubic term added to the Hamiltonian leads
to a dynamical behaviour that shows a difference between the case that two
masses are equal and the more general case of four different masses. For
two equal masses the normalized system is integrable and chaotic behaviour
is small-scale. In the transition to four different masses we find a
Hamiltonian-Hopf bifurcation of one of the normal modes leading to complex
instability and Shilnikov-Devaney bifurcation. The other families of
short-periodic solutions can be localized from the normal forms together
with their stability characteristics. For illustration we use action
simplices and the behaviour with time of the $H_2$ integral of the normal
forms.
\end{abstract}

\section{Introduction} The Fermi-Pasta-Ulam (FPU) chain or lattice is an $n$
degrees-of-freedom (dof) Hamiltonian system that models a chain of
oscillators with nearest-neighbour interaction, see \cite{FPU55} and
\cite{F92}. We will describe the model in section \ref{model}, see also
\cite{J91}. There exists a huge amount of literature on the FPU chain but
nearly always regarding the case of equal masses, sometimes called the
mono-atomic case. In this paper we will outline a research program to study
the inhomogeneous case where the masses are different. An inhomogeneous
nonlinear lattice with nearest neighbour interaction is studied in
\cite{UM} with emphasis on energy control. It is understandable that only a
few results were obtained for inhomogeneous lattices as the choice of
inhomogeneities, the masses of the lattice, seems to be arbitrary. We will
solve this arbitrariness by focusing on the presence of resonances induced
by the choice of masses. After referring to some basic material on
Hamiltonians and normal forms we formulate in section 2 the periodic FPU
$\alpha$ chain with arbitrary (positive) masses. In such a $n$
degrees-of-freedom system there exists a momentum integral that enables us
to reduce to a $n-1$ dof system. An inverse problem is considered in
section 3: how do we find mass distributions producing prominent resonances
in the spectrum induced by $H_2(p, q)$? This involves the analysis of the
inverse map of the vector of mass distribution to the vector of positive
eigenvalues of an associated coefficient matrix. This problem is solved in
section 3 for the cases of 3 and 4 particles; in the latter case it turns
out that of the four 1st order resonances that exist in general (see for
the terminology \cite{SVM})
3 exist, of the 12 possible 2nd order resonances 10 exist in this FPU
problem. In section 4 we focus on the $1:2:3$ resonance that arises for a
one-dimensional variety of mass ratios. It turns out that for one
particular combination of mass ratios, the normal form of the nonlinear
system is integrable. Moving from this particular case into the variety of
mass ratios, one of the periodic solutions shows Hamilton-Hopf bifurcation
that corresponds with Shilnikov-Devaney bifurcation in this Hamiltonian
system and produces a chaotic normal form.

The appendix contains general statements on the relation between mass ratios
and the spectrum induced by $H_2(p, q)$ that can be useful for future
research. Table 3 summarises the instructions for the case of 4 particles.
It is shown that for a given $n$-dimensional eigenvector characterizing the
FPU chain, all positive solutions of an $n$-dimensional mass distribution
are in a compact subset of ${\mathbb{R}}^n$. This subset is empty in some
cases, for instance the important $1:1: \ldots :1$ resonance does not arise
for the periodic FPU $\alpha$ chain with four or more particles.

\subsection{Hamiltonian formulation} For an autonomous Hamiltonian system
with $n$ degrees-of-freedom (dof), $n$ independent integrals suffice for
integrability, in that case there will be no chaotic motion in such a
system. However, in general, Hamiltonian systems with two or more
degrees-of-freedom (dof). are non-integrable. In many cases, this
 phenomenon was identified with homoclinic chaos as predicted by Poincar\'e
 in the nineteenth century, see \cite{PMC}, vol. 3; for a description see
\cite{FV12}, sections 5.4 and 9.3.

In the seventies of last century, a number of scientists started with the
computation and analysis of normal forms of general Hamiltonian systems
near equilibrium. Introductions and surveys of results can be found in
\cite{SVM}, chapter 10 and \cite{FV15}. One starts with an $n$
degrees-of-freedom system with Hamiltonian $H(p, q)$ that can be expanded
near equilibrium to a certain order as:
\[ H(p, q) = H_2(p, q) + H_3(p, q) + \ldots + H_m(p, q) + \ldots \]
The index $m$ indicates the degree in the variables $(p, q)$ of the
homogeneous polynomial $H_m$. Sometimes, other coordinate systems are
useful, for instance action-angle variables $\tau_i, \phi_i, i=1 \ldots n$.
The normal form technique was developed by Poincar\'e, Birkhoff and modern
scientists using analytic and algebraic tools. A basic element is that the
resonances that exist near equilibrium produce resonant terms that are kept
in the normal form while the non-resonant terms are averaged away. Such a
normal form $\bar{H}(p, q)$ does generally not converge when
$m \rightarrow \infty$, but a finite expansion contains already a lot of
quantitative and qualitative information. The respective polynomials
$H_m(p, q),  m= 3, 4, \ldots$ are after normalization indicated by
$\bar{H}_m(p, q)$. \\
Usually, consideration of a neighbourhood of stable equilibrium is made
explicit by scaling with a small positive parameter $ \varepsilon$
$p \rightarrow \varepsilon p, \, q \rightarrow \varepsilon q$ and dividing
the resulting Hamiltonian by $\varepsilon^2$. The terms $H_m,\, m \geq 2$
have the coefficient $\varepsilon^{m-2}$. As the normalization is
 canonical,
\[ \bar{H}(p, q) = H_2(p, q) + \varepsilon \bar{H}_3(p, q) + \ldots +
\varepsilon^{m-2} \bar{H}_m(p, q) \]
is the Hamiltonian integral of the normal form to degree $m$, whereas,
because of the normal form technique, also $H_2(p, q)$ is an integral of
the normal form system. This means that two degrees-of-freedom Hamiltonian
normal forms are always integrable, they contain no chaos.

\begin{enumerate}
\item Hamiltonian normal forms of three or more dof are generally
non-integrable; for a recent survey see \cite{FV15}. In the present paper
we will explore to some extent the presence of first and second order
resonances for the inhomogeneous FPU problem. The results for the
occurrence of resonances will be summarized in table \ref{tab-res}.

\item The presence of prominent (first and second order) resonances suggests
a research programme outlined in subsection \ref{rp}. For illustration and
as a start we will study the $1:2:3$ resonance for the inhomogeneous FPU
problem in the case of four oscillators in a so-called periodic
$\alpha$-chain.

\end{enumerate}

It is standard to use action-angle variables
$\tau_i, \phi_i, i= 1, \ldots, n$ near stable equilibrium:
\begin{equation} \label{aa}
 p_i= \sqrt{2 \tau_i} \cos \phi_i,\, q_i= \sqrt{2 \tau_i} \sin \phi_i,\, i=
 1, 2, \dots, n.
 \end{equation}
The equations of motion in action-angle variables are after transforming
$p, q \rightarrow \tau, \phi$:
\[ \dot{\tau}_i = - \frac{\partial H}{\partial \phi_i},\, \dot{\phi}_i =
\frac{\partial H}{\partial \tau_i},\,i= 1, 2, \dots, n.;\, n \geq 3 \]
However, using action-angle variables, special care is needed near the
normal modes. After giving arguments in the next subsection, we will use
co-moving coordinates in the coordinate planes or whenever an action is
near zero. Also, we will often use polar coordinates instead of
action-angle variables for orbits in general position; although such
transformations are not canonical, they preserve the energy, are easier to
establish the effect of resonances and most importantly, they produce
qualitatively and quantitatively mathematically equivalent results to
action-angle variables (for estimates see \cite{SVM}).

In the sequel, a periodic solution should be understood as a periodic
solution for a fixed value of the energy (iso-energetic solution), so
actually it corresponds for the full Hamiltonian system with a family of
periodic solutions parameterized by the energy.

\subsection{On normal forms and Floquet exponents} Normal form computations
for Hamiltonian systems can be carried out in various ways. Apart from
efficiency, the main point is to keep the system energy-preserving and
preferably canonical. Using for instance action-angle coordinates
(\ref{aa}) or amplitude-phase coordinates one can perform averaging over
the angles or explicitly time to obtain a first-order normal form. One may
consult \cite{SVM} for more details. An introductory text is \cite{V05},
chapters 11 and 12.

In section \ref{123} we will analyze periodic $\alpha$-chains (FPU chains
where the Hamiltonian is truncated after the cubic terms), containing the
$1:2:3$ resonance with main objective to investigate the stability of the
short-periodic solutions on the energy manifold and the integrability of
the normal form. This is highly relevant for the characterization of the
chaotic dynamics of the system but, as mentioned above, it raises special
problems. In the cases of vanishing actions or amplitudes, for instance
when studying normal modes, the procedure will be as follows (see also
section \ref{bal}).

Starting with the equations of motion, we will use co-moving coordinates
(see for instance transformation
(11.9-10) in chapter 11 of \cite{V05}) to obtain a first order normal form.
This normal form is used to localize the short-periodic solutions; the
normal form conserves the energy but the transformation is not canonical.
We will use averaging-normalization as it yields rigorous approximation
results (see \cite{SVM}), the results are qualitatively and quantitatively
precise. The same holds when we use polar coordinates outside the
coordinate planes.\\
In section \ref{123}, the short-periodic solutions can be computed
explicitly. The next step is then to linearize near the periodic solutions
and to determine the Floquet exponents for which we have to study coupled
Mathieu equations. This is still a formidable task, but we can obtain a
first order approximation of the exponents by normalizing the coupled
Mathieu equations. This will give a number of stability results in section
\ref{123}.

\subsection{Outline of a research programme} \label{rp}
The original Fermi-Pasta-Ulam chain \cite{FPU55} consists of $n$ oscillators
of equal mass with nearest-neighbour interaction; the chain will be
described in the next section. Thousands of papers and a number of
conferences were devoted to FPU chains, its stimulus for nonlinear science
has been enormous. Among the various problem formulations there was one
(nearly always) constant element: the masses of the chain were taken equal.
We will present here arguments for considering other mass distributions.

In a neighbourhood of equilibrium, the spectrum of the linear part of the
equations of motion plays a crucial part regarding the nature of the
ensuing dynamics, see for instance \cite{SVM} or \cite{FV15}. Considering
inhomogeneous mass distributions in FPU chains, one can produce a great
many different spectra induced by $H_2$. Each of these cases may produce
different dynamics in the corresponding FPU chain. In section \ref{br} we
will consider resonant spectra for the case of three and more extensively
four particles with periodic boundary conditions i.e. chains where the
first and the last oscillator are identified. For the case of four
particles we will focus on the rich dynamics of the $1:2:3$ resonance. An
outline of possible further research follows here:
\begin{enumerate}
\item According to table \ref{tab-res} regarding the case of four particles,
we also have to study two first order resonances ($1:2:1$ and $1:2:4$) and
ten second order resonances. Also, higher order resonances may be
worthwhile to investigate. Special attention should be given to the $1:1:2$
and $1:1:3$ cases as only four special mass ratios produce these
 resonances. In such a case degenerations may arise so that we have to
 consider detuning phenomena, see~\cite{SVM}.

\item Cases of five and more particles will present many more problems.
\item The present study is restricted to so-called periodic $\alpha$-chains.
Including quartic terms in the Hamiltonian ($\beta$-chains) and considering
lattices with fixed begin- and end-point will produce new results.
\item The study presented here and possibly future studies will throw light
on qualitative and quantitative differences between systems in
nearest-neighbour interaction and non-local interaction, a topic that is
relevant for plasma physics and stellar dynamics.
\end{enumerate}

\section{The Fermi-Pasta-Ulam chain} \label{model}
The FPU-chain with periodic boundary conditions has been a topic for many
studies. It describes a model for nonlinear interaction of identical point
masses moving on a circle with nearest-neighbour coupling. Numerical
integrations in the early 1950s showed that the expectation by physicists
of thermalisation by energy transport was not correct. Putting all the
energy originally in one mode, it was observed that this energy was shared
by only a few other modes. Nice introductions can be found in \cite{F92}
and \cite{J91}.

For the mono-atomic case of the original periodic FPU-problem (all masses
equal) it was shown in \cite{RV} for up to six degrees-of-freedom (dof) and
much more general in \cite{BR}, that the corresponding normal forms are
governed by $1:1$ resonances and that these Hamiltonian normal forms are
integrable. This explains the recurrence phenomena near equilibrium.

We will drop the original assumption of identical (mono-atomic) particles to
consider the periodic FPU-problem again. For $n$ particles with mass
$m_j>0$, position $q_j$ and momentum $p_j= m_j \dot{q}_j, j= 1 \ldots n$,
$\varepsilon \geq 0$ a small parameter, the Hamiltonian is of the form:
\begin{equation} \label{HAMFPU}
H(p, q)= \sum_{j=1}^n (\frac{1}{2m_j}p_j^2 + V(q_{j+1}-q_j)) \,\,
{\rm{with}}\,\, V(z)= \frac{1}{2}z^2+ \varepsilon \frac{\alpha}{3}z^3 +
\varepsilon^2 \frac{\beta}{4}z^4.
\end{equation}
 The quadratic part of the Hamiltonian is not in diagonal form; for
 $n=3, 4 \ldots$ the linearized equations of motion can be written as:
\begin{eqnarray} \label{lineqs}
\begin{cases}
m_1 \ddot{q}_1 +2q_1-q_2 -q_n & = 0,\\
m_2 \ddot{q}_2 +2q_2-q_3 -q_1 & = 0,\\
m_3 \ddot{q}_3 +2q_3-q_4 -q_2 & = 0,\\
\ldots \,\, \,\,\, \ldots & = 0, \\
m_n \ddot{q}_n +2q_n-q_1 -q_{n-1} & = 0.
\end{cases}
\end{eqnarray}
We can write for the quadratic part of $H(p, q)$:
\begin{equation} \label{quadr}
 H_2 = \frac{1}{2}p^TA_n p + \frac{1}{2} q^TC_n q,
\end{equation}
with $A_n$ the $n\times n$ diagonal matrix with at position $(i,i)$ the
value $m_i^{-1} =: a_i$, $C_n$ is an $n \times n$ matrix. For an analysis
of the quadratic term $H_2(p, q)$ we need to know the eigenvalues of
$A_nC_n$. The relation between the eigenvalues of $A_nC_n$ and the
eigenvalues of the matrix of coefficients of system \eqref{lineqs} will be
given below. Since the null space of $C_n$ has dimension one, the matrix
$A_nC_n$ has an eigenvalue $0$ corresponding to a (translational) momentum
integral. It will turn out that the other eigenvalues of $A_nC_n$ are
positive, as expected. For a given set of masses, the calculation of the
remaining eigenvalues corresponding with the frequencies of the linearized
system is easy, but we are faced with another, an inverse problem. To focus
ideas, suppose that $n=4$. The presence of the momentum integral implies
that we have to consider a three degrees-of-freedom (dof) Hamiltonian
problem. We know, see for instance \cite{SVM} chapter 10 or \cite{FV15},
that the first order resonances are $1:2:1, 1:2:2, 1:2:3$ and $1:2:4$. The
question is then if and how we can choose the masses so that these
prominent resonances are present. Of course, this problem will be more
formidable if $n>4$. In the next section we determine for $n=4$ the ratios
of masses that produce the resonance $1:2:3$. The approach works equally
well for other prescribed rations of eigenvalues, as we discuss in the
Appendix. Prominent resonances for $n > 4$ can be found but a systematic
study of these cases poses a difficult open algebraic problem.

\section{The spectrum induced by $H_2$} \label{br}

After a number of general considerations we will give details for the cases
of three and four particles. The first case is rather trivial as far as the
spectrum goes, the case of four particles is already quite complicated.
 Here we mention the main facts that we need in the later sections. In the
Appendix we will give more details.

\subsection{The matrix for inhomogeneous FPU-lattices and its
eigenvalues}\label{sect-mev}
The linear system (\ref{lineqs}) can be written as
\begin{equation}\label{MAC}
\begin{pmatrix}
\dot q \\ \ddot q
\end{pmatrix} = M
\begin{pmatrix}
 q \\ \dot q
\end{pmatrix}\,, \text{ where } M=
\begin{pmatrix} 0 & I_n \\
-A_nC_n & 0
\end{pmatrix}\,,
\end{equation}
where the matrix $A_n$ is a diagonal matrix with the inverse masses
$m_j^{-1}=:a_j$ on the diagonal, and where the matrix $C_n$ has elements
$2$ on the diagonal, and $-1$ at positions $(i,i+1)$ and $(i,i-1)$, with
the indices taken modulo~$n$. For instance,
\[ C_5= \begin{pmatrix}
2&-1&\hphantom{-}0&\hphantom{-}0&-1\\
-1&\hphantom{-}2&-1&\hphantom{-}0&\hphantom{-}0\\
\hphantom{-}0&-1&\hphantom{-}2&-1&\hphantom{-}0\\
\hphantom{-}0&\hphantom{-}0&-1&\hphantom{-}2&-1\\
-1&\hphantom{-}0&\hphantom{-}0&-1&\hphantom{-}2
\end{pmatrix}\,. \]
(This matrix turns up elsewhere in mathematics. It is the affine Cartan
matrix of the completed root system $\bar A_n$. See eg.\ \cite[D�f. 3 in
1.5 of Chap. 6, and Planche I]{Bo}.)

The $(2n)\times(2n)$ matrix $M$ has a double eigenvalue $0$, corresponding
to the momentum integral
\begin{equation} \label{moment}
\sum_{j=1}^n m_j \dot{q}_j = \rm{constant}.
\end{equation}
In the sequel we will choose the case of vanishing momentum integral which
is not a restriction of generality. If $\ld$ is a positive eigenvalue of
$A_nC_n$, then $i\sqrt{\ld}$ and $-i\sqrt{\ld}$ are eigenvalues of $M$,
corresponding to frequences of eigenmodes of the linearized system. So it
is useful to collect results concerning the eigenvalues of $A_nC_n$.

\begin{prop}\label{prop-ev}For $n=  3, 4 \ldots$ the matrix $A_nC_n$ has one
eigenvalue $0$ and $n-1$ positive eigenvalues
$\ld_1,\ldots,\allowbreak\ld_{n-1}$, possibly coinciding. If eigenvalues
coincide the corresponding eigenspace has maximal dimension.
\end{prop}
\prf Since the $a_j=m_j^{-1}$ are positive, the matrix $A_n^{1/2}$ is
well-defined. The symmetric matrix $A_n^{1/2}C_nA_n^{1/2}$ has real
eigenvalues, and the algebraic and geometric multiplicities of eigenvalues
coincide.

If $y$ is an eigenvector of $A_nC_n$ with eigenvalue $\ld$, then
\[ \ld \sum_i a_i^{-1} \, y_i^2 = 2\sum_i y_i^2 - 2 \sum_i y_i\,
y_{i-1}\,.\]
(Indices taken modulo $n$.) So
\[ \sum_i (2-\ld/a_i) y_i^2 =2 \sum_i y_1 y_{i-1}\,.\]
With Schwarz's inequality this implies $\ld \geq 0$. Equality occurs only if
the vectors $(y_i)$ and $(y_{i-1})$ are positive multiples of each other,
which occurs only for multiples of $(1,1,\ldots,1)$. \qed

\medskip
For the investigation of the linearized problem we need to understand the
map $\RR_{>0}^n \rightarrow \RR_{>0}^{n-1}$, from a vector
$(a_1,\ldots, a_n)$ of inverse masses to a vector
$(\ld_1,\ldots,\ld_{n-1})$ of positive eigenvalues. The order of the
eigenvalues is not determined, so we have, more precisely, a map
 $\rho_n: \RR_{>0}^n \rightarrow S_{n-1}\backslash \RR_{>0}^{n-1}$, with the
 action of the symmetric group $S_{n-1}$ on the coordinates. For the
 linearized inhomogeneous FPU-chain described by system (\ref{lineqs}), the
 dihedral group $D_n$ with $2n$ elements permutes the coordinate $q_j$
 (generated by a shift and a reflection). This transforms system
 (\ref{lineqs}) into an equivalent system. Another symmetry is by scaling:
 $\rho\bigl( t(a_1,\ldots,a_n)\bigr) = t \,\rho(a_1,\ldots,a_n)$ for $t>0$.

To investigate the correspondence between eigenvalues and inverse masses we
use the equality
\[ \det\bigl( A_n C_n - \ld \,I_n) = -\ld \prod(\ld_j-\ld)\,,\]
for $(\ld_1,\ldots,\ld_{n-1})=\rho(a_1,\ldots,a_n)$. This leads to
equalities
\begin{equation}\label{pldrel} p_j (A_n) =
e_{n-j}(\{\ld_1,\ldots,\ld_{n-1}\})\qquad
(1\leq j \leq n-1)\,,
\end{equation}
with the elementary symmetric functions $e_k$ and homogeneous polynomials
$p_j(A_n)$ in the $a_j$ of degree $n-j$. This describes the structure of
the set of diagonal matrices $A_n$ for a prescribed spectrum of $A_n C_n$.
 It is the set of points with positive coordinates in an algebraic set in
  $\mathbb{C}^n$ which is the intersection of $n-1$ hyperplanes given by
  equations of degree $1, \, 2, \ldots, n-1$.

In particular,
\begin{equation}\label{pe-eq} p_{n-1}(A_n) = 2\sum_i a_i\,,\quad
p_{n-2}(A_n) = \sum_{1\leq i <j\leq n} c_{i,j} \, a_ia_j\,,\end{equation}
with $c_{i,j}=3$ if $i-j=\pm 1\bmod n$, and $c_{i,j}=4$ otherwise. All
$p_j(A_n)$ are invariant under the action of the dihedral group $D_n$ on
the coordinates~$a_j$.

In subsection \ref{sect-ell} of the Appendix we'll prove relation
\eqref{pe-eq}. We will also show that all real solutions $(a_1,\ldots,a_n)$
for a given eigenvalue vector $(\ld_1,\ldots,\ld_{n-1})$ are in a compact
subset of $\RR^n$. This subset may be empty. For all $n\geq 4$ the
$1:1:\cdots:1$
 resonance does not occur for any mass distribution.

\subsection{The case of three particles} For $n=3$ the determination of the
eigenvalues for given inverse masses amounts to solving the quadratic
equation
\[ \ld^2 -2(a_1+a_2+a_3) \ld + 3(a_1a_2+a_1a_3+a_2a_3)=0,\]
which has positive solutions.

\begin{figure}[bf,tp]
\begin{center}
\includegraphics[width=6cm]{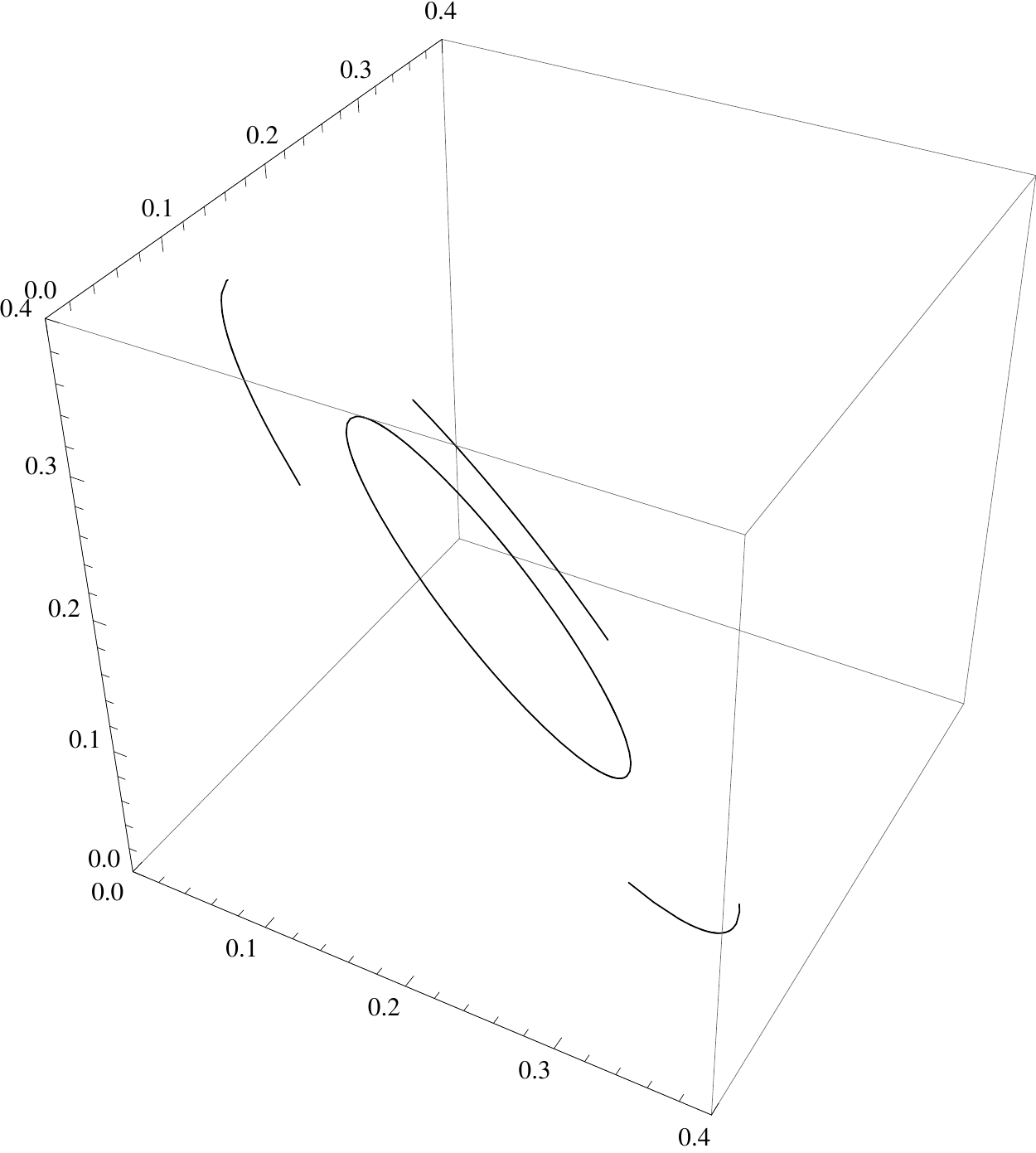}
\end{center}
 \caption{Solutions sets of inverse masses $(a_1,a_2,a_3)$ for the FPU-chain
with three particles. For the eigenvalue ratio $\ld_1/\ld_2=2$ of matrix
$A_nC_n$ the solution set is compact; it is an ellipse in $\RR_{>0}^3$. For
the eigenvalue ratio $\ld_1/\ld_2=4$ the solutions are on a larger ellipse
in $\RR^3$, which intersects $\RR_{>0}^3$ in three open curves.}
\label{ell3}\end{figure}

Conversely, for all choices $(\ld_1,\ld_2)$ of positive eigenvalues, values
of $a_1,a_2,a_3$ can be found such that $A_3C_3$ has eigenvalues $\ld_1$,
$\ld_2$ and $0$. If $\ld_1=\ld_2$ there is exactly one solution
$a_1=a_2=a_3=\frac 13 \ld$ (equal masses). If the eigenvalues have ratio
$\ld_1/\ld_2>1$ then the corresponding points $(a_1,a_2,a_3)$ in $\RR^3$
form an ellipse. This ellipse may or may not be contained in the positive
octant. See figure~\ref{ell3}.

\subsection{The case of four particles} In the case $n=4$ we use the scaling
to restrict our further investigation to eigenvalues satisfying
$\ld_1+\ld_2+\ld_3=1$. From \eqref{pe-eq} we obtain three equations for a
given vector $(\ld_1,\ld_2,\ld_3)\in \RR_{>0}^3$:
\begin{equation}\label{eq4}
\begin{aligned}
4(a_1a_2a_3+a_2a_3a_4+a_3a_4a_1+a_4a_1a_2) &= \ld_1\ld_2\ld_3 =: \xi\,,\\
3(a_1a_2+a_2a_3+a_3a_4+a_4a_1)+4(a_1a_3+a_2a_4)&= \ld_1\ld_2
+\ld_2\ld_3+ \ld_3\ld_1 =: \eta\,,\\
2(a_1+a_2+a_3+a_4)&= 1\,.
\end{aligned}
\end{equation}
We call the set of $(a_1,\ldots,a_4)\in \RR_{>0}^4$ satisfying these
relations the \emph{fiber} of $(\xi,\eta)\in \RR_{>0}^2$. In subsection
\ref{sect-ell} we will give a precise characterization of the set of
$(\xi,\eta)$ for which the fiber is non-empty.

The resonances deserve special attention. A resonance $(n_1:n_2:n_3)$ in the
linearized system \eqref{lineqs} corresponds to an eigenvalue vector of
$A_4C_4$ with the ratios $(n_1^2 : n_2^2 : n_3^2)$. We considered all
resonances of order one and two, and obtained the results in
Table~\ref{tab-res}. As noted in subsection \ref{rp}, the resonances
$1:1:2$ and $1:1:3$ need special attention.
\begin{table}[bf, htp]
\renewcommand\arraystretch{1.2}
\begin{center}
\begin{tabular}{|cl|cl|}\hline
\text{ratio}&\text{fiber}& \text{ratio}&\text{fiber}
\\ \hline
$\bigl(1:1:\sqrt 2\bigr)$& \multicolumn{3}{l|} {one point
(classical case $A_4=I_4$)}\\
\hline
\multicolumn{4}{|c|}{resonances of order $1$}\\
$(1:1:2)$& four points& $(1:2:2)$& empty\\
$(1:2:3)$& four open curves & $(1:2:4)$&12 open curves\\ \hline
\multicolumn{4}{|c|}{resonances of order $2$}\\
$(1:1:1)$&empty& $(1:1:3)$& four points \\
$(1:2:5)$ & 12 open curves & $(1:2:6)$ & 12 open curves \\
$(1:3:3)$&empty& $(1:3:4)$& four open curves \\
$(1:3:5)$& four open curves & $(1:3:6)$ & 12 open curves \\
$(1:3:7)$& 12 open curves & $(1:3:9)$ & 12 open curves \\
$(2:3:4)$ & two compact curves & $(2:3:6)$ & two compact curves \\ \hline
\end{tabular}
\end{center}
\caption{Fibers of resonances} \label{tab-res}
\end{table}

\subsection{The resonance $(1:2:3)$}\label{sect-comp123}Here we consider the
resonance that is the subject of study in the next section.

By scaling we arrange $\ld_1=\frac 9{14}$, $\ld_2=\frac 27 $,
$\ld_3=\frac1{14}$ to satisfy the last equation in~\eqref{eq4}. The first
equation
$e_3\bigl( \{a_1,a_2,a_3,a_4\}\bigr) = \frac 14 \xi =\frac 9{2744}$
contains the third elementary symmetric polynomial in the $a_j$. The middle
equation has the form
\begin{equation} 3 q_2(a_1,a_2,a_3,a_4) + 4 q_1(a_1,a_2,a_3,a_4)=\eta
=\frac14\,,
\end{equation}
with polynomials
\begin{equation}
q_1(a_1,a_2,a_3,a_4)=a_1a_3+a_2a_4\,,\quad q_2(a_1,a_2,a_3,a_4)=
a_1a_2+a_2a_3+a_3a_4+a_4a_1 \,,\end{equation}
which are invariant under the dihedral group~$D_4$. We solve the system of
equations~\eqref{eq4} by prescribing values for these two polynomials. So
we work with $q_1(a_1,a_2,a_3,a_4)=\eta_1$, $q_2(a_1,a_2,a_3,a_4)=\eta_2$,
where $\eta_1$ and $\eta_2$ satisfy $4\eta_1+3\eta_2=\eta=\frac 14$. Since
both are to be positive this requires $0<\eta_2<\frac 1{12}$.

Now we have four equations for the four unknown quantities $a_j$, and may
expect a discrete set of solutions for each appropriate value of the
parameter $\eta_2$.

First we consider the positive quantities $s_{13}=a_1+a_3$, and
$s_{24}=a_2+a_4$. They satisfy
\[ s_{13}+s_{24} = \frac12\,, \quad s_{13}s_{24}=\eta_2\,. \]
Hence we may take
\begin{equation}\label{sch}
s_{13},\,s_{24}= \frac14\mp \frac14 \sqrt{1-16\eta_2}\,.
\end{equation}
This has positive values only if $0<\eta_2<\frac1{16}$. In the case of the
resonance $(1:2:3)$ it turns out to be convenient to write
$\sqrt{1-16 \eta_2}= \frac{5-u}7$, and to use $u\in (-2,5]$ as the
parameter. We take
\begin{equation}
s_{13}= \frac14-\frac14\,\frac{5-u}7 = \frac{2+u}{28}\,,\qquad s_{24} =
\frac{12-u}{28}\,.
\end{equation}

For $p_{13}=a_1a_3$ and $p_{24}=a_2a_4$ we find the relations
\[ p_{13}+p_{24}=\eta_1\,, \quad s_{13}p_{24}+s_{24}p_{13}=\frac\xi4\,. \]
If $u=5$ we have $s_{13}=s_{14}$. Then $\eta_2=\frac 1{16}$, and
$\eta_1=\frac14\left( \frac 14 - \frac3{16}\right)=\frac1{64}$. Since
$\xi=\frac 9{686}\neq \frac1{64}$ this does not lead to a solution. So we
can proceed with $-2<u<5$ and find solutions
\begin{equation}
p_{13}= \frac{\xi/4-s_{13}\eta_1}{\frac12\sqrt{1-16\eta_2}} \,,\quad
p_{24}=\eta_1-p_{13}\,.
\end{equation}
These quantities should be positive. To have $p_{13}>0$ we need to restrict
$u$ to the interval $(0,u_1)$ with
\begin{equation}
u_1= \frac 83 - \frac 23\sqrt[3]{19} \approx .887732\,.
\end{equation}

Now we have $a_1+a_3=s_{13}$ and $a_1a_3=p_{13}$. This gives a quadratic
equation for $a_1$ and $a_3$, with discriminant
\begin{equation}
\frac{(16-u)(6-u)u}{1568(5-u)}\,.
\end{equation}
So there are real solutions that coincide for $u=0$.
\begin{equation}\label{a13}
a_1,a_3= \frac{2+u}{56}\mp \frac{\sqrt2}{112}\,
\sqrt{\frac{u(6-u)(16-u)}{5-u}} \,,
\end{equation}
where we take the minus sign for $a_1$. Both functions are positive for
$u\in [0,u_1)$. At $u=0$ they have the same value. The limit
$\lim_{u\rightarrow u_1}a_1(u)$ is zero, corresponding to the extreme case
of an infinite mass.

\begin{figure}[bf,htp]
\begin{center}
\includegraphics[width=10cm]{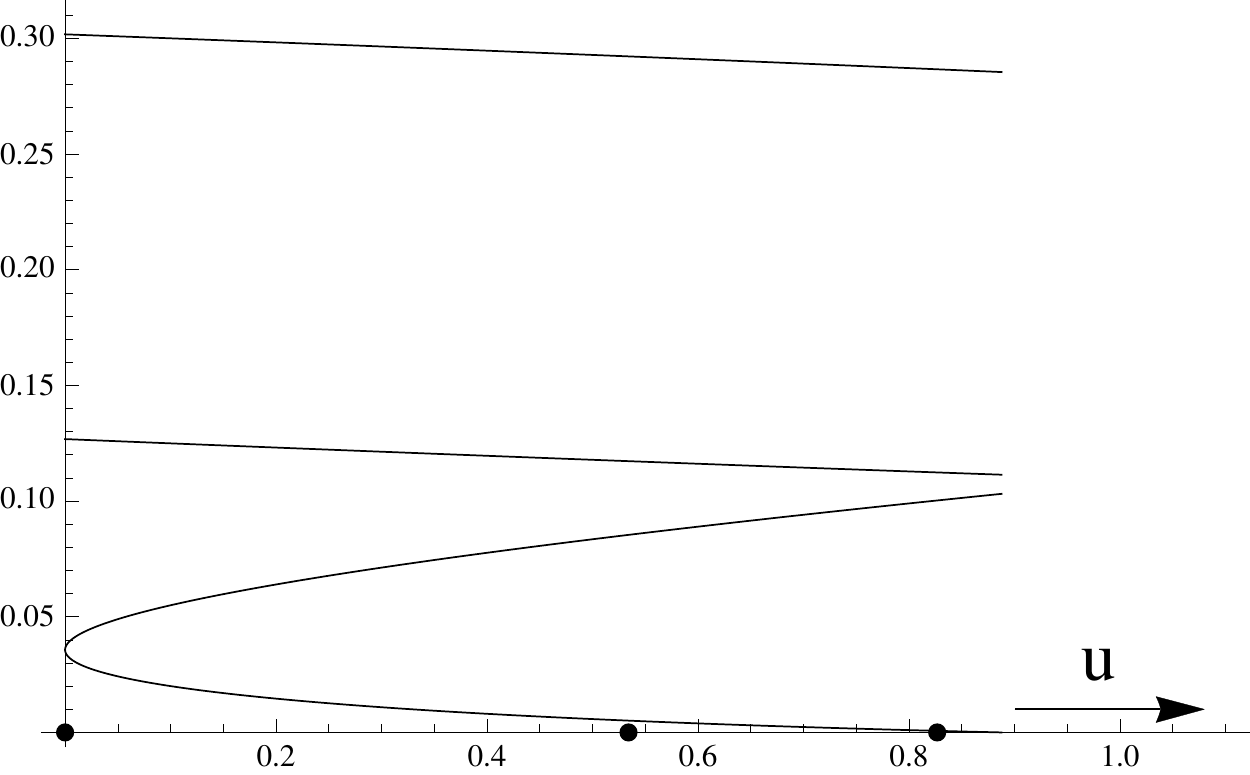}
\end{center}
\caption{One branch of the fiber for the resonance $(1:2:3)$ is given by the
functions $a_1\leq a_3<a_2<a_4$ on the interval $[0,u_1)$. (Horizontal
axis: parameter $u$; vertical axis: values of $a_j(u)$.) The three dots on
the horizontal axis correspond to the values $0,\, 0.534105$ and $0.826713$
of the parameter $u$ for which we will carry out simulations in the next
section (the cases $0$, $1$ and $2$). }\label{fig-a}\end{figure}

The discriminant $s_{24}^2-4p_{24}$ of the equation for $a_2$ and $a_4$ is
positive for all $u\in [0,u_1)$. and leads to two solutions that are
unequal for all $u\in [0,u_1)$.
\begin{equation}\label{a24}
a_2,a_4= \frac{12-u}{56} \mp \frac1{56\sqrt 2}\,
\sqrt{\frac{(6+u)(4-u)(10-u)}{5-u}}\,,
\end{equation}
with the minus sign for $a_2$. Figure~\ref{fig-a} gives a plot.

In the course of the proof we have made three sign choices, in \eqref{a13}
and \eqref{a24}, and in \eqref{sch}. We get all solutions when we let the
dihedral group $D_4$ act on the solutions that we constructed.
\subsection{Illustration of the fiber}The equations \eqref{a13} and
\eqref{a24} describe a curve
$u\mapsto \bigl(a_1(u),a_2(u),a_3(u),a_4(u)\bigr)$ in $\RR_{>0}^4$
corresponding to a one-parameter family of solutions for the inverse
masses. To illustrate it we use the second and last equation in
\eqref{eq4}, which describe an ellipsoid in the hyperplane
$a_1+a_2+a_3+a_4=\frac12$. In subsection \ref{sect-spher} in the appendix
we'll describe this ellipsoid in a more explicit way. The first equation in
\eqref{eq4} produces an intersection with this ellipsoid in some curves.
The points with positive coordinates in this intersection form the fiber.

\begin{figure}[bf,htp]
\begin{center}
\includegraphics[width=10cm]{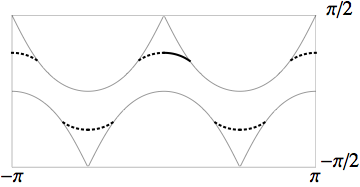}
\end{center}
 \caption[]{The fiber for the resonance $(1:2:3)$ is contained in an
ellipsoid, which we describe in spherical coordinates. (Horizontally the
azimuth $\phi$, and vertically the inclination $\psi$. See \eqref{pc1}.)
 The thick curve corresponds to the branch of the fiber in
figure~\ref{fig-a}. The dotted curves correspond to the translates of this
branch under the dihedral group $D_4$. The thin curves indicate the
boundary of the region corresponding to coordinates in $\RR_{>0}^4$.

The picture illustrates that the fiber for $(1:2:3)$ consists of four open
curves, and that \eqref{a13} and \eqref{a24} give a fundamental domain for
the action of the dihedral group on the fiber.
 }\label{fig-ell123}\end{figure}

On the ellipsoid we can use a system of spherical coordinates, mapping the
ellipsoid to the rectangle $[-\pi,\pi]\times[-\frac12\pi,\frac12\pi]$, with
boundary identifications. The image of the fiber under this map is given in
figure~\ref{fig-ell123}.

\subsection{Transformation of the Hamiltonian}\label{sect-hamtr}
We form the diagonal matrix $A_4(u)$ with diagonal elements $a_j(u)$,
$1\leq j \leq 4$. In the proof of Proposition~\ref{prop-ev} we noted that
$A_4(u)^{1/2}C_4 A_4(u)^{1/2}$ is a symmetric matrix (as long as
$u\in [0,u_1)$), so we can find an orthogonal matrix $U(u)$ such that
$A_4(u)^{1/2}C_4 A_4(u)^{1/2}= U(u) \Ld U(u)^T$, where $\Ld$ is the
diagonal matrix with diagonal elements $\frac 9{14}$, $\frac 27$,
$\frac1{14}$, and $0$. Then the transformation matrices
\begin{equation} K(u) = A_4(u)^{-1/2} U(u)\,, \quad L(u) = A_4(u)^{1/2}U(u)
\end{equation}
determine a symplectic transformation
\begin{equation}\label{pqyx} p=K(u) y\,,\quad q=L(u)x\,,
\end{equation}
which transforms the quadratic part in \eqref{quadr} of the Hamiltonian into
\begin{equation}\label{H2trf} H_2 = \frac12 y^T y + \frac12 x^T \Ld x =
\frac12\sum_{j=1}^4 \Bigl( y_j^2 + \ld_j x_j^2\Bigr)\,.
\end{equation}
This will produce the so-called quasi-harmonic form of the equations of
motion. To see that $H_2$ takes the form \eqref{H2trf} we need the
existence of an orthogonal matrix $U(u)$ diagonalizing
$A_4(u)^{1/2} \, C_4\, A_4(u)^{1/2}$. We do not need to know $U(u)$, $K(u)$
or $L(u)$ explicitly.

To transform the cubic and higher order terms of the Hamiltonian to
coordinates corresponding to the eigenmodes of the linearized system we
need to know the transformation matrix $L(u)$ explicitly. For any given
$u\in [0,u_1)$ it is no problem to do this numerically with {\sc
Mathematica} or {\sc Matlab}. It is nicer to have $U(u)$, and hence $L(u)$
and $K(u)$, symbolically in terms of the parameter $u$; see subsection
\ref{app-trm}. The explicit description of our choice of $L(u)$ will be
given in Table~\ref{L-expl} in the appendix.

For the cubic term we note that (with indices modulo~$4$)
\begin{equation} \frac 13 \sum_j \bigl(q_{j+1}-q_j\bigr)^3 = \sum_j
\bigl(q_{j+1}-q_{j-1}\bigr)q_j^2\,.
\end{equation}
The substitution $(q_1,\ldots,q_4)^T = L(u) \bigl( x_1,\ldots,x_4)^T$ gives
\begin{equation}
\begin{aligned}\label{H3dj}
H_3 &= \varepsilon\bigl( d_1(u)x_1^3+d_2(u)x_1^2x_2+d_3(u)x_1^2x_3
+d_4(u)x_2^2x_1+d_5(u)x_3^2x_1\\
&\quad
+d_6(u)x_1x_2x_3
+d_7(u)x_2^3+d_8(u)x_3^3+ d_9(u)x_3^2x_2+d_{10}(u)x_2^2x_3\bigr)\,,
\end{aligned}
\end{equation}
with the functions $d_j$ as indicated in Table~\ref{tab-dj}.

\begin{table}
\begin{align*}
d_1(u)&=\sqrt{u}\,\frac{27 \sqrt{4-u} \sqrt{6-u} \sqrt{10-u} (16-u)
   (u+6)}{35840 \sqrt{35} (5-u)} \,,\\
d_2(u)&= -\sqrt{u}\, \frac{3 \sqrt3 \sqrt{10-u} \sqrt{16-u} \sqrt{u+6}
   \left(3 u^2-30 u+52\right)}{4480 \sqrt{70}
   (5-u)}\,,\\
   d_3(u)&=-\frac{3 \sqrt{3} \sqrt{4-u} \sqrt{6-u} \sqrt{16-u} \sqrt{u+6}
   \left(3 u^2-30 u+160\right)}{35840\sqrt 7 (5-u)}\,,\\
   d_4(u)&=- \sqrt{u}\,\frac{\sqrt{4-u} \sqrt{6-u} \sqrt{10-u} \left(-3
   u^2+30 u+68\right)}{1120 \sqrt{35} (5-u)}\,,\\
d_5(u)&=-\sqrt{u}\,\frac{\sqrt{4-u} \sqrt{6-u} \sqrt{10-u} \left(3 u^2-30
   u+64\right)}{7168 \sqrt{35} (5-u)}\,,\\
 d_6(u)&= -\frac{-3 u^4+60 u^3-352 u^2+520 u+960}{2240\sqrt{14}
   (5-u)}\,,\\
 d_7(u)&= \sqrt{u}\,\frac{\sqrt{10-u} \sqrt{16-u} (6-u) (4-u)
   \sqrt{u+6}}{420 \sqrt{210} (5-u)}\,,\\
 d_8(u)&=u\,\frac{\sqrt{4-u} \sqrt{6-u} \sqrt{16-u} (10-u)
   \sqrt{u+6}}{21504 \sqrt{21} (5-u)}\,,\\
 d_9(u)&=- \sqrt{u} \, \frac{\sqrt{10-u} \sqrt{16-u}\sqrt{u+6} \left(u^2-10
   u+28\right)}{869 \sqrt{210} (5-u)}\,,\\
 d_{10}(u)&= \frac{\sqrt{4-u} \sqrt{6-u} \sqrt{16-u} \sqrt{u+6} \left(u^2-10
   u+20\right)}{1120 \sqrt{21} (5-u)}\,.
\end{align*}
\caption[ ]{Coefficients of the cubic term of the Hamiltonian, transformed
to eigenmodes, the so-called quasi-harmonic form. For $u=0$ only $d_3$,
$d_6$ and $d_{10}$ are non-zero. } \label{tab-dj}
\end{table}

\section{The $1:2:3$-resonance for the periodic $\alpha$-lattice ($n=4$)}
\label{123}
For any possible inhomogeneous FPU $\alpha$-chain with four dof we have the
system:
\begin{eqnarray} \label{nleqs}
\begin{cases}
\dot{q_1} &= v_1,\\
\dot{v}_1 & = [-2q_1+q_2 +q_4 - \varepsilon \alpha ((q_1-q_4)^2-(q_2-q_1)^2)
] a_1,\\
\dot{q_2} &= v_2,\\
\dot{v}_2 & = [-2q_2+q_3 +q_1 - \varepsilon \alpha ((q_2-q_1)^2-(q_3-q_2)^2
)] a_2,\\
\dot{q_3} &= v_3,\\
\dot{v}_3 & = [-2q_3+q_4 +q_2 - \varepsilon \alpha ((q_3-q_2)^2
-(q_4-q_3)^2)] a_3,\\
\dot{q_4} &= v_4,\\
\dot{v}_4 & = [-2q_4+q_1 +q_3- \varepsilon \alpha ((q_4-q_3)^2
-(q_1-q_4)^2)] a_4,\\
\end{cases}
\end{eqnarray}

The coefficient $\alpha$ has been retained for reference to the literature;
here we will take $\alpha=1$. If $a_1= \ldots =a_4=1$, we have the
classical periodic FPU problem with four particles; it was shown in
\cite{RV}, that in this case the normal form is integrable. The implication
is that for $\varepsilon$ small, chaos is negligible in this classical
case.

Apart from the Hamiltonian we have from \eqref{moment} as a second
(momentum) integral:
\begin{equation} \label{M}
m_1v_1+m_2v_2+m_3v_3+m_4v_4 = \rm{constant}.
\end{equation}
The presence of the momentum integral results in two zero eigenvalues of the
matrix M in eq. (\ref {MAC}), so by reduction we have to deal essentially
with a three dof system. \\

\medskip
According to table \ref{tab-res} the $1:2:3$ resonance is present among the
possible inhomogeneous FPU lattices. Fig. \ref{fig-a} gives one branch of
values of inverse masses $a_1, \ldots, a_4$ producing this resonance. All
vectors $(a_1, \ldots, a_4)$ are obtained by the action of the dihedral
group $D_4$ on the coordinates and the scaling
$(a_1, \ldots, a_4) \mapsto (ta_1, \ldots, ta_4)$ with $t>0$.

 Table \ref{tab-res} and fig. \ref{fig-a} show that the $1 : 2 : 3$
 -resonance appears in one case with relatively well-balanced masses, two of
which are equal. We denote this by case $0$; it will turn out in subsection
\ref{bal} that this case is quite special dynamically. The other cases are
less balanced regarding the masses. Case $0$ corresponds to $u = 0$; as $u$
increases (we have $0 \leq u < u_1$ with $u_1=.887732$), the masses get
less well-balanced, one of them tending to infinity. We study the dynamical
behaviour in subsection \ref{unbal}. For numerical simulations we have
singled out two more cases indicated in fig. \ref{fig-a}.

  The expression for the quadratic part of the Hamiltonian $H_2$ is:
\begin{equation} \label{H2}
H_2 = \frac{1}{2} \sum_{i=1}^4 \frac{v_i^2}{a_i}+ \frac{1}{2}
[(q_2-q_1)^2+(q_3-q_2)^2+(q_4-q_3)^2+(q_1-q_4)^2].
\end{equation}
$H_2$ is a first integral of the linear system (\ref{lineqs}), it is also a
first integral of the normal form of the full system (\ref{nleqs}). When
using $H_2$ from the solutions of the truncated normal form
\[ \bar{H}(p, q) = H_2(p, q) + \varepsilon \bar{H}_3(p, q), \]
we obtain an $O(\varepsilon)$ approximation of the (exact) $H_2(p(t), q(t))$
valid for all time; for a proof see \cite{SVM} chapter 10. Note that in the
equations we use the velocities instead of the momenta. Using the
expression $H_2(p(t), q(t))$ for the solutions of the full system
(\ref{nleqs}) shows the accuracy of the normal form and gives an impression
of the nature of the dynamics.

The normal form $\bar{H}_3(p, q)$, written in action-angle coordinates
\eqref{aa} or amplitude-phase coordinates
(see below), will contain certain combination angles corresponding with the
resonance. If $\bar{H}_3$ contains only one combination angle, we have an
additional integral of motion and the normal form $H_2+ \bar{H}_3$ is
integrable. In the case of two or more independent combination angles, we
have to investigate the
(non-)integrability of the normal form.

\medskip
To display the quantitative aspects of the solutions we have the possibility
of drawing an energy- or action-simplex or as an alternative to produce a
time series for explicit solutions or integrals of the normal forms. Both
techniques will be used.

\begin{figure}[htp]
\begin{center}
\resizebox{12cm}{!}{
\includegraphics{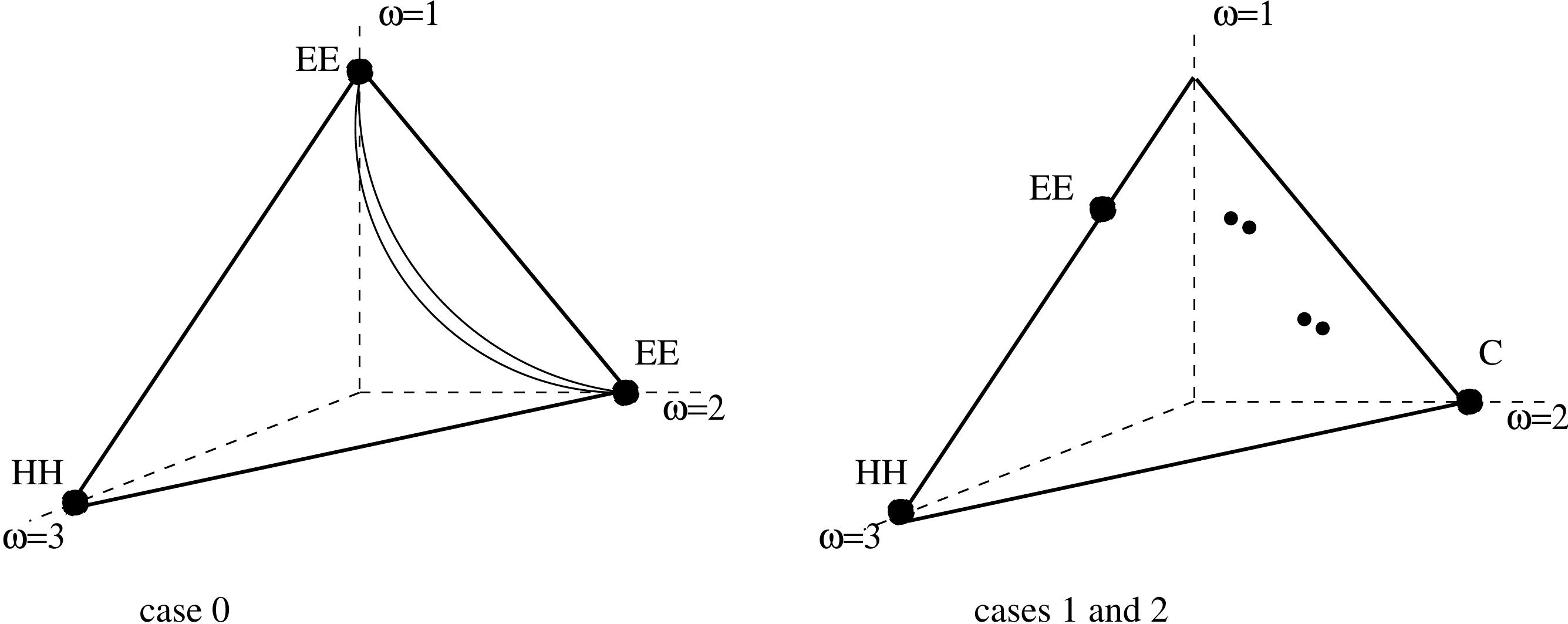}}
\end{center}
\caption{The action simplices of the cases $0, 1$ and $2$; the cases $1$ and
$2$ are typical for the family of Hamiltonians where $0<u<u_1$. The actions
$\tau_i$ (related to $r_i^2$) form a triangle for fixed values of $H_2$
which is an integral of the normal forms. The frequencies have been
normalized to $1, 2, 3$ to indicate the $x_3, x_2, x_1$ normal mode
positions at the vertices. The black dots indicate periodic solutions, the
indicated stability types are HH (hyperbolic-hyperbolic), EE
(elliptic-elliptic) and C (complex). The two (roughly sketched) curves
connecting the $x_2$ and $x_3$ normal modes in the left simplex correspond
with two tori consisting of periodic solutions, respectively with
combination angles $\chi=0$ and $\pi$. The tori break up into $4$ general
position periodic solutions if $u>0$ (cases $1$ and $2$). \label{sim}}
\end{figure}

As the short-periodic solutions have constant actions (or constant radii in
polar coordinates), the integral $H_2$ of the normal form produces for
fixed energy an action-simplex with short-periodic solutions represented by
points; the actions $\tau_i$ and the polar coordinates $r_i$ are related by
the transformations
(\ref{aa}) and (\ref{pc}). One way of displaying the position of
short-periodic solutions and their stability on the $5$-dimensional energy
manifold is the use of this action-simplex with normal modes at the
vertices and solutions in the coordinate planes at the sides. The interior
of the faces may contain short-periodic solutions in general position.
Their stability is indicated by $E$ (elliptic i.e. imaginary eigenvalues),
$H$ (hyperbolic i.e.real eigenvalues)
and $C$ (complex eigenvalues with real parts non-zero). See for instance for
 the action simplices displaying periodic solutions fig. \ref{sim}.

\subsection{Case $0$: the FPU problem with well-balanced masses} \label{bal}

In this case we have the $1:2:3$ resonance with mass values that are as much
as possible similar; we have with $u = 0$ in (\ref{a13}) and (\ref{a24}):
\[ a_1= 0.0357143, a_2= 0.126804, a_3= 0.0357143, a_4= 0.301767. \]
Note that $a_1=a_3$. We checked numerically that the time series
$H_2(p(t), q(t))$ based on the original formulation of system (\ref{nleqs})
and the time series obtained from the transformed Hamiltonian (\ref{Ham1})
produce the same result as it should.


To put system (\ref{nleqs}) in the standard form of quasi-harmonic equations
we have to apply the symplectic transformation $p = K(0)y, q = L(0)x$ in
\eqref{pqyx}. This leads with \eqref{H3dj} and table \ref{tab-dj} to the
transformed Hamiltonian
\begin{equation} \label{Ham1}
H(y, x) = \frac{1}{2} \sum_{i=1}^4 (y_i^2+ \omega_i^2 x_i^2)+
\varepsilon(d_3x_1^2+d_{10}x_2^2+d_6x_1x_2)x_3,
\end{equation}
with
\[ \omega_1^2= \frac{9}{14}, \omega_2^2= \frac{4}{14}, \omega_3^2=
\frac{1}{14}, \omega_4^2=0, d_3= -9 \frac{\sqrt{21}}{490}, d_{10}= 2
\frac{\sqrt{21}}{245}, d_6= -3 \frac{\sqrt{14}}{490}. \]
Rescaling time $t/ \sqrt{14} \rightarrow t$, the equations of motion for the
three dof system become:
\begin{eqnarray} \label{nleqsHam1}
\begin{cases}
\ddot{x}_1 + 9 x_1 & = - \varepsilon 14(2d_3x_1x_3 + d_6x_2x_3), \\
\ddot{x}_2 + 4 x_2 & = - \varepsilon 14(2d_{10}x_2x_3 + d_6 x_1x_3), \\
\ddot{x}_3 + x_3 & = - \varepsilon 14(d_3x_1^2 + d_{10}x_2^2 +d_6 x_1x_2).
\end{cases}
\end{eqnarray}
According to the Weinstein \cite{W} result there exist at least three
families of short-periodic solutions of system
(\ref{nleqsHam1}). Inspection of the equations provides us directly with one
family given by:
\begin{equation} \label{p1Ham1}
x_1(t)=\dot{x}_1(t)=x_2(t)=\dot{x}_2(t)=0,\, \ddot{x}_3 + x_3=0.
\end{equation}
For fixed energy we refer to this periodic solution as the $x_3$ normal
mode; to find such an exact solution explicitly is slightly unusual, the
solution is harmonic. Additional periodic solutions are obtained as
approximations from normal forms as in \cite{HV}. In general, when
normalizing a three dof system, one recovers the three actions (introduced
in (\ref{aa})) and one expects to find the angles in combinations according
to the actual resonances. For the $3:2:1$ resonance these are to first
order after normalization the so-called combination angles
$\phi_1- \phi_2 - \phi_3$ and $2 \phi_3 - \phi_2$. At second order the
combination angle $\phi_1 - 3 \phi_3$ will arise etc., for details see
section 10.2.1 of \cite{SVM}; for instance the term `genuine resonance'
associated with the so-called `annihilators' of $H_2$ can be found in
definition 10.2.2 of \cite{SVM}.

Computing the normal form of system (\ref{nleqsHam1}) to $O(\varepsilon)$
($H_2+ \varepsilon \bar{H}_3$)
as in \cite{HV} or \cite{SVM} and as we shall explicitly show below, only
the $d_6$ term survives in $\bar{H}_3$; this makes the Hamiltonian
(\ref{Ham1}) non-generic. An {\em intermediate normal form} of the
equations of motion becomes:
\begin{eqnarray} \label{Hcase0inter}
\begin{cases}
\ddot{x}_1 + 9 x_1 & = - \varepsilon 14d_6x_2x_3, \\
\ddot{x}_2 + 4 x_2 & = - \varepsilon 14d_6x_1x_3, \\
\ddot{x}_3 + x_3 & = - \varepsilon 14d_6x_1x_2 .
\end{cases}
\end{eqnarray}
As discussed in the Introduction, there is a lot of freedom in choosing
coordinate systems to compute the normal form of the equations of motion.
Near the coordinate planes, in particular to study the stability of the
normal modes, we will use co-moving coordinates. Away from the coordinate
planes (solutions in general position), action-angle variables or polar
coordinates are easier to handle than co-moving coordinates. Some authors
frown upon the use of polar coordinates anyway, because they do not
conserve the canonical character of the normal forms; however, they
preserve the energy and as normal forms they still present a mathematical
precise normal form approximation of the solutions. For general position
orbits we will use in system (\ref{nleqsHam1}) transformations
$x_i, \dot{x}_i \rightarrow r_i, \psi_i$ of the form:
\begin{equation} \label{pc}
x_i = r_i \cos(\omega_i t + \psi_i),\,\, \dot{x}_i = - r_i \omega_i
\sin(\omega_i t + \psi_i).
\end{equation}
The actions $\tau_i$ are related to the $r_i^2$, the angles $\phi_i$ to the
arguments $(\omega_i t + \psi_i)$. Putting $\chi = \psi_1- \psi_2- \psi_3$
and averaging over time $t$, the averaging-normal form equations outside
the coordinate planes become:
\begin{eqnarray} \label{case0nf}
\begin{cases}
\dot{r}_1 & = \varepsilon \frac{7}{6}d_6r_2r_3 \sin \chi, \\
\dot{r}_2 & = - \varepsilon \frac{7}{4}d_6r_1r_3 \sin \chi, \\
\dot{r}_3 & = - \varepsilon \frac{7}{2}d_6r_1r_2 \sin \chi, \\
\dot{\chi} & = \varepsilon \frac{7}{2}d_6 \frac{\cos \chi}{r_1r_2r_3}
\left(\frac{r_2^2r_3^2}{3} - \frac{r_1^2r_3^2}{2} -
\frac{r_1^2r_2^2}{1}\right).
\end{cases}
\end{eqnarray}
 The integral $H_2$ of the normal form equations becomes:
\begin{equation}
9r_1^2+4r_2^2+r_3^2 = 2E_0,
\end{equation}
with $E_0$ a positive (energy) constant. The combination angle
$2 \phi_3 - \phi_2$ is missing; another integral of the normal form
\eqref{case0nf} is:
\begin{equation}\label{int}
2r_2^2 - r_3^2 = C\, (\rm{constant}).
\end{equation}
In the original variables this integral is:
\[ 2x_2^2+ \frac{1}{2}\dot{x}_2^2-x_3^2- \dot{x}_3^2 = \rm{constant}. \]
As we have three independent integrals of the normal form equations
\eqref{case0nf}, the normal form is integrable. Because of the
approximative character of the normal form, this means that chaotic motion
in the original system (\ref{nleqsHam1}) is restricted to $O(\varepsilon)$.

Periodic solutions in general position exist if $\sin \chi =0, t \geq 0$ for
certain values of the $r_i$. From the 4th equation of system
(\ref{case0nf}) we find the requirement:
\[ \frac{r_2^2r_3^2}{3} - \frac{r_1^2r_3^2}{2} - \frac{r_1^2r_2^2}{1} =0 .\]
Eliminating $r_1$ by the $H_2$ integral we find after some rearrangements
the condition
\begin{equation} \label{2fam}
2r_2^2r_3^2+ \frac{4}{3}r_2^4+ \frac{1}{6}r_3^4 =
\frac{1}{3}E_0(2r_2^2+r_3^2),\, 0 < r_2 < \sqrt{\frac{E_0}{2}},\,0 < r_3 <
\sqrt{2E_0}.
\end{equation}
Both for $\chi =0$ and for $\chi = \pi$ we find from condition (\ref{2fam})
tori imbedded in the energy manifold. The two tori consist of periodic
solutions in general position connecting the $x_2$ and $x_3$ normal modes.
Their period is $O(\varepsilon)$ modulated by their position on the tori.
The relation between the presence of a continuous family of periodic
solutions on the energy manifold and the existence of another integral
(\ref{int}) is an example of a more general theory on characteristic
 exponents of periodic solutions developed by Poincar\'e in \cite{PMC}, vol.
 1.\\

\medskip
{\bf Periodic solutions in the coordinate planes}\\
It is clear from the intermediate normal form (\ref{Hcase0inter}) that the
normalized equations of motion will contain all three normal modes. We will
use co-moving coordinates to study the stability:
\begin{eqnarray} \label{co}
\begin{cases}
x_1 & = y_1 \cos 3t + \frac{1}{3}y_2 \sin 3t, \dot{x}_1 = -3y_1 \sin 3t +
y_2 \cos 3t,\\
x_2 & = z_1 \cos 2t + \frac{1}{2}z_2 \sin 2t, \dot{x}_2 = -2z_1 \sin 2t +
z_2 \cos 2t,\\
x_3 & = u_1 \cos t + u_2 \sin t, \dot{x}_3 = -u_1 \sin t + u_2 \cos t.
\end{cases}
\end{eqnarray}
The normalized variables are obtained by averaging over time $t$ and are
satisfying the system:
\begin{eqnarray} \label{nfd6}
\begin{cases}
\dot{y}_1 & = \varepsilon \frac{7}{6} d_6 (z_1u_2 + \frac{1}{2}z_2u_1),\\
\dot{y}_2 & = - \varepsilon \frac{7}{2} d_6 (z_1u_1 - \frac{1}{2}z_2u_2),\\
\dot{z}_1 & = \varepsilon \frac{7}{4} d_6 (-y_1u_2 + \frac{1}{3}y_2u_1),\\
\dot{z}_2 & = - \varepsilon \frac{7}{2} d_6 (y_1u_1 + \frac{1}{3}y_2u_2),\\
\dot{u}_1 & = \varepsilon \frac{7}{2} d_6 (-\frac{1}{2}y_1z_2 +
\frac{1}{3}y_2z_1),\\
\dot{u}_2 & = - \varepsilon \frac{7}{2} d_6 (y_1z_1 + \frac{1}{6}y_2z_2).
\end{cases}
\end{eqnarray}
The generic picture for the existence of short-periodic solutions in the
Hamiltonian $1:2:3$ resonance is given in \cite{HV}. As stated above we
recover three normal modes instead of generically two; this is caused by
the already mentioned degenerate form of Hamiltonian (\ref{Ham1}).

The three normal modes of the normalized system are harmonic functions:
\[ A \cos mt + B \sin mt,\,m=3,2,1,\,A^2+B^2 >0. \]
To study their stability we linearize around the normal modes to obtain
coupled Mathieu equations; we approximate the characteristic exponents by
normalizing these coupled systems. We find:

\begin{figure}[tp]
\begin{center}
\resizebox{14cm}{!}{
\includegraphics{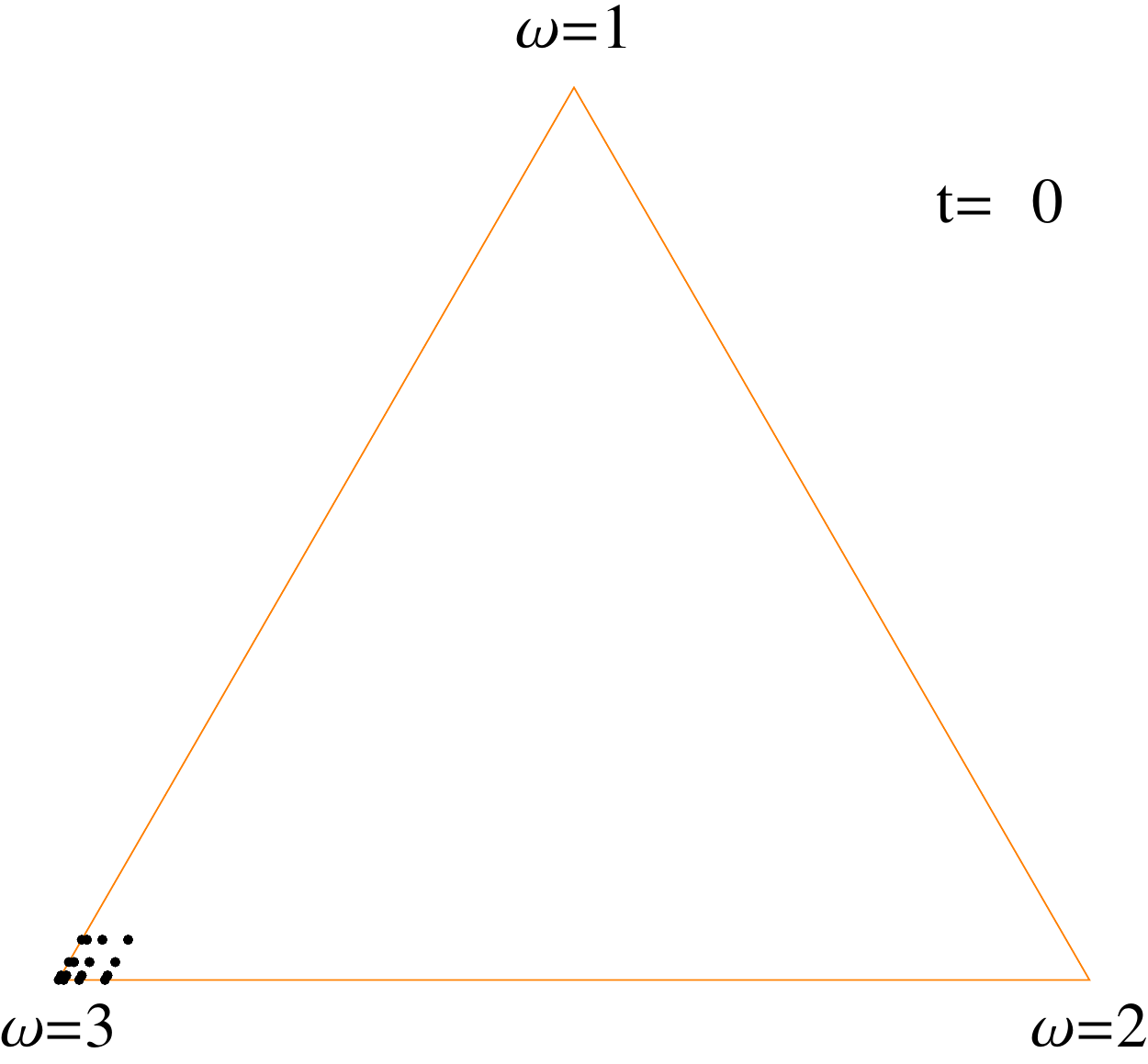}\,\,\includegraphics{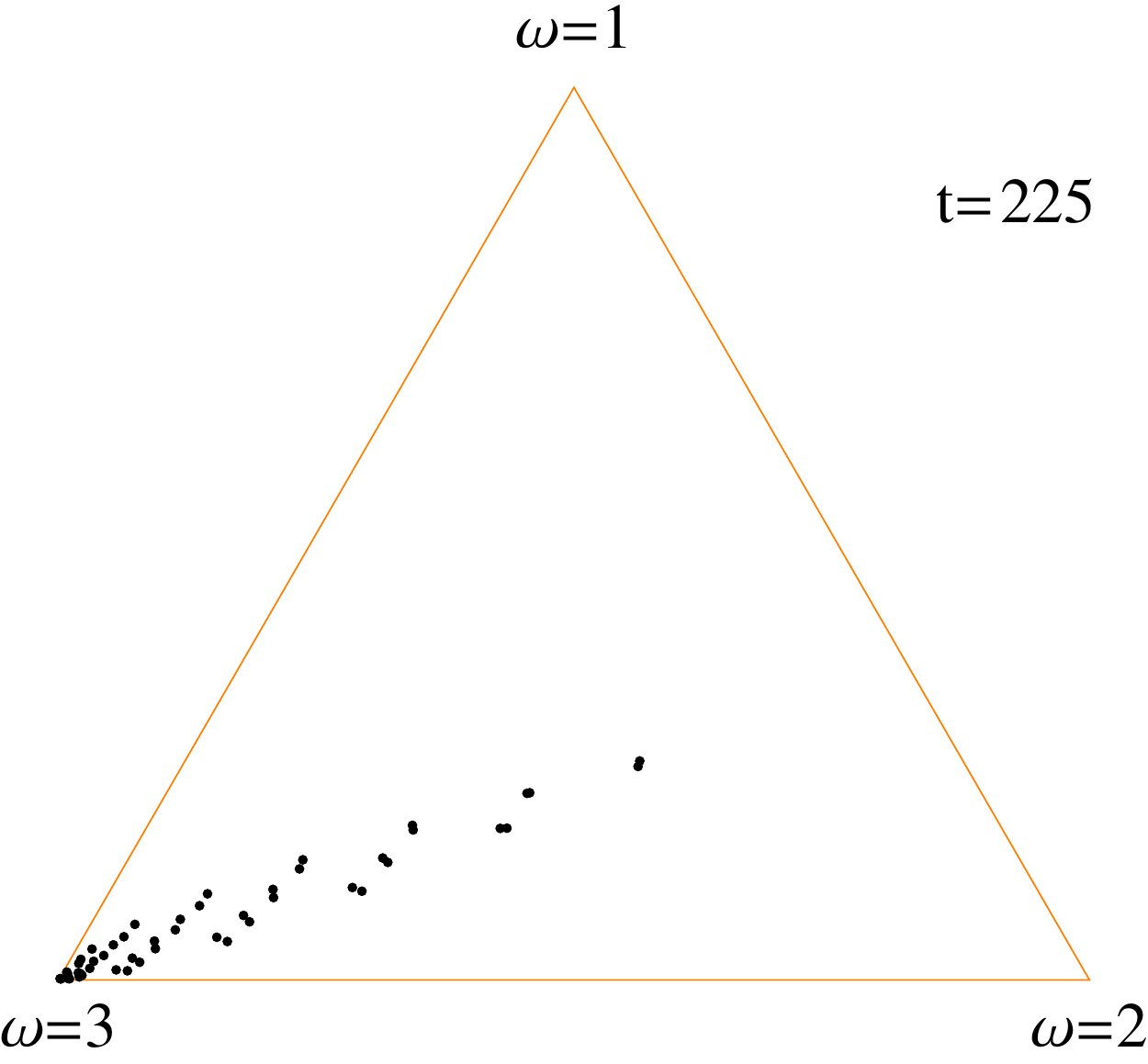}\,\,\includegraphics{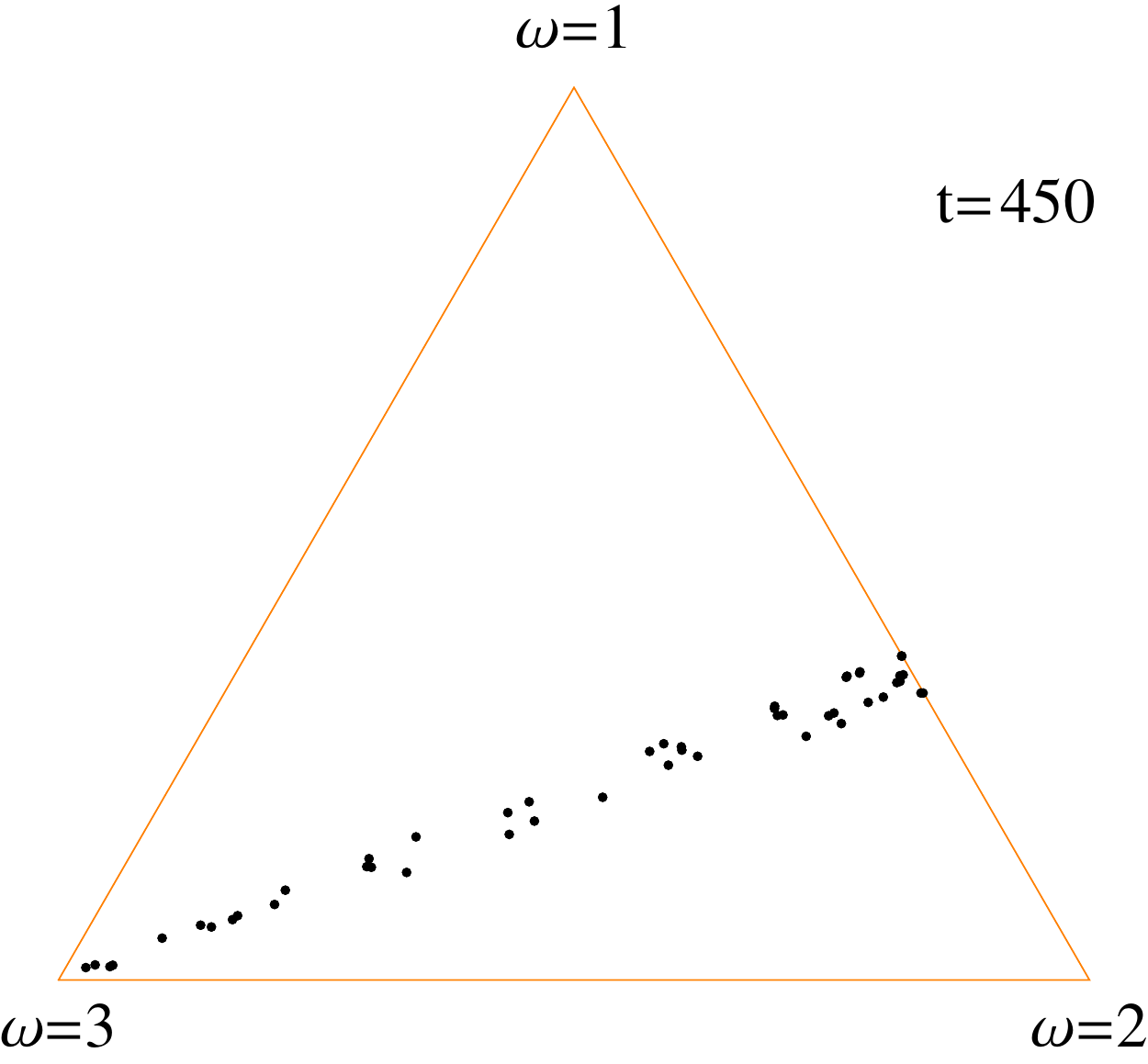}}
\end{center}
\caption{The $\omega =3$ normal mode ($x_1$) exists in the case $0$ and is
unstable (see also fig. \ref{sim}). We consider the time evolution of $98$
initial positions near this normal mode by displaying the actions in the
action-simplex at $t=0,\, 225,\, 450$. The evolution is based on
Hamiltonian \eqref{Ham1}; $ \varepsilon = 0.2$. The unstable manifold is
two-dimensional after which the action points remain near a line in the
action simplex. The inclination is explained by the expression of the third
 integral \eqref{int} of the normal form. \label{simcase0-1}}
\end{figure}

\begin{figure}[tp]
\begin{center}
\resizebox{14cm}{!}{
\includegraphics{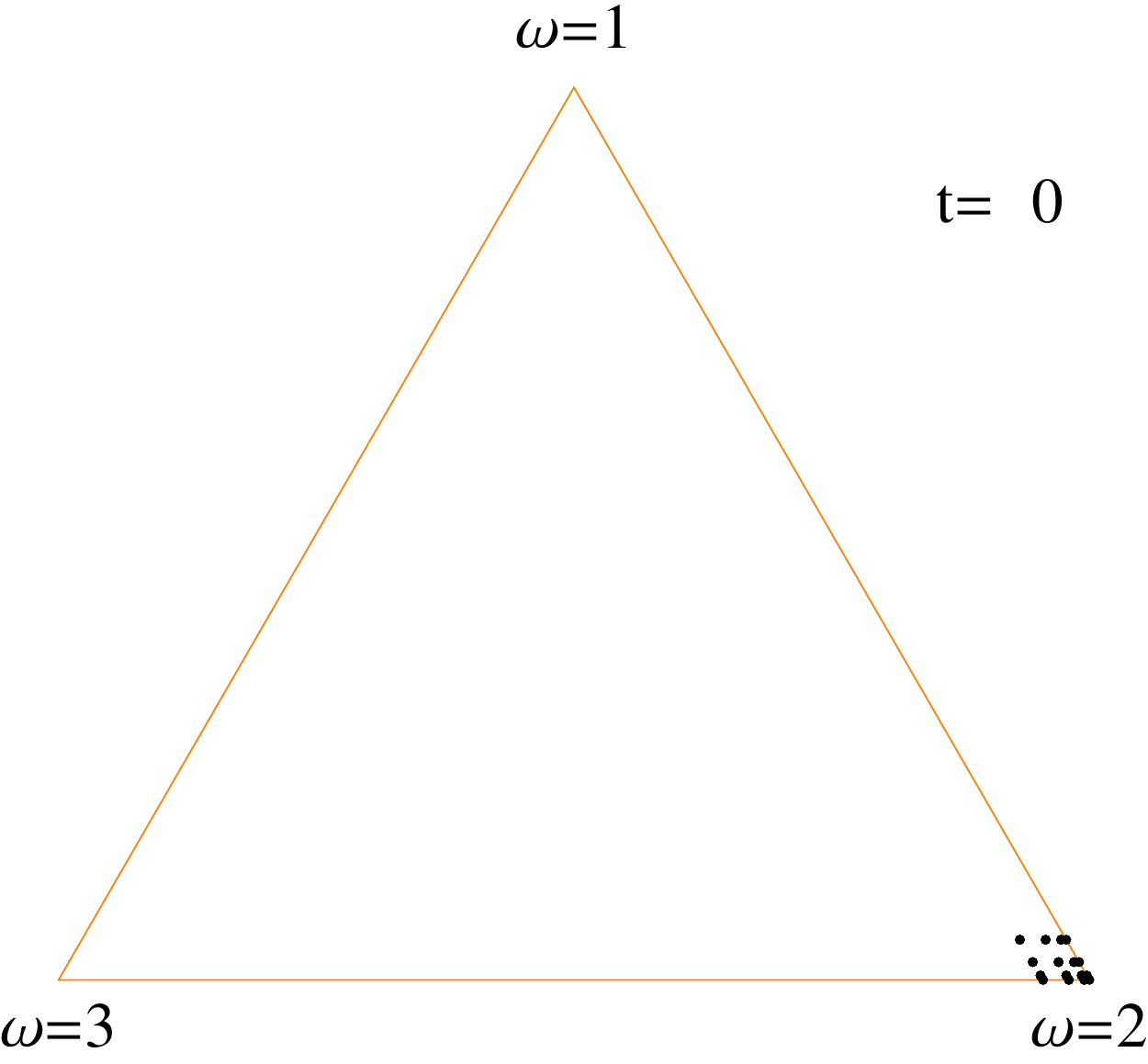}\,\,\includegraphics{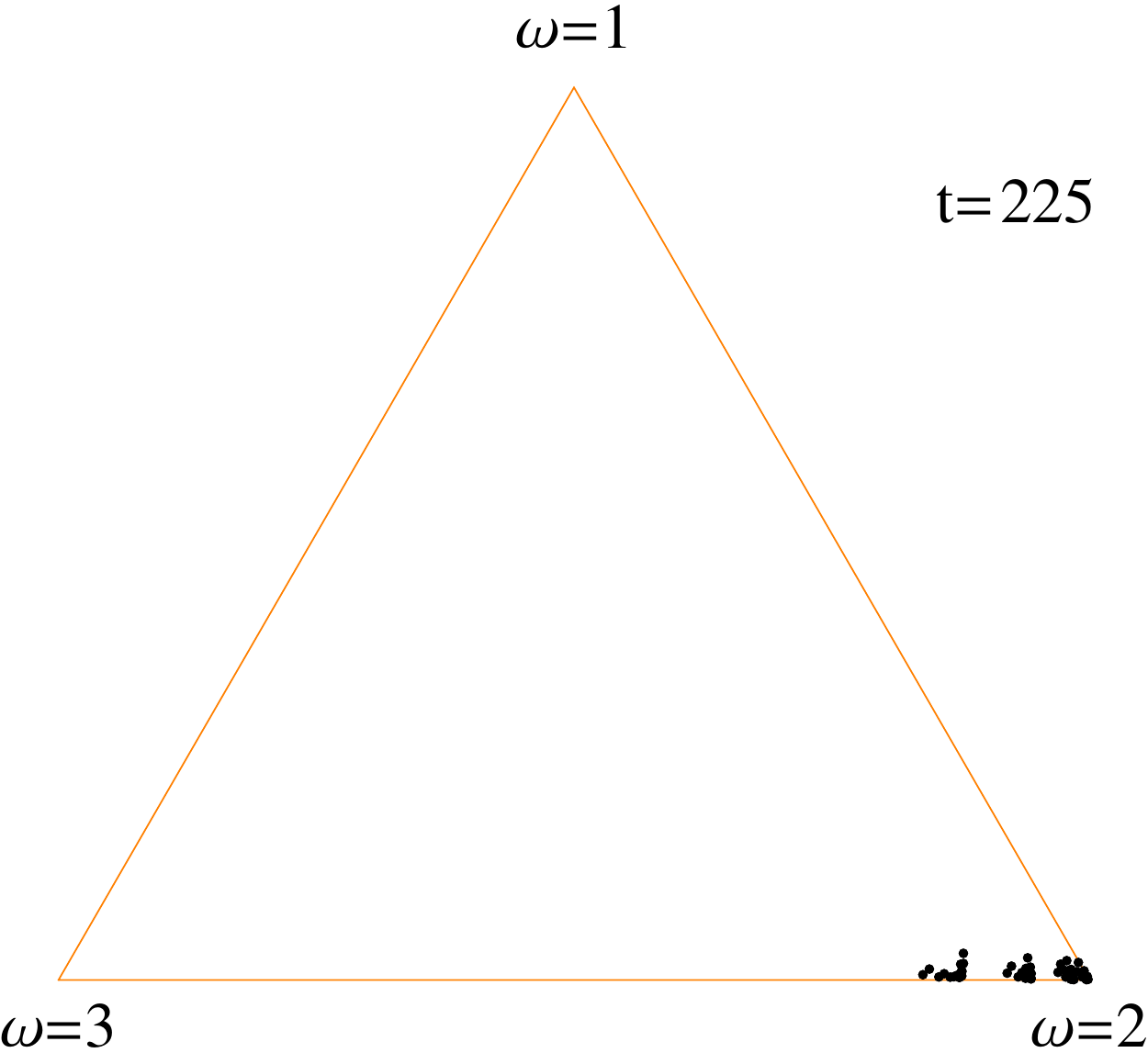}\,\,\includegraphics{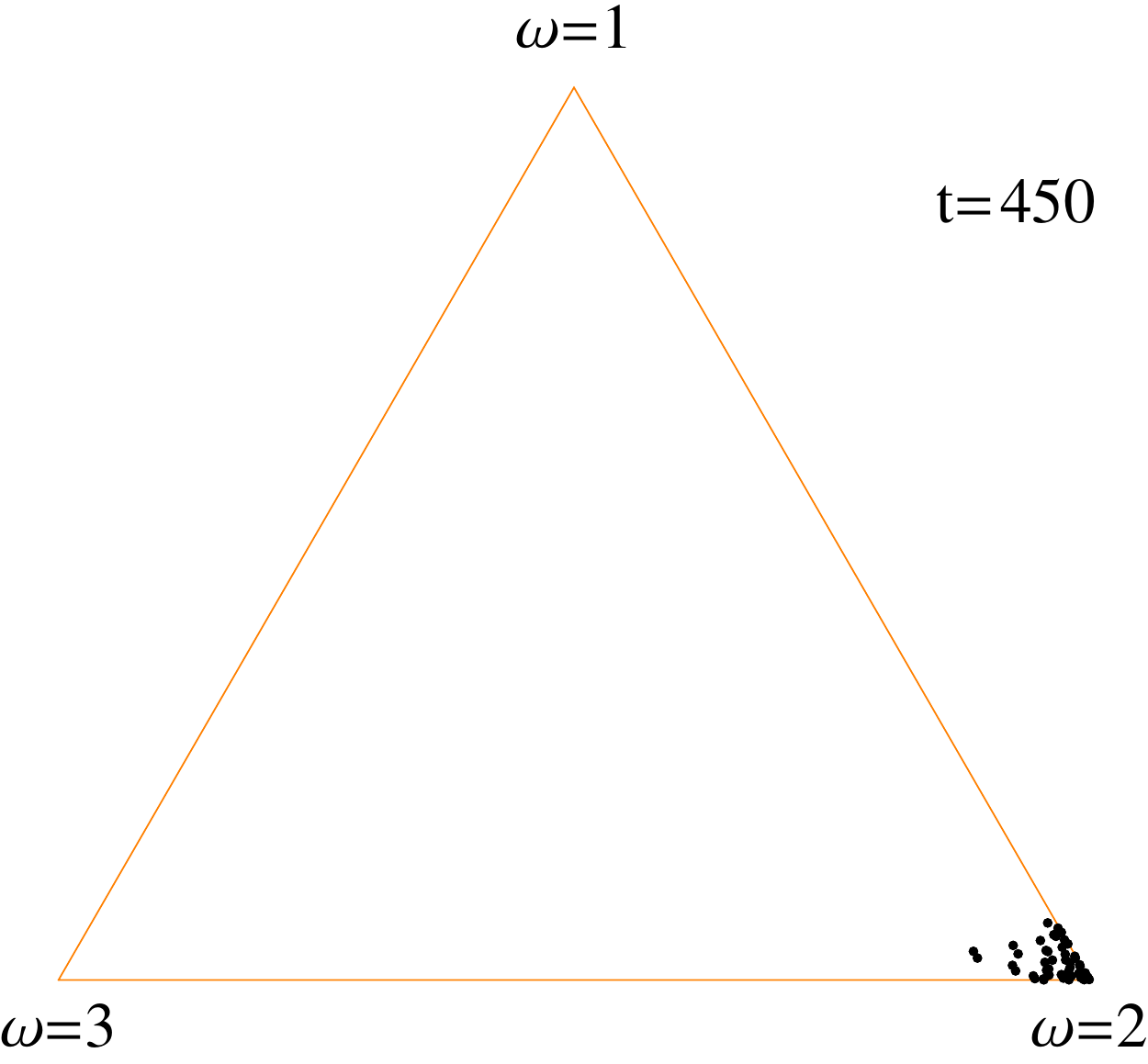}}
\end{center}
\caption{The $\omega =2$ normal mode ($x_2$) is stable in the case $0$. We
consider the time evolution based on Hamiltonian \eqref{Ham1} of $98$
initial positions near this normal mode by displaying the actions in the
action-simplex at $t=0,\, 225,\, 450$; $ \varepsilon = 0.2$.
\label{simcase0-2}}
\end{figure}

\medskip
1. Normal mode $x_1$: put
$x_1= A\cos 3t + B \sin 3t + w_1, x_2=w_2,x_3=w_3$. \\
Transforming in the linearized system by (\ref{co}) and normalization we
find:
\begin{eqnarray*}
\dot{z}_1 & = - \varepsilon \frac{7}{4} d_6(Bu_1 - Au_2),\\
\dot{z}_2 & = \varepsilon \frac{7}{2} d_6(Au_1 + Bu_2),\\
\dot{u}_1 & = - \varepsilon \frac{7}{2} d_6(Bz_1 - \frac{1}{2}Az_2),\\
\dot{u}_2 & = \varepsilon \frac{7}{2} d_6(Az_1 + \frac{1}{2}Bz_2).
\end{eqnarray*}
The eigenvalues of the matrix describing this linear system have
multiplicity $2$ and are multiples of:
\[ \pm \sqrt{A^2+B^2}. \]
In the nomenclature of \cite{SVM} section 10.7.3 this is the unstable case
HH. \\
It is interesting to consider the action-simplex with a number of initial
conditions near the $x_1$ normal mode, see fig. \ref{simcase0-1}. The
unstable manifold of the normal mode is two-dimensional but the solutions,
displayed by dots in the simplex, remain in a narrow strip extending to the
edge where $x_1=0$. This is caused by the third integral \eqref{int} of the
normal form which tells us that the action corresponding with $x_1$ is
proportional to the action of $x_2$.

\medskip
2. Normal mode $x_2$: put
$x_1= w_1, x_2= A\cos 2t + B \sin 2t +w_2,x_3=w_3$. \\
Transforming in the linearized system by (\ref{co}) and normalization by
averaging we find:
\begin{eqnarray*}
\dot{y}_1 & = - \varepsilon \frac{7}{6} d_6(Bu_1 + Au_2),\\
\dot{y}_2 & = \varepsilon \frac{7}{2} d_6(Au_1 - Bu_2),\\
\dot{u}_1 & = \varepsilon \frac{7}{2} d_6(By_1 - \frac{1}{3}Ay_2),\\
\dot{u}_2 & = \varepsilon \frac{7}{2} d_6(Ay_1 + \frac{1}{3}By_2).
\end{eqnarray*}
The eigenvalues have multiplicity $2$ and are multiples of:
\[ \pm i \sqrt{A^2+B^2}. \]
In the nomenclature of \cite{SVM} this is the marginally stable case EE, but
with both positive and negative imaginary eigenvalues coincident. A
numerical calculation confirms the stability in the sense that the
solutions remain near the normal mode during a finite time. \\
When varying $u$, this will produce a Hamiltonian-Hopf bifurcation, see the
next subsection. \\
As the normal mode is marginally stable, it is of interest to display the
behaviour of the actions of solutions starting near this normal mode. In
fig. \ref{simcase0-2} we show that for a limited time interval, the actions
stay nearby.

\medskip
\medskip
3. Normal mode $x_3$: put $x_1= w_1, x_2= w_2,x_3=A\cos t + B \sin t +w_3$.
\\
Transforming in the linearized system by (\ref{co}) and normalization we
find:
\begin{eqnarray*}
\dot{y}_1 & = - \varepsilon \frac{7}{6} d_6(Bz_1 + \frac{1}{2}Az_2),\\
\dot{y}_2 & = \varepsilon \frac{7}{2} d_6(Az_1 - \frac{1}{2}Bz_2),\\
\dot{z}_1 & = \varepsilon \frac{7}{4} d_6(By_1 - \frac{1}{3}Ay_2),\\
\dot{z}_2 & = \varepsilon \frac{7}{2} d_6(Ay_1 + \frac{1}{3}By_2).
\end{eqnarray*}
The eigenvalues have multiplicity $2$ and are multiples of:
\[ \pm i \sqrt{A^2+B^2}. \]
In the nomenclature of \cite{SVM} section 10.7.3, this is the marginally
stable case EE, but again with both positive and negative imaginary
eigenvalues coincident. The numerical behaviour (not shown) looks similar
to fig. \ref{simcase0-2}.

Our choice of well-balanced masses involves the symmetry $a_1=a_3$. In the
sequel we will see that other choices of masses producing $1:2:3$ resonance
give qualitatively different results. It is interesting to compare the
dynamics of case $0$ ($u=d_9=0$) with the dynamics for $u>0$. Such a
comparison will be given in the next subsections.

\subsection{The Hamiltonian-Hopf bifurcation} \label{unbal}
In the preceding subsection we considered a rather symmetric case,
$a_1=a_3$, corresponding with $u=0$, producing an integrable normal form;
see subsection \ref{sect-hamtr} and table \ref{tab-dj}. We will now
consider the cases $0<u<u_1 (= 0.887732\ldots)$; as $u$ increases through
the interval $(0, u_1)$ the masses will differ more and more, producing
generic Hamiltonians. To put system (\ref{nleqs}) in the standard form of
perturbed harmonic equations we have to apply again a symplectic
transformation, i.e. (\ref{pqyx}) from subsection \ref{sect-hamtr} . This
leads to a transformed Hamiltonian (with rescaled frequencies) of the form
$H_2+ \varepsilon H_3$ with
\[ H_2 = \frac{1}{2} (\dot{x}_1^2+ \frac{9}{14}x_1^2 + \dot{x}_2^2+
\frac{4}{14}x_2^2 +\dot{x}_3^2+ \frac{1}{14}x_3^2) \]
and
\begin{eqnarray} \label{Ham2}\begin{cases}
H_3 = & d_1x_1^3+ d_2x_2x_1^2
+d_3x_3x_1^2+d_4x_2^2x_1+d_5x_3^2x_1+d_6x_1x_2x_3+ d_7x_2^3 + d_8x_3^3+\\
& d_9x_2x_3^2 + d_{10}x_2^2x_3,
\end{cases}
\end{eqnarray}
with all coefficients non-zero, see table \ref{tab-dj}. After rescaling time
$t \rightarrow t/ \sqrt{14}$, the equations of motion for the three dof
system can be written as:
\begin{eqnarray} \label{nleqsHam2}
\begin{cases}
\ddot{x}_1 + 9x_1 & = - \varepsilon 14(3d_1x_1^2+2d_2x_1x_2
+2d_3x_1x_3+d_4x_2^2+ d_5x_3^2+ d_6x_2x_3), \\
\ddot{x}_2 + 4 x_2 & = - \varepsilon
14(d_2x_1^2+2d_4x_2x_1+d_6x_1x_3+3d_7x_2^2+ d_9x_3^2+ 2d_{10}x_2x_3), \\
\ddot{x}_3 + x_3 & = - \varepsilon 14(d_3x_1^2+2d_5x_3x_1+d_6x_1x_2
+3d_8x_3^2+ 2d_9x_2x_3 + d_{10}x_2^2).
\end{cases}
\end{eqnarray}
The size of the coefficients of $H_3$ are comparable with the size of $d_6$
or smaller, we will give them explicitly as examples for the cases $1$ and
$2$ in subsection~\ref{exp} with less balanced masses.

In the cubic part of the normalized Hamiltonian we retain of the cubic part
only the terms with $d_6$ and $d_9$; the other terms are, after
normalization, active only at higher order. So, anticipating this, an
intermediate normal form of the equations of motion becomes:
\begin{eqnarray} \label{nleqsHaminter}
\begin{cases}
\ddot{x}_1 + 9 x_1 & = - \varepsilon 14d_6x_2x_3, \\
\ddot{x}_2 + 4 x_2 & = - \varepsilon 14(d_6x_1x_3+d_9x_3^2), \\
\ddot{x}_3 + x_3 & = - \varepsilon 14(d_6x_1x_2 + 2d_9x_2x_3),
\end{cases}
\end{eqnarray}

{\bf The normal form and periodic solutions outside the coordinate planes}\\
Using transformation (\ref{pc}) and putting
$\phi_1- \phi_2 - \phi_3 = \chi_1$, $2 \phi_3 - \phi_2 = \chi_2$, we find
after averaging-normalization:
\begin{eqnarray} \label{casegennf}
\begin{cases}
\dot{r}_1 & = \varepsilon \frac{7}{6}d_6r_2r_3 \sin \chi_1, \\
\dot{r}_2 & = - \varepsilon \frac{7}{4} \left(d_6r_1r_3 \sin \chi_1 +d_9
r_3^2 \sin \chi_2 \right), \\
\dot{r}_3 & = - \varepsilon \frac{7}{2}\left(d_6r_1r_2 \sin \chi_1 -
2d_9r_2r_3 \sin \chi_2 \right), \\
\dot{\chi_1} & = \varepsilon \frac{7}{2}\left[d_6 \frac{\cos
\chi_1}{r_1r_2r_3} \left(\frac{r_2^2r_3^2}{3} - \frac{r_1^2r_3^2}{2} -
\frac{r_1^2r_2^2}{1}\right) - d_9 \frac{\cos \chi_2}{r_2}\left(
\frac{1}{2}r_3^2+2r_2^2\right) \right],\\
\dot{\chi_2} & = \varepsilon \frac{7}{4} \left(d_6 \frac{r_1 \cos
\chi_1}{r_2r_3} (4r_2^2-r_3^2) + d_9 \frac{\cos \chi_2}{r_2} (8r_2^2 -
r_3^2) \right).
\end{cases}
\end{eqnarray}
 The integral $H_2$ of the normal form equations becomes again:
\begin{equation}
9r_1^2+4r_2^2+r_3^2 = 2E_0,
\end{equation}
Periodic solutions in general position with constant amplitude have to
satisfy $\sin \chi_1 = \sin \chi_2 =0$ or $\chi_1= 0, \pi$ and
$\chi_2= 0, \pi$. We have
\[ \cos \chi_1 \cos \chi_2 = \pm 1,\,\, q = \frac{d_9}{d_6} >0. \]
From the last two equations of system (\ref{casegennf}) we have the
conditions:
\begin{equation} \label{c1}
\frac{r_2^2r_3^2}{3} - \frac{r_1^2r_3^2}{2} - \frac{r_1^2r_2^2}{1}= \pm
qr_1r_3 \left( \frac{1}{2}r_3^2+2r_2^2\right),
\end{equation}
\begin{equation} \label{c2}
4r_2^2-r_3^2 = \pm q \frac{r_3}{r_1} (r_3^2-8r_2^2).
\end{equation}
Eliminating $r_1$ from \eqref{c1} using \eqref{c2} we obtain two equations
that are quadratic in $r_2^2$ and $r_3^2$. Eliminating $r_1$ from the $H_2$
integral we find one equation that is quadratic in $r_2^2$ and $r_3^2$.
These expressions have to be handled for the range of $q$ determined by
$u \in (0, u_1)$. Using {\sc Mathematica} and corresponding plots we find
four positive solutions corresponding with four periodic solutions
characterized by two different phases.

We omit the stability analysis, but note that the generic case of the
$1:2:3$ resonance was studied in \cite{HV} that produces four general
position periodic solutions with the stability types $EE$ and $EH$.

\medskip
{\bf Periodic solutions in the coordinate planes}\\
Inspection of the intermediate normal form system (\ref{nleqsHaminter})
shows that the $x_1$ and $x_2$ normal modes exist as solutions of this
system, the $x_3$ normal mode does not. It is shown in \cite{HV} that the
normal mode $x_2$ is unstable. If the instability is of class $C$ (complex
eigenvalues), a Shilnikov-Devaney bifurcation \cite{D76} may take place
resulting in chaotic dynamics originating from a neighborhood of the
complex unstable normal mode. To avoid singularities near the normal modes
we use again the co-moving variables from transformation (\ref{co}). The
normalized variables satisfy the system:
\begin{eqnarray} \label{nfd6d9}
\begin{cases}
\dot{y}_1 & = \varepsilon \frac{7}{6} d_6 (z_1u_2 + \frac{1}{2}z_2u_1),\\
\dot{y}_2 & = - \varepsilon \frac{7}{2} d_6 (z_1u_1 - \frac{1}{2}z_2u_2),\\
\dot{z}_1 & = \varepsilon \frac{7}{2} [ \frac{1}{2}d_6 (-y_1u_2 +
\frac{1}{3}y_2u_1)+d_9u_1u_2],\\
\dot{z}_2 & = - \varepsilon \frac{7}{2} [d_6 (y_1u_1 +
\frac{1}{3}y_2u_2)+d_9(u_1^2-u_2^2)],\\
\dot{u}_1 & = \varepsilon \frac{7}{2}[ d_6 (-\frac{1}{2}y_1z_2 +
\frac{1}{3}y_2z_1)+ d_9(-2z_1u_2+z_2u_1)],\\
\dot{u}_2 & = - \varepsilon \frac{7}{2}[ d_6 (y_1z_1 +
\frac{1}{6}y_2z_2)+d_9(2z_1u_1+z_2u_2)].
\end{cases}
\end{eqnarray}
 We find three families of short-periodic solutions; the constants $A, B$
 are real, $A^2+B^2>0$.
\begin{enumerate}
\item $x_1(t)=A \cos 3t + B \sin 3t, x_2=x_3=0$.
\item $x_2(t)= A \cos 2t + B \sin 2t, x_1=x_3=0.$
\item If $x_2(t)=0, d_6 \neq 0$:
\begin{eqnarray}\label{spec}
\begin{cases}
x_1(t) & = \frac{d_9}{d_6}\left(\frac{A}{A^2+B^2}(3B^2-A^2) \cos 3t -
\frac{B}{A^2+B^2}(3A^2-B^2) \sin 3t\right), \\
x_3(t) & = A \cos t + B \sin t.
\end{cases}
\end{eqnarray}
If $d_9$ differs from zero, this family of periodic solutions moves along
the $x_2=0$ edge of the simplex in fig. \ref{sim} starting from the $x_3$
normal mode that exists if $d_9=0$.
\end{enumerate}
To evaluate the stability of the periodic solutions we will linearize system
(\ref{nleqsHaminter}) near these solutions; this produces coupled Mathieu
equations which we will analyze by normalization.\\

\medskip
{\bf The $x_2$ normal mode}\\
 Put:
\[ x_1=w_1, x_2 = A \cos 2t + B \sin 2t +w_2, x_3=w_3, \]
with real constants $A, B, A^2+B^2 >0$ and corresponding expressions for the
derivatives. We find after linearization
\begin{eqnarray} \label{nleqslinx2}
\begin{cases}
\ddot{w}_1 + 9 w_1 & = - \varepsilon 14 d_6 (A \cos 2t + B \sin 2t)w_3, \\
\ddot{w}_2 + 4 w_2 & = 0, \\
\ddot{w}_3 + w_3 & = - \varepsilon 14 [d_6w_1(A \cos 2t + B \sin 2t) + 2d_9
(A \cos 2t + B \sin 2t)w_3],
\end{cases}
\end{eqnarray}
We study the stability of this system by normalization to find the
eigenvalues of the matrix (omitting the factor $7 \varepsilon/2$)
\[ \left( \begin{array}{cccc}
0 & 0 & \frac{d_6}{3}B & \frac{d_6}{3}A \\
0 & 0 & -d_6A & d_6B \\
-d_6B & \frac{d_6}{3}A & 2d_9B & - 2d_9A \\
-d_6A & -\frac{d_6}{3}B & - 2d_9A & -2 d_9B
\end{array} \right) \]
produce first order approximations of the characteristic exponents of system
(\ref{nleqslinx2}). For the eigenvalues we find apart from the factor
$7 \varepsilon/2$:
\[ \lambda^2 = -(A^2+B^2)\left[(\frac{1}{3}d_6^2-2d_9^2) \pm 2d_9
\sqrt{d_9^2- \frac{1}{3}d_6^2}\right]. \]
A sufficient condition for the complex case $C$ to arise is
\begin{equation} \label{condC}
d_6^2 > 6d_9^2.
\end{equation}
This condition corresponds with the condition in table 1 of \cite{HV}.
Condition (\ref{condC}) is satisfied for $0<u <u_1$ so that the complex
case $C$ arises for $u>0$. \\

\begin{figure}[htp]
\begin{center}
\resizebox{10cm}{!}{
\includegraphics{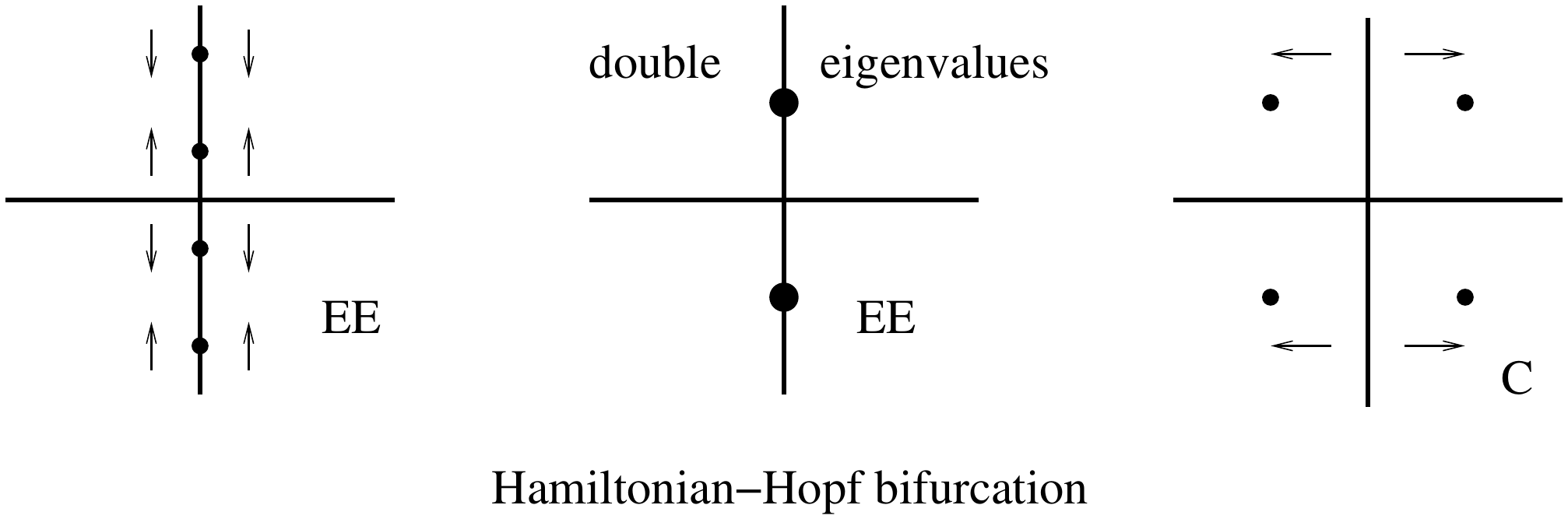}}
\end{center}
\caption{The Hamiltonian-Hopf bifurcation of a periodic solution in a three
dof system as takes place for the $x_2$ normal mode in $1:2:3$ resonance of
\cite{HV}. \label{HnaarC}}
\end{figure}

Another view of the eigenvalues is obtained by realizing that in subsection
\ref{bal} we had $u=0$ resulting in $d_6 \neq 0, d_9=0$; $u=0$ gives for
the $x_2$ normal mode purely imaginary eigenvalues with multiplicity two.
As $u$ increases ($d_9 \neq 0$), the eigenvalues move from the imaginary
 axis into the complex domain. This is part of the Hamiltonian-Hopf
 bifurcation, see fig. \ref{HnaarC}.

\begin{figure}[htp]
\begin{center}
\resizebox{14cm}{!}{
\includegraphics{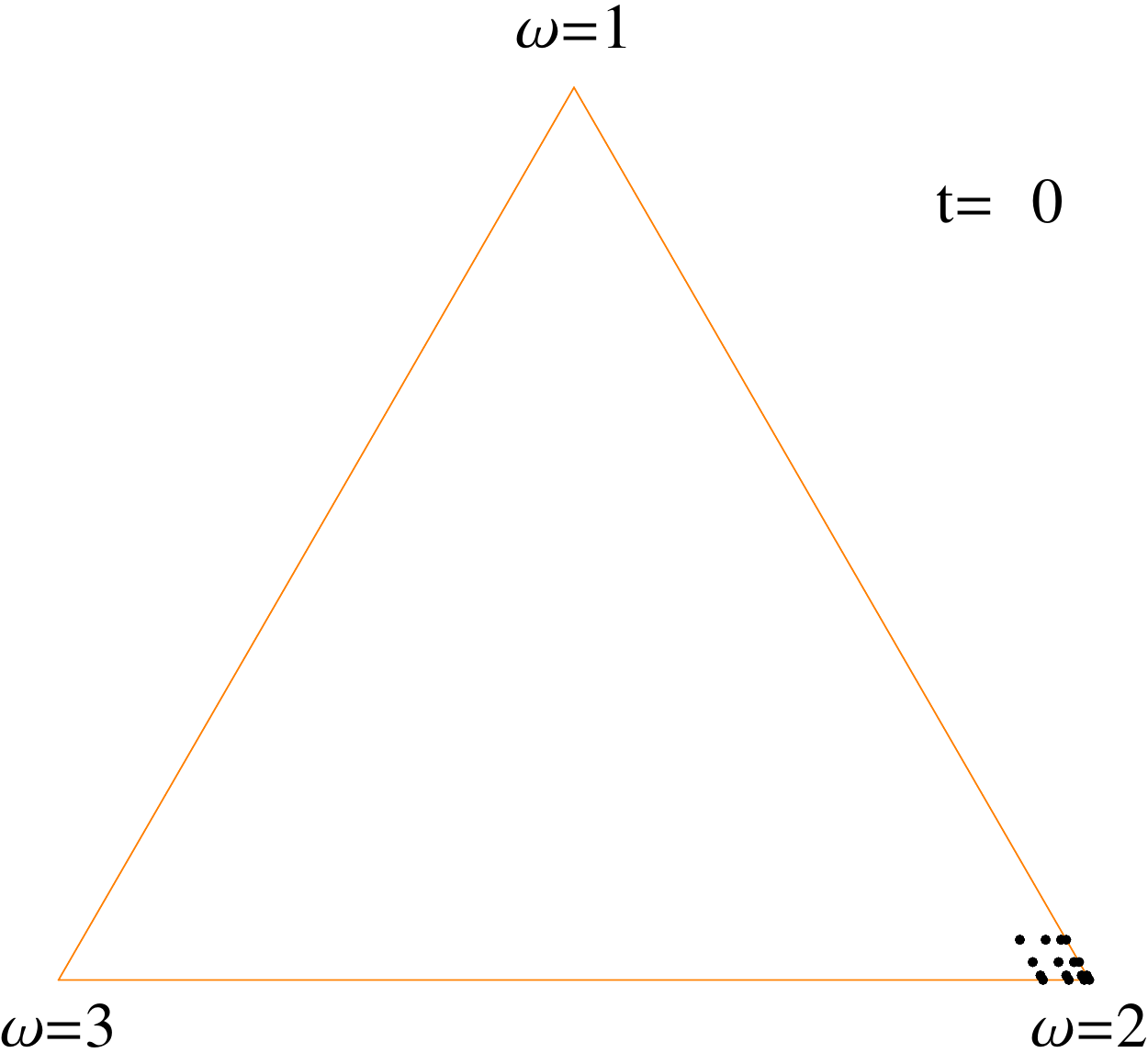}\,\,\includegraphics{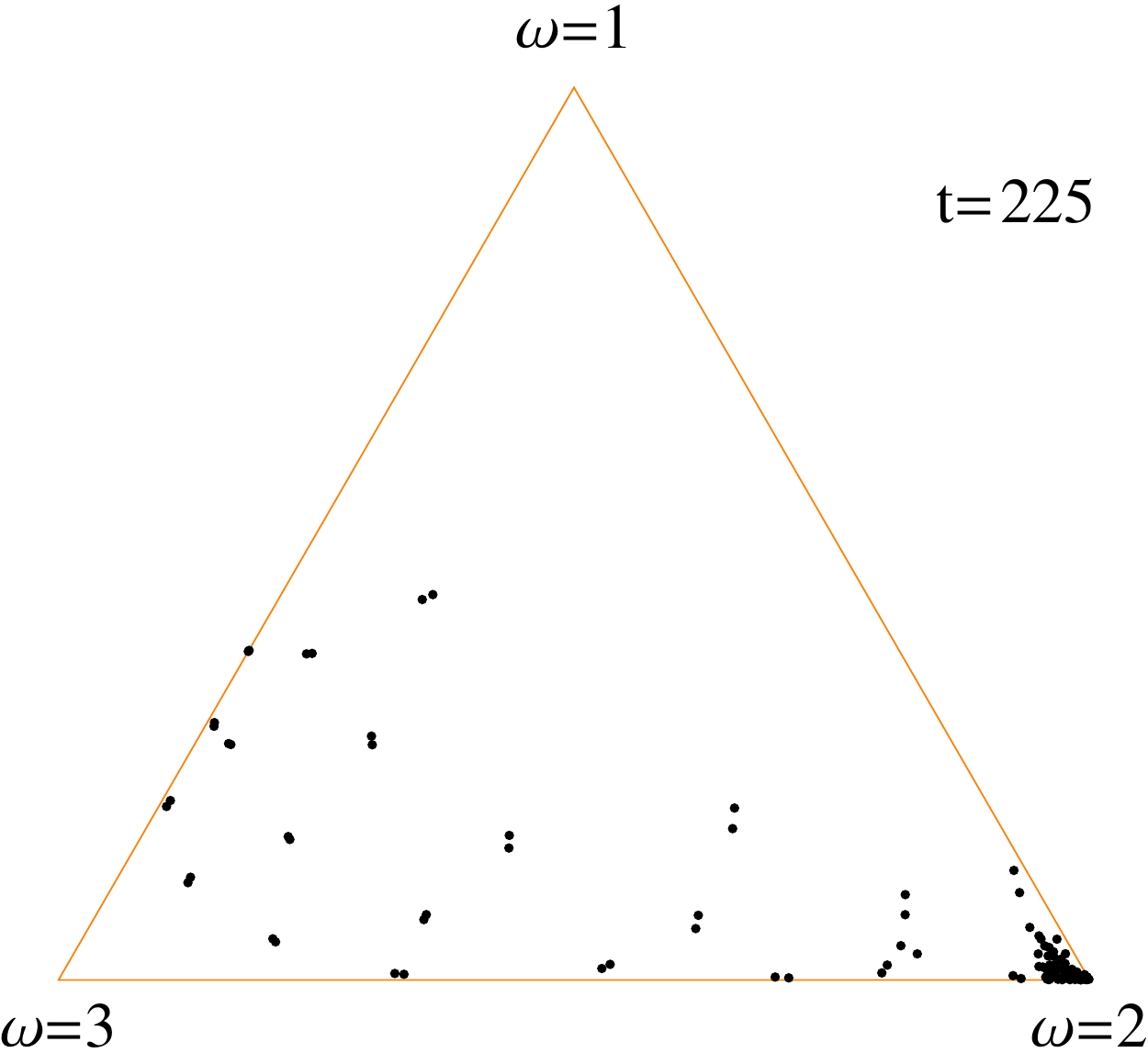}\,\,\includegraphics{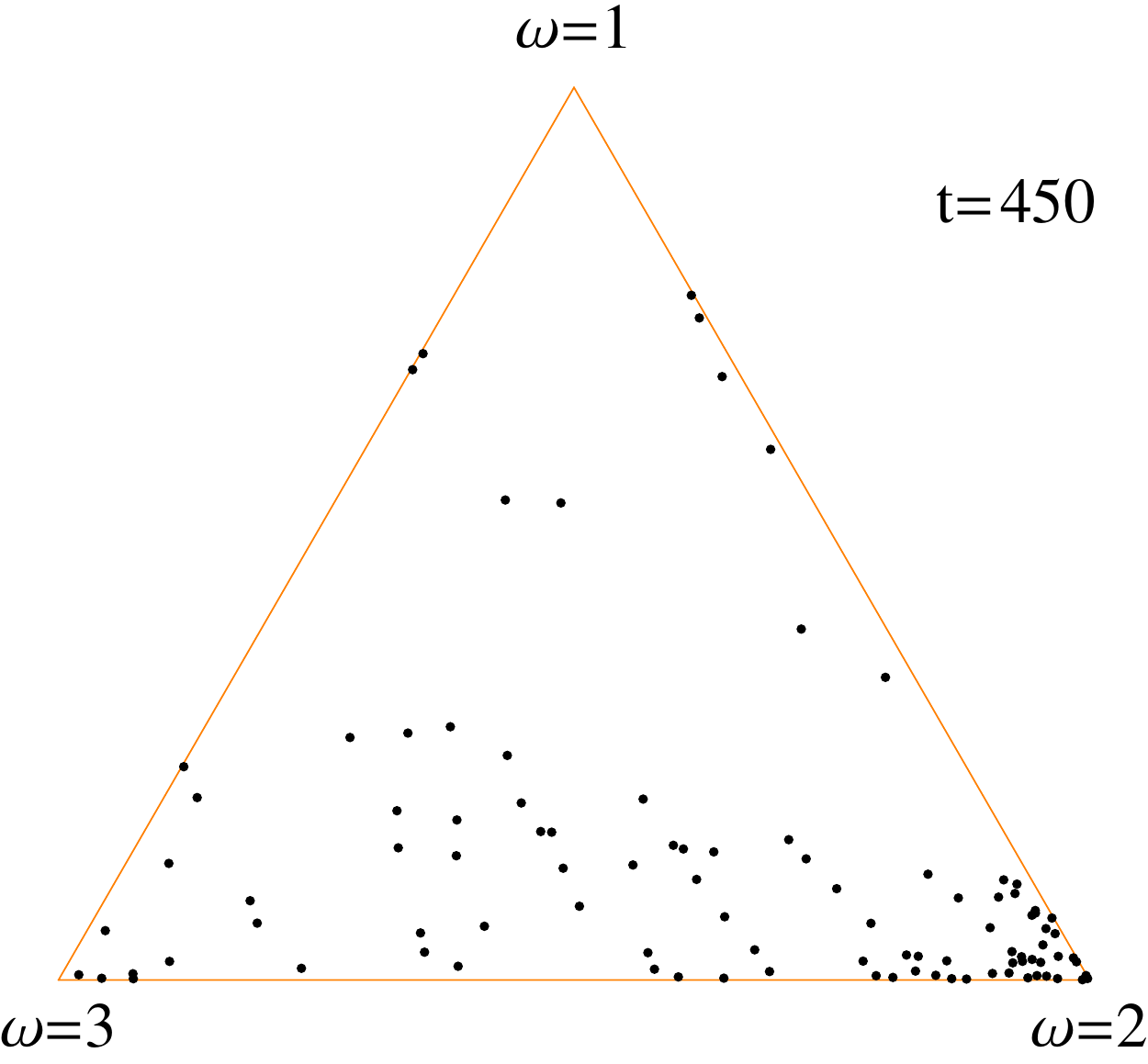}}
\end{center}
\caption{The $\omega =2$ normal mode ($x_2$) exists in the case $2$ and is
complex unstable (see also fig. \ref{sim}). We consider the time evolution
of $98$ initial positions near this normal mode by displaying the actions
in the action-simplex at $t=0,\, 225,\, 450$; $\varepsilon = 0.2$.
\label{simcase2-2}}
\end{figure}

For case $2$ (see subsection \ref{exp}) we show in the action-simplex of
fig. \ref{simcase2-2} the behaviour of solutions starting near this complex
unstable normal mode.

\begin{figure}[htp]
\begin{center}
\resizebox{14cm}{!}{
\includegraphics{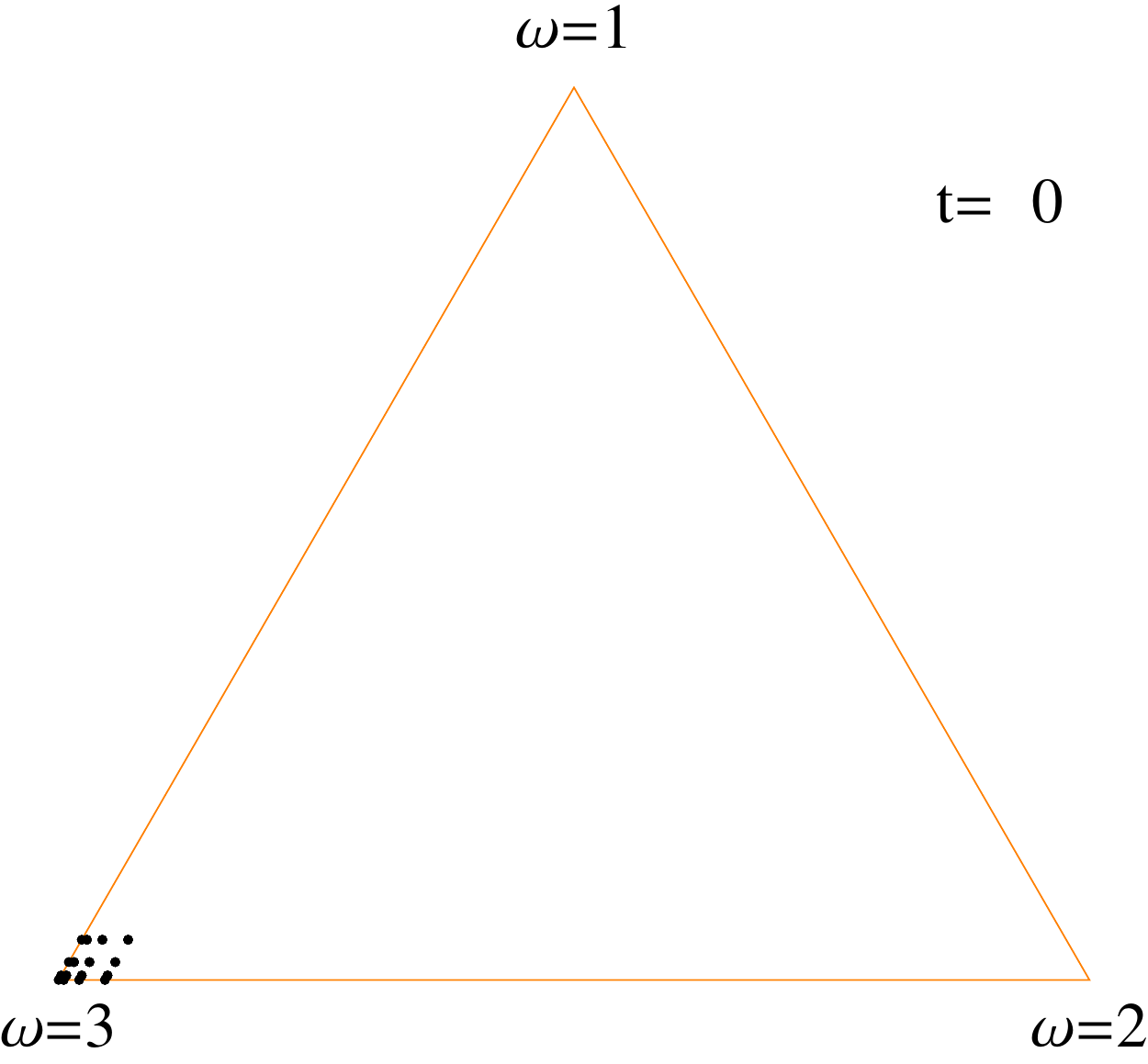}\,\,\includegraphics{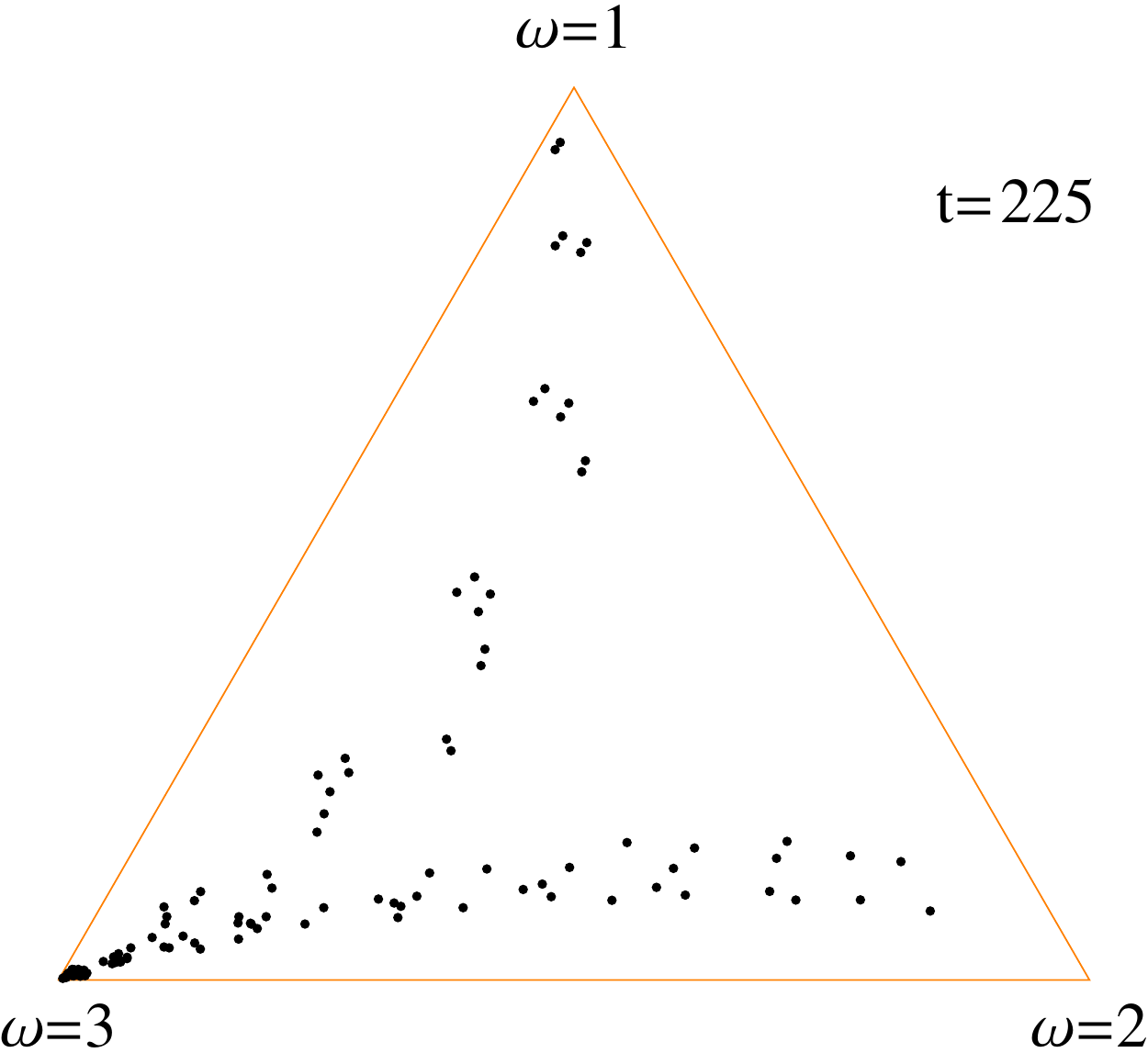}\,\,\includegraphics{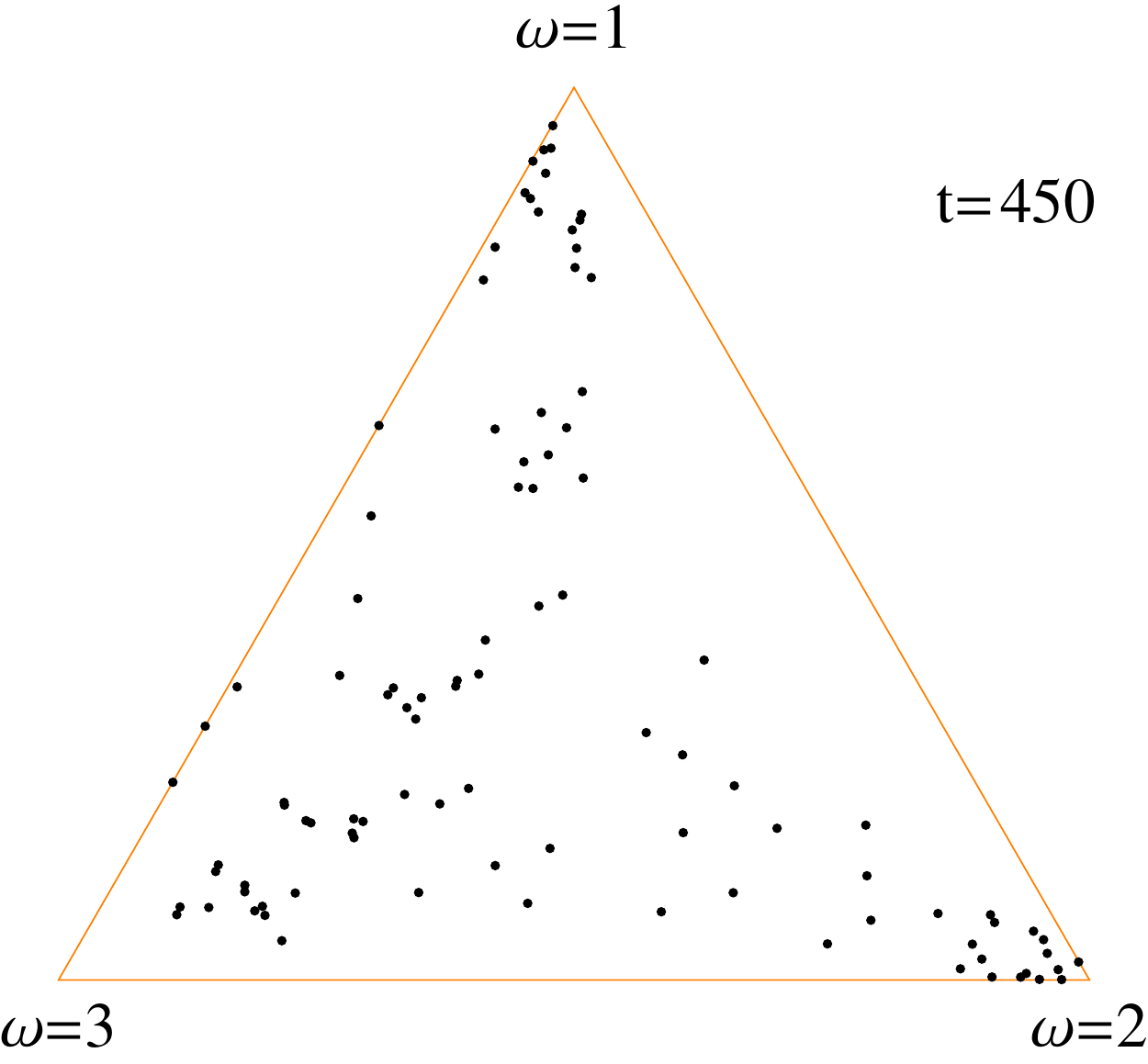}}
\end{center}
\caption{The $\omega =3$ normal mode ($x_1$) exists in the case $2$ and is
unstable (see also fig. \ref{sim}). We consider the time evolution of $98$
initial positions near this normal mode by displaying the actions in the
action-simplex at $t=0,\, 225,\, 450$, $\varepsilon = 0.2$. The behaviour
is different from the case $0$, see fig. \ref {simcase0-1}, as in this case
the normal form is not integrable. \label{simcase2-1}}
\end{figure}

\medskip
{\bf The $x_1$ normal mode}\\
 For $d_9=0$ we have found in the preceding subsection the case HH. This is
 a generic case of eigenvalues, so for $d_9$ small enough the nature of the
instability will not change but the dynamics is very different as the
normal form is not integrable. \\
For case $2$ (see subsection \ref{exp}) we show in the action-simplex of
fig. \ref{simcase2-1} the behaviour of solutions starting near this
unstable normal mode.

\medskip
{\bf The periodic solution for $x_2(t)=0$ .}\\
For the periodic solution \eqref{spec} we put:
\[ x_1= C \cos 3t + D \sin 3t,\,\, x_3 = A \cos t + B \sin t. \]
Transforming
\[ x_1= C \cos 3t + D \sin 3t + w_1, x_2=w_2, x_3= A \cos t + B \sin t +
w_3, \]
and substitution into system \eqref{nleqsHaminter}, we find after
linearization:
\begin{eqnarray} \label{nleqslin2}
\begin{cases}
\ddot{w}_1 + 9 w_1 & = - \varepsilon 14d_6 (A \cos t + B \sin t)w_2, \\
\ddot{w}_2 + 4 w_2 & = - \varepsilon 14[d_6(C \cos 3t +D \sin 3t)w_1 +\\
& d_6 ( A \cos t + B \sin t)w_1 +2d_9( A \cos t + B \sin t)w_3], \\
\ddot{w}_3 + w_3 & = - \varepsilon 14[d_6(C \cos 3t +D \sin 3t)w_2 +2d_9 (A
\cos t + B \sin t)w_2].
\end{cases}
\end{eqnarray}
To investigate stability we normalize near the periodic solution; apart from
a factor $7 \varepsilon /2$, this produces the matrix:
\[ \left( \begin{array}{cccccc}
0 & 0 & \frac{d_6}{3}B & \frac{d_6}{6}A & 0 & 0 \\
0 & 0 & -d_6A & \frac{d_6}{2}B & 0 & 0 \\
-\frac{d_6}{2}B & \frac{d_6}{6}A & 0 & 0 & d_6 \frac{D}{2} +d_9B & - d_6
\frac{C}{2}+d_9A \\
-d_6A & - \frac{d_6}{3}B & 0 & 0 & -d_6C- 2d_9A & -d_6D+2 d_9B \\
0 & 0 & d_6D-2d_9B & - d_6 \frac{C}{2}+d_9A & 0 & 0 \\
0 & 0 & -d_6C - 2d_9A & -d_6 \frac{D}{2}-d_9B & 0 & 0
\end{array} \right). \]
Using the values of $C$ and $D$ given in \eqref{spec}, we find purely
imaginary eigenvalues with multiplicity two. The results have been
summarized in fig. \ref{sim}.

\subsection{Experiments for two cases with $u>0$} \label{exp}
We consider a few experiments for two cases that are typical for the
dynamics when $u>0$.

\subsection*{Case $1$ with less-balanced masses} 
\begin{figure}[htp]
\begin{center}
\includegraphics[ viewport=0 90 580 750 , clip, height=7cm,
width=7.2cm]{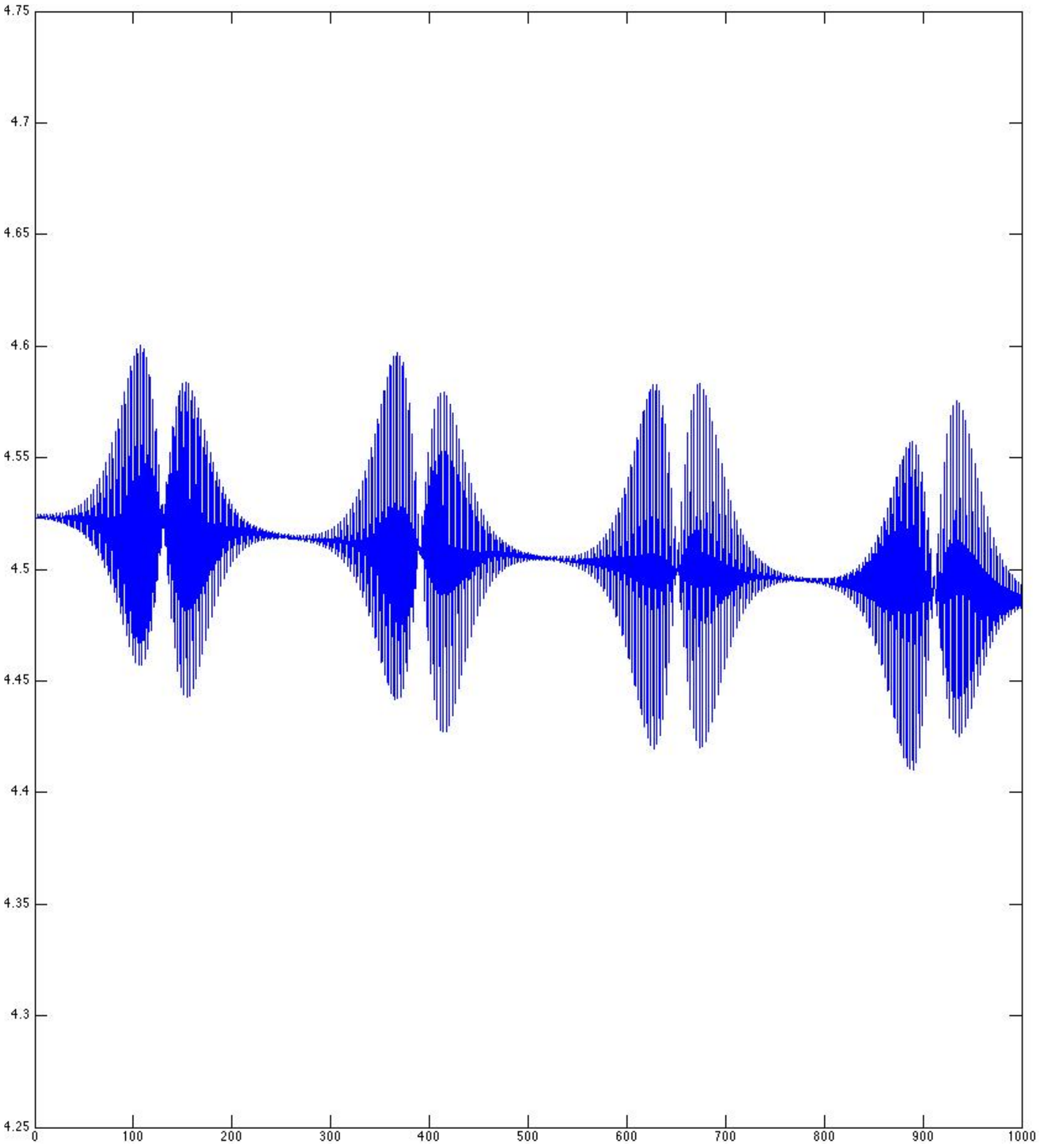} \quad
 \includegraphics[viewport=45 195 585 615 ,clip, 
height=7cm,width=7cm]{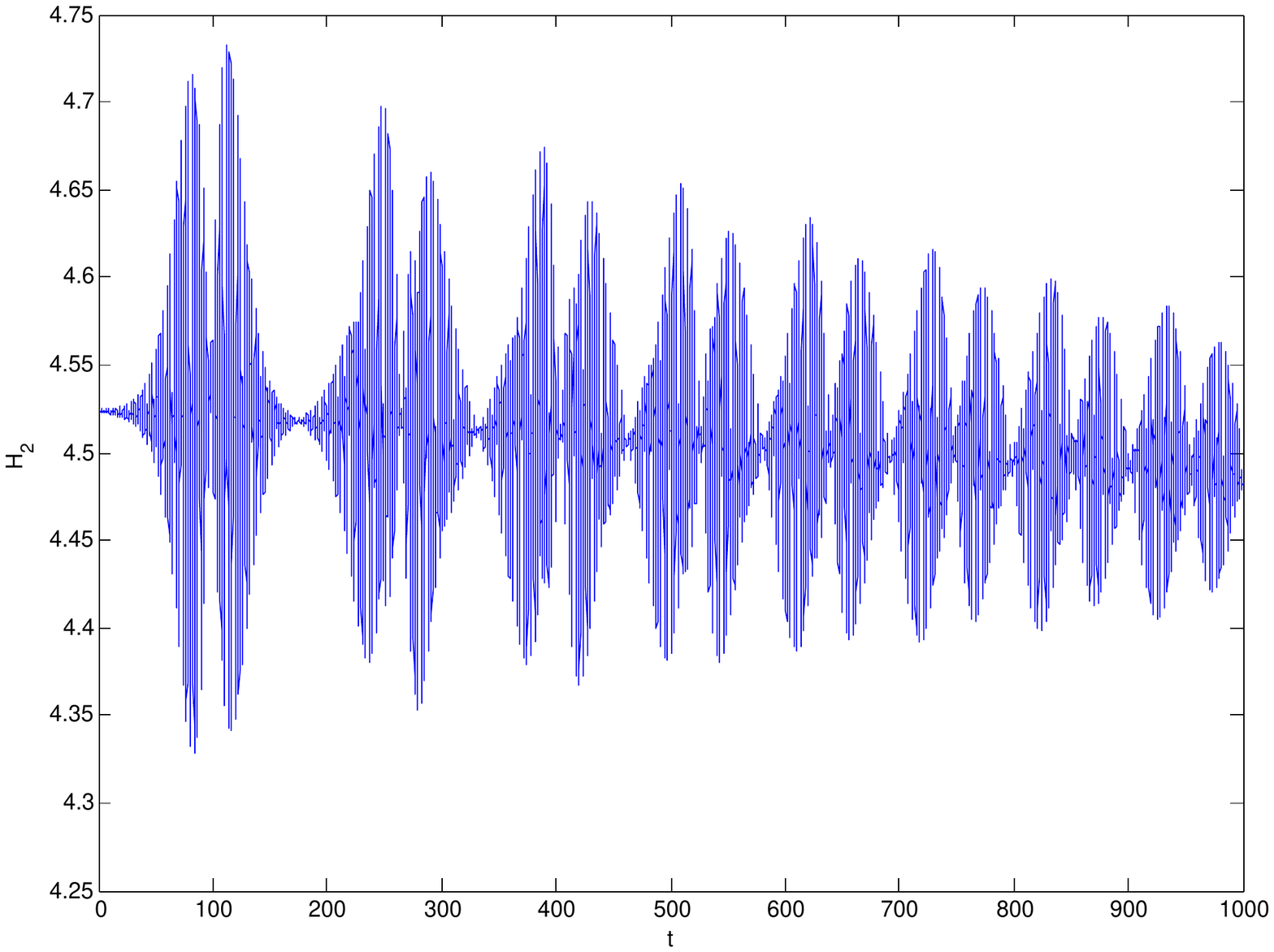}
\end{center}
\caption{Left the $H_2(t)$ time series based on system \eqref{Hcase0inter},
case $0$, with initial values
$x_1=1, x_2=0.1, x_3= 0.1, \dot{x}_i=0, i=1, 2, 3$;
$\varepsilon = 0.5, H_2(0) \approx 4.52$. Right the $H_2(t)$ time series
for case $1$ based on system \eqref{nleqsHaminter} with the same initial
conditions. Horizontal scales: time in $[0,1000]$, vertical scales: energy
in $[4.25,4.75]$. \label{H2case0-1}}
\end{figure}

\begin{figure}[htp]
\begin{center}
\includegraphics[ viewport=0 90 580 748 , clip, 
height=7cm, width=7.2cm]{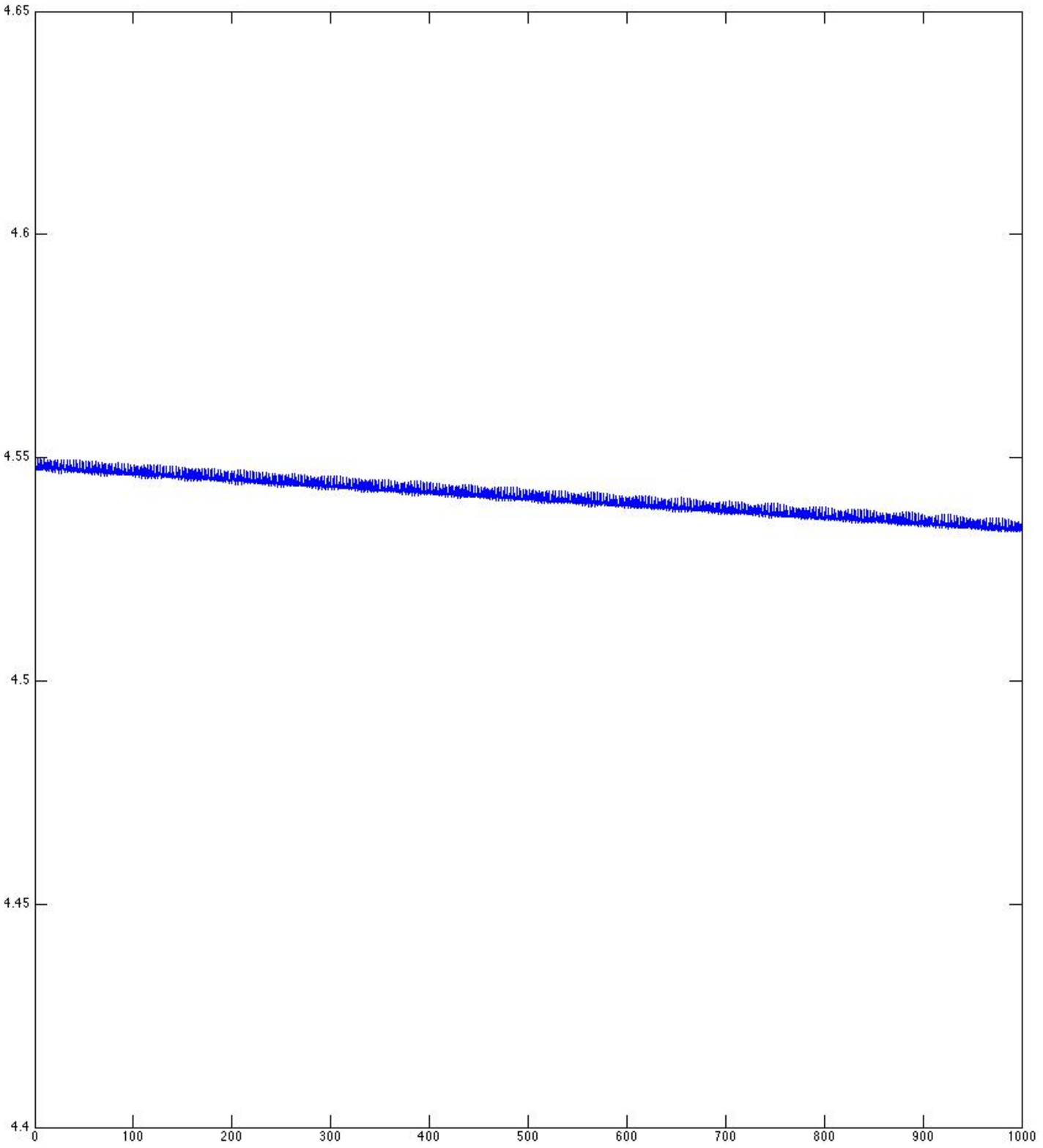} \quad
 \includegraphics[viewport=45 195 585 615 ,clip,
height=7cm,width=7cm]{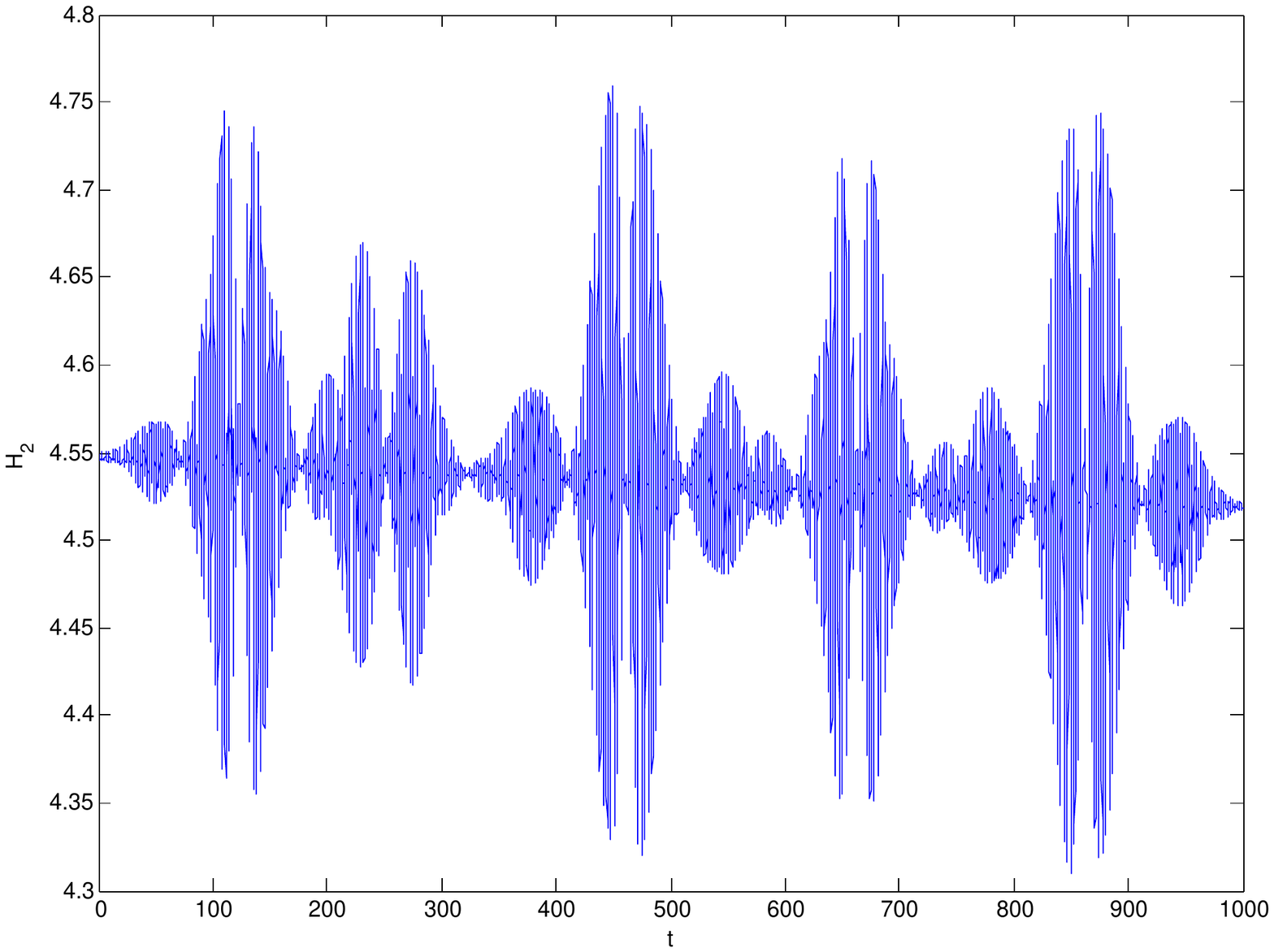}
\end{center}
\caption{Left the $H_2(t)$ time series based on system \eqref{Hcase0inter},
case $0$, with initial values
$x_1=0.1, x_2=1.5, x_3= 0.1, \dot{x}_i=0, i=1, 2, 3$;
$\varepsilon = 0.5, H_2(0) \approx 4.52$. Right the $H_2(t)$ time series
for case $1$ based on system \eqref{nleqsHaminter} with the same initial
conditions i.e. near the complex unstable $x_2$ normal mode. Horizontal
scales: time in $[0,1000]$, vertical scales: energy in $[4.4,4.65]$ (left)
and in $[4.3,4.8]$ (right). \label{H2case0-1x2}}
\end{figure}

We choose for
$u = 0.534105$ from eqs. (\ref{a13}) and (\ref{a24}):
 \[ a_1= 0.00510292, a_2= 0.117265, a_3= 0.0854008, a_4= 0.292231 \]
In this case we have $m_1 > m_3 > m_2> m_4$. With these mass ($a_i$) values
the symplectic transformation of subsection \ref{sect-hamtr} to system
(\ref{nleqsHam2}) produces the expression:
\begin{eqnarray*}
H_3 & = 0.0281999 x_1^3 -0.0258437 x_1^2x_2 -0.0777574x_1^2x_3-
0.0275058x_1x_2^2- 0.00252349x_1x_3^2 \\
& -0.0306229x_1x_2x_3+ 0.0157538x_2^3 + 0.000502655x_3^3- 0.0089438x_2x_3^2
+ 0.028527x_2^2x_3.
\end{eqnarray*}
We have the case:
\[ d_6= -0.0306229,\, d_9= -0.0089438. \]
so that the $x_2$ normal mode is complex unstable; see fig. \ref{sim}.
$H_2(t)$ time series are shown in figs. \ref{H2case0-1} and
\ref{H2case0-1x2}.

Note that $d_9$ is still fairly small with the implication that the
expansion of the flow near the $x_2$ normal mode will not be very
explosive. This may reduce the amount of chaos present in the system. We
will compare with case $0$ and give a few more details for different
initial conditions based on integration of system \eqref{Hcase0inter} and
system \eqref{nleqsHaminter}. We established that in all cases the $x_1$
normal mode is unstable ($HH$), see also fig. \ref{sim}. Starting near the
$x_1$ normal mode in case $0$, the solutions move away, guided by the
two-dimensional unstable manifold of the normal mode; the integrability of
the normal form produces a fairly regular $H_2(t)$, see fig.
\ref{H2case0-1}. Also in this figure we display $H_2(t)$ for case $1$ with
the same initial conditions; its behaviour is influenced by the chaotic
character of the normal form. On this interval of time $[0, 1000]$, energy
is clearly pumped into $H_3$ but the recurrence of the Hamiltonian system
will return this on a much longer timescale.

The chaos in case $1$ (and $2$) is strongly influenced by the complex
instability of the $x_2$ normal mode. In case $0$ this mode is stable so
that $H_2(t)$ will vary even less. Using the same initial conditions for
case $1$ we find strong variations of $H_2(t)$, but always within the
limits of the error estimates; see fig. \ref{H2case0-1x2}.

\subsection*{Case $2$ with less-balanced masses} 
\begin{figure}[htp]
\begin{center}
\resizebox{14cm}{!}{
\includegraphics{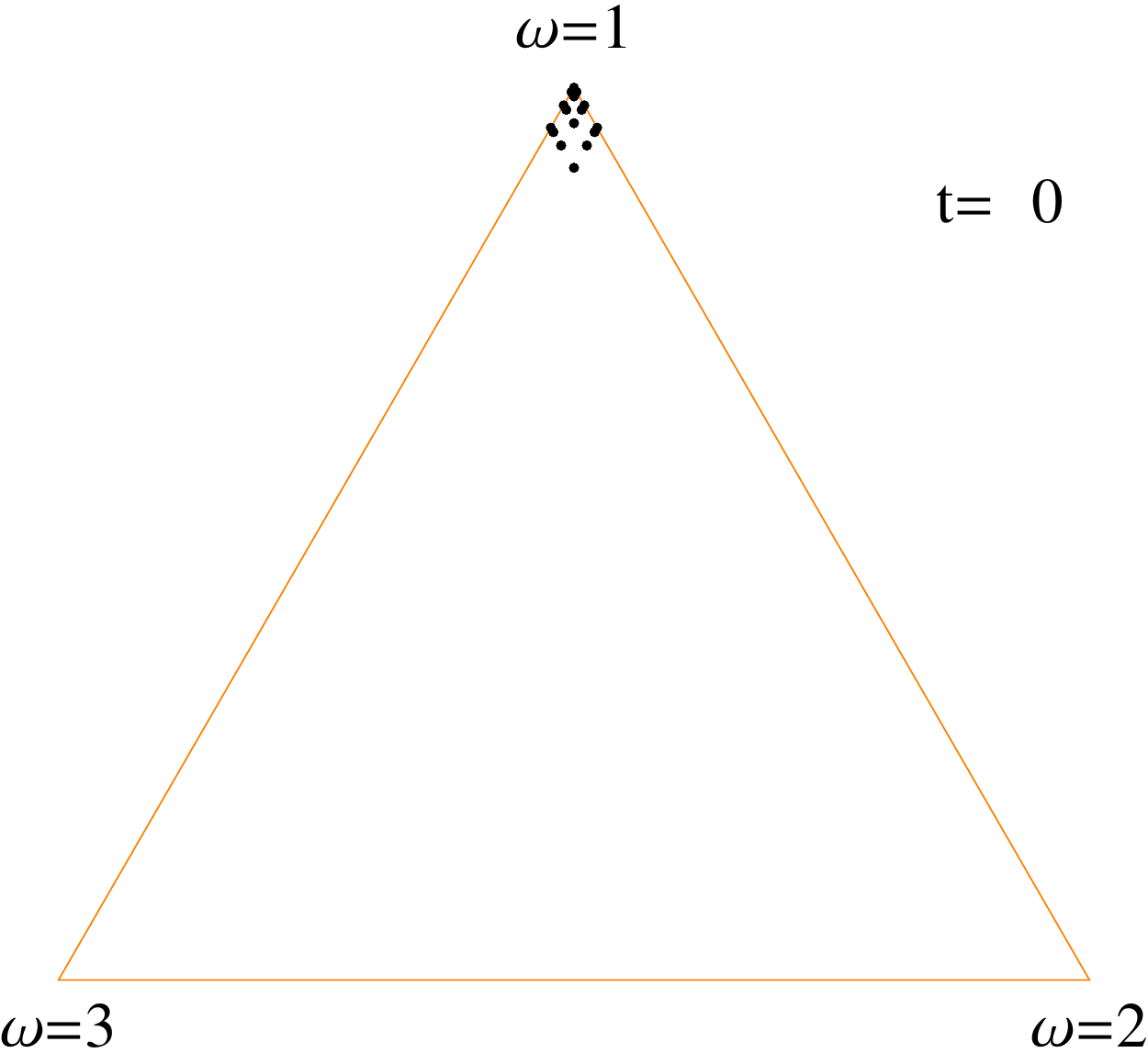}\,\,\includegraphics{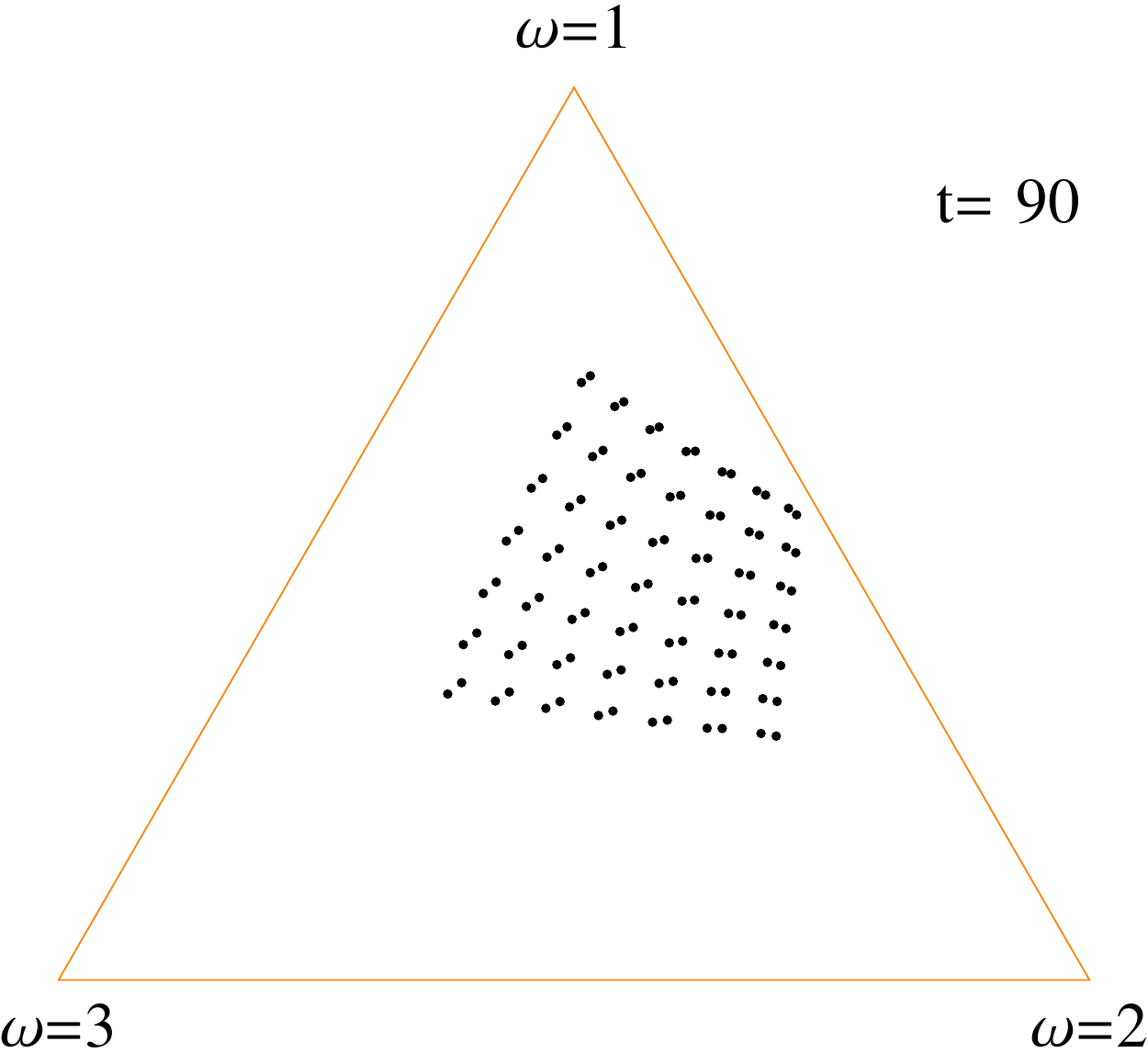}\,\,\includegraphics{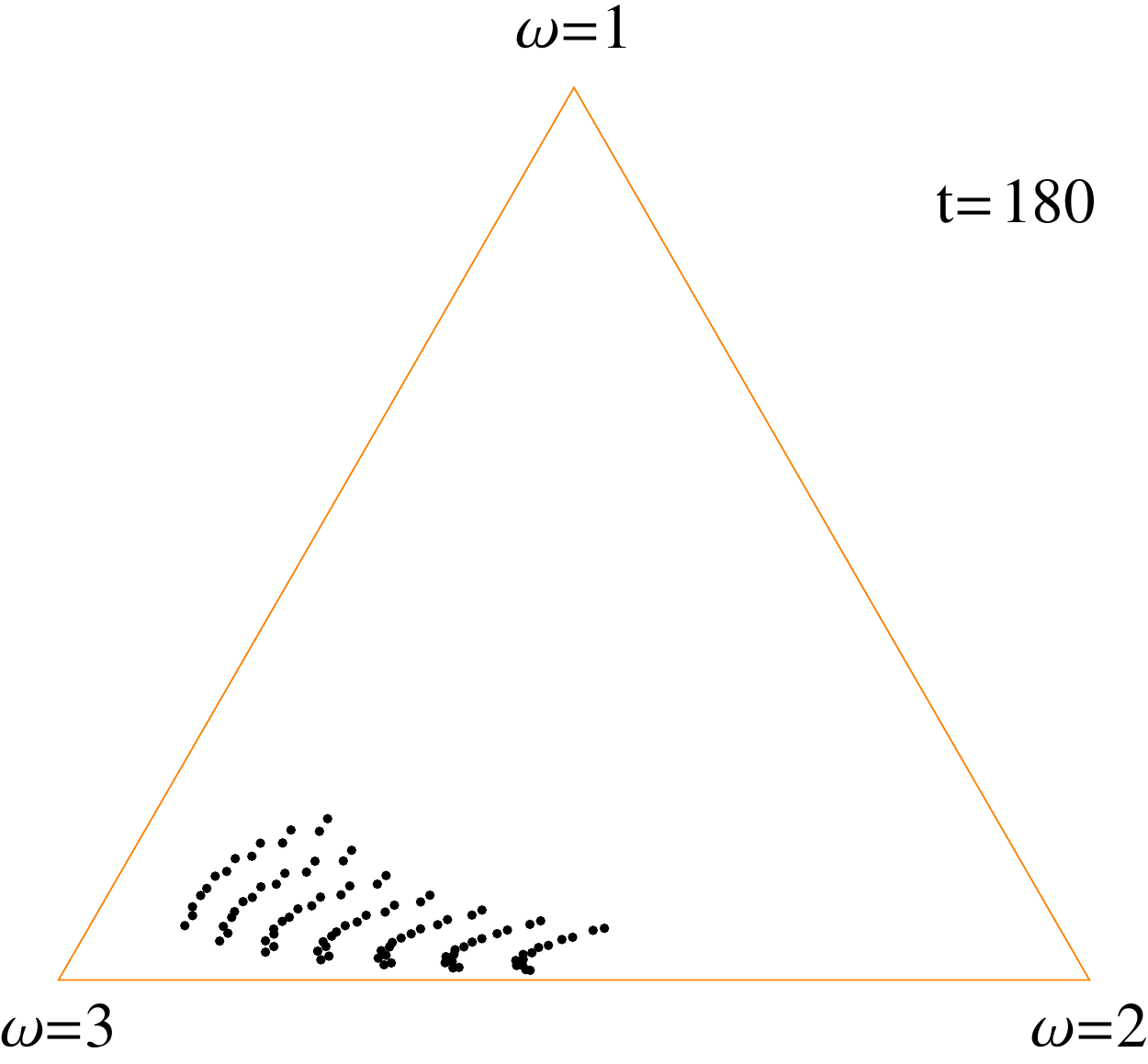}}\\
\resizebox{14cm}{!}{
\includegraphics{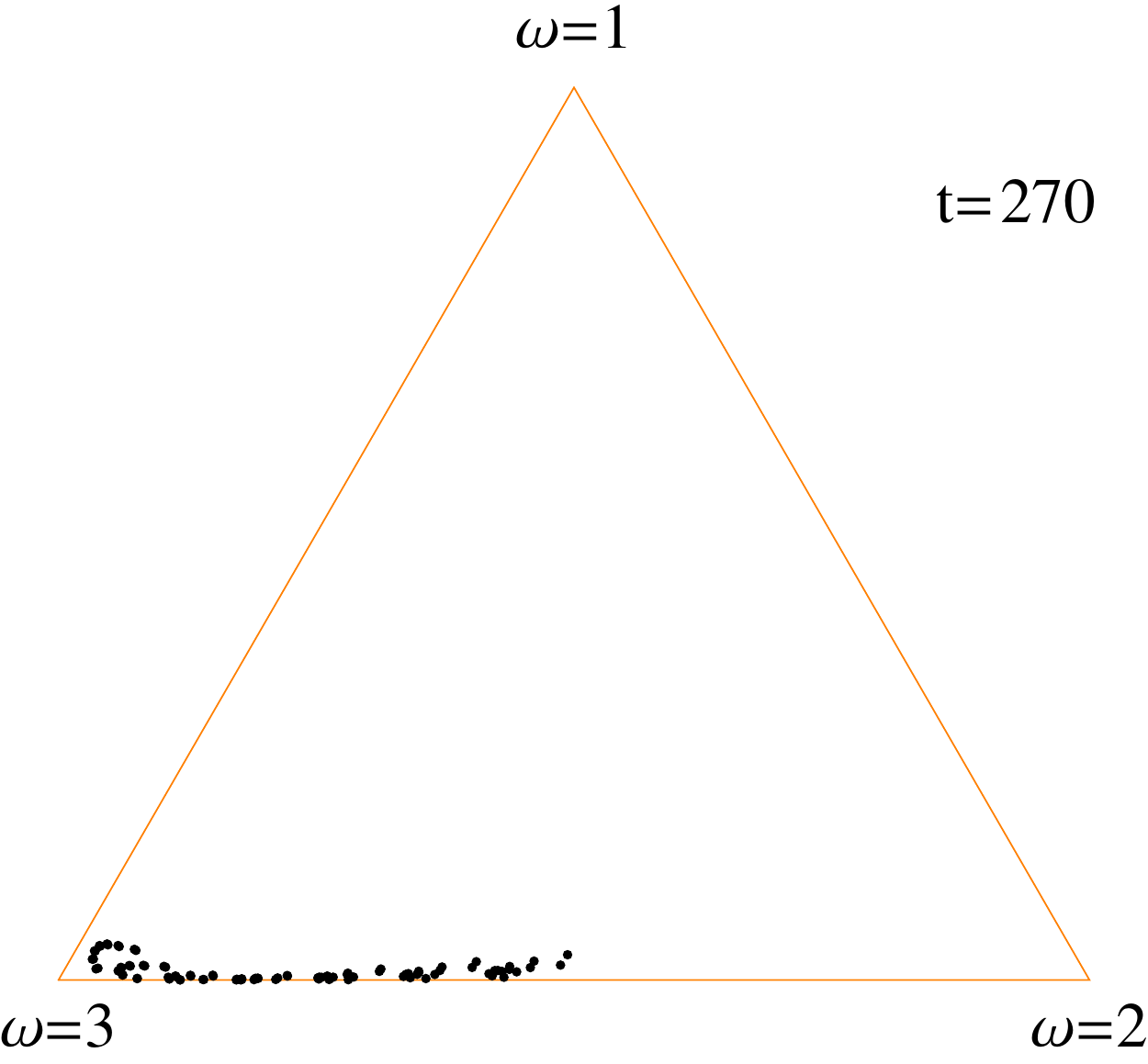}\,\,\includegraphics{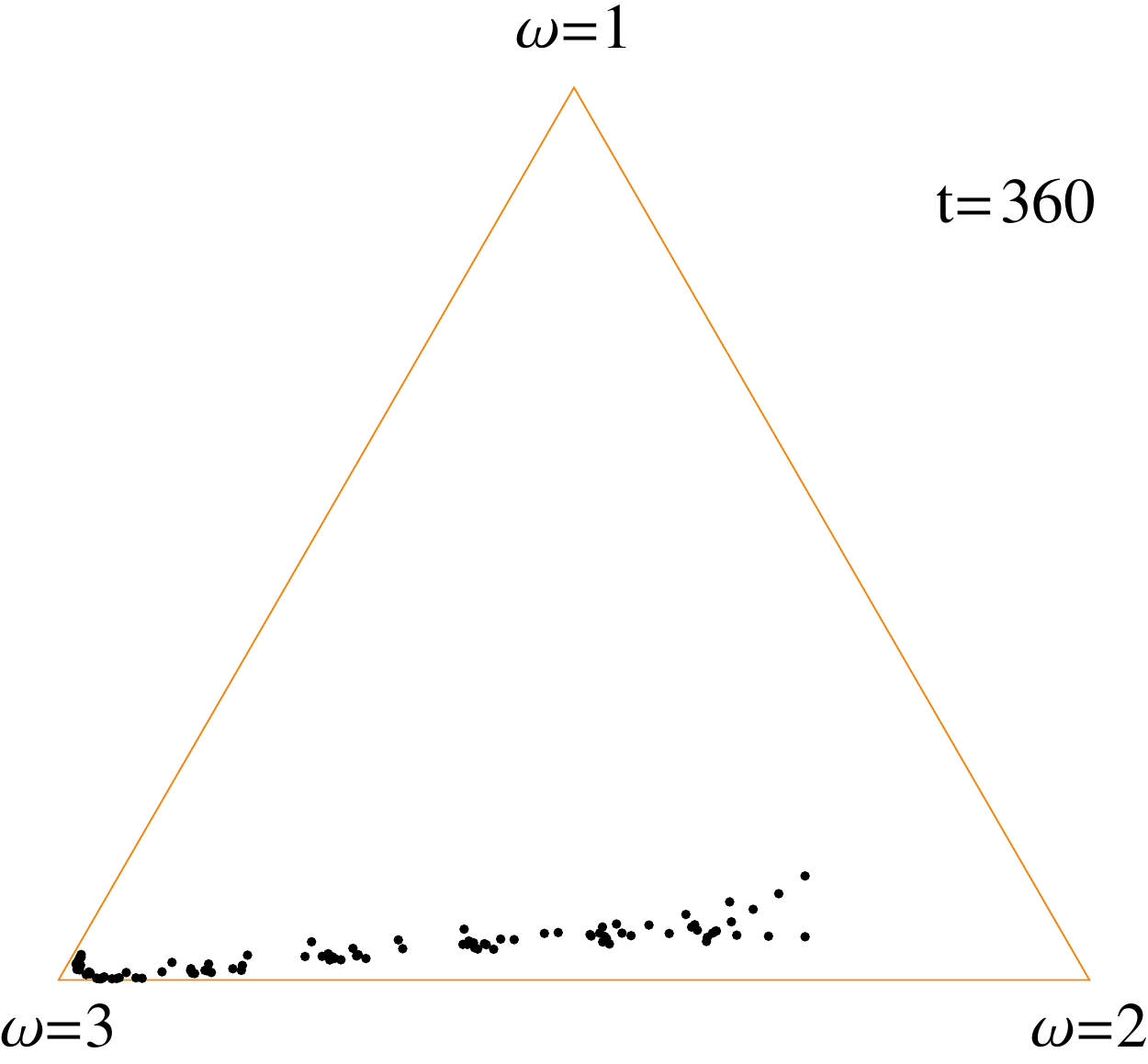}\,\,\includegraphics{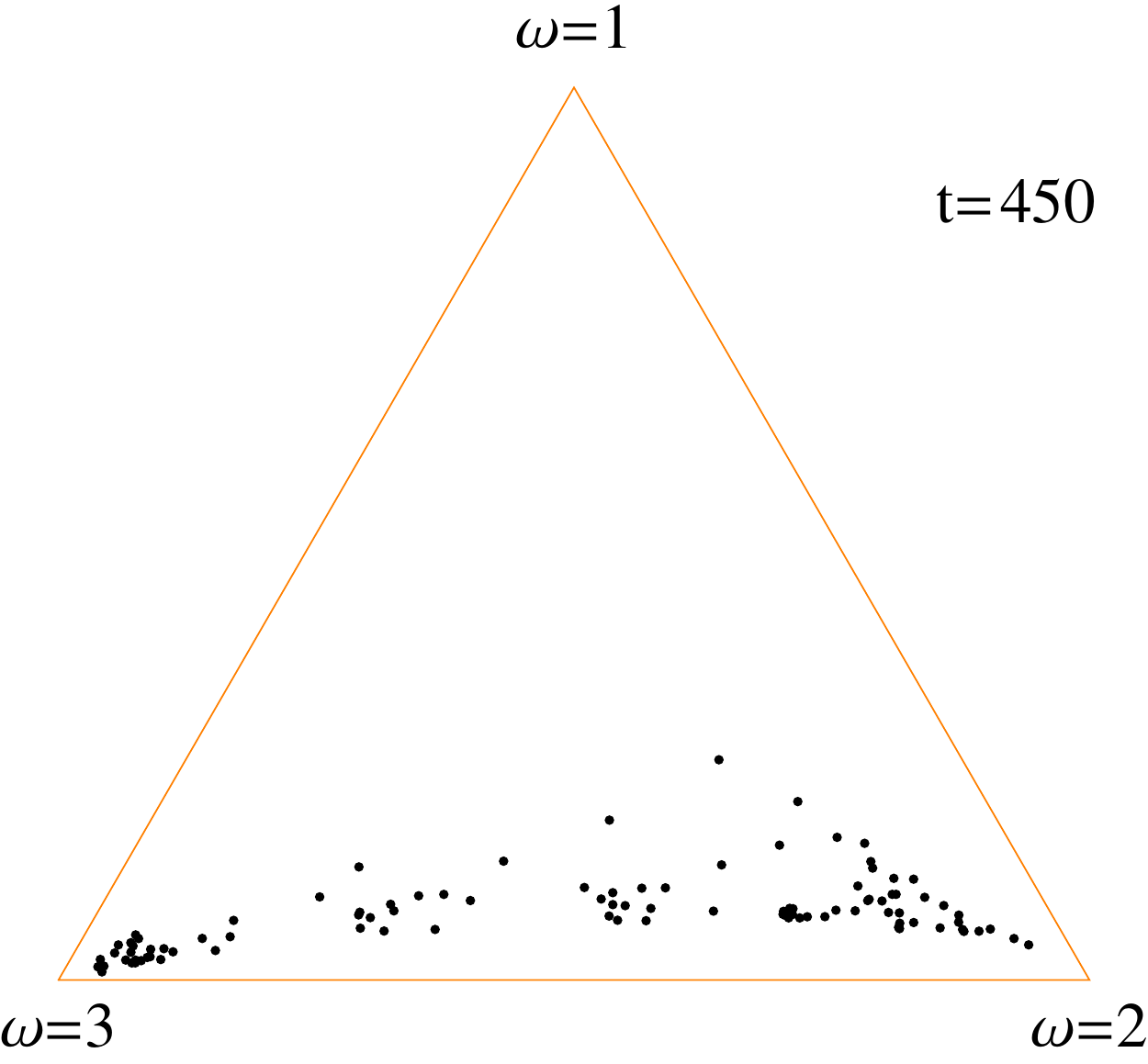}}
\end{center}
\caption{ We consider for case $2$ the time evolution of $98$ starting
points near the $\omega =1$ vertex
by displaying the action-simplex at various times. \label{simcase2-3}}
\end{figure}

We choose for
$u = 0.826713$ from eqs. (\ref{a13}) and (\ref{a24}) a case with even less
balanced masses; in this case $m_1$ is quite massive. We have:
\[ a_1= 0.000685158, a_2= 0.11239, a_3= 0.100269, a_4= 0.286656. \]
With these mass ($a_i$) values the symplectic transformation of subsection
\ref{sect-hamtr} to system (\ref{nleqsHam2})
produces the expression:
\begin{eqnarray*}
H_3 & = 0.0352657 x_1^3 -0.0272316 x_1^2x_2 -0.0743155x_1^2x_3-
0.0366184x_1x_2^2- 0.00260064x_1x_3^2 \\ & -0.0337877x_1x_2x_3+
0.0181144x_2^3 + 0.000760425x_3^3- 0.0105601x_2x_3^2 + 0.023904x_2^2x_3.
\end{eqnarray*}
We have the case:
\[ d_6= -0.0337877,\, d_9= -0.0105601 \]
If $d_9 \neq 0$ (the cases $1$ and $2$), the $x_3$ normal mode does not
exist. In fig. \ref{simcase2-3} we show the action-simplex for solutions
starting near the $x_1=x_2=0$ position, so near the $\omega =1$ vertex.

\begin{figure}[htp]
\begin{center}
\resizebox{14cm}{!}{
\includegraphics{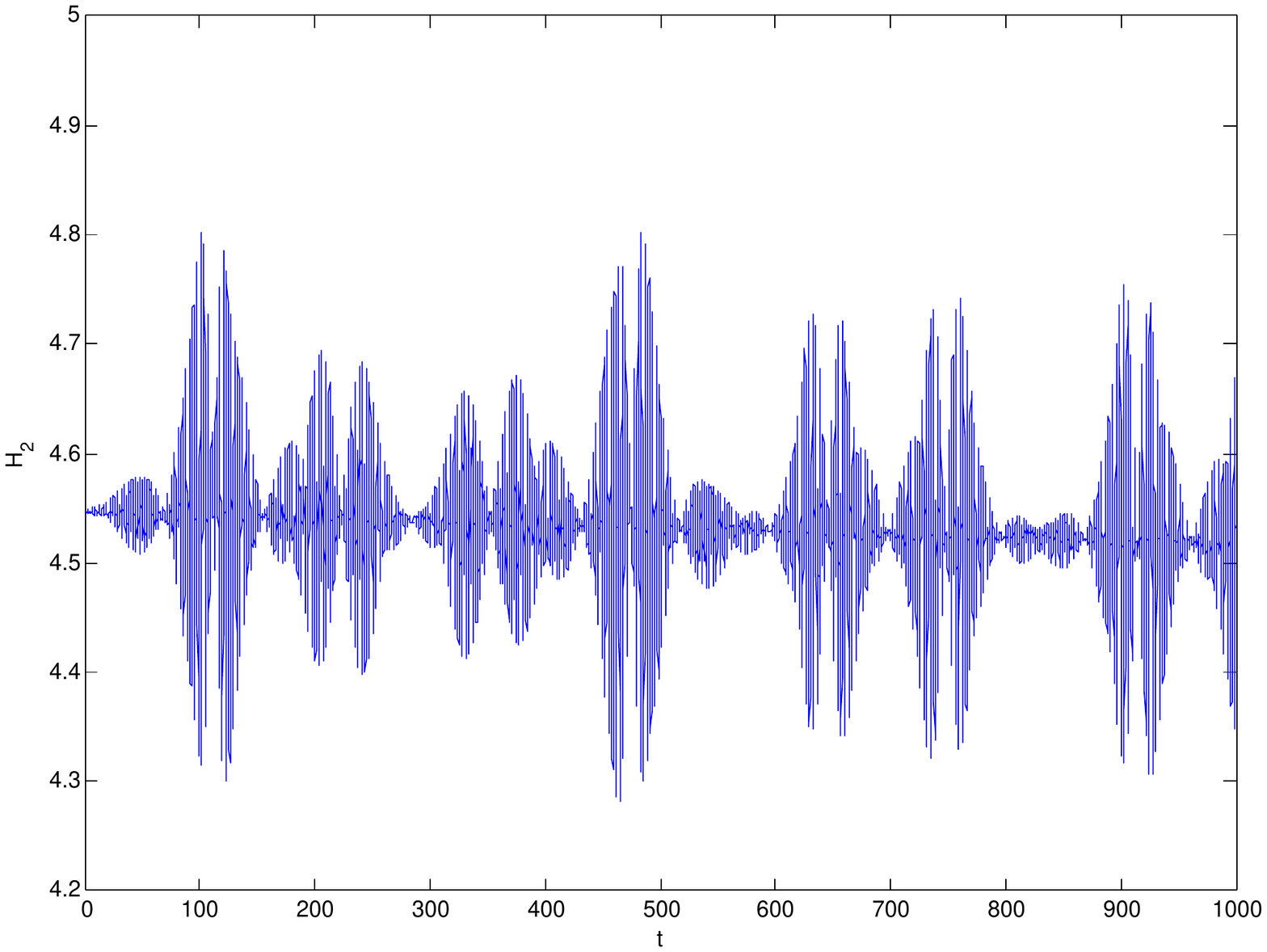}}
\end{center}
\caption{The $H_2(t)$-time series $0-1000$ based on system
\eqref{nleqsHaminter}, case $2$, with initial values
$x_1= 0.1, x_2=1.5 ,x_3= 0.1, 
\dot{x}_1= \dot{x}_2 = \dot{x}_3=0, \varepsilon = 0.5$. \label{H2eps3}}
\end{figure}

 We present the $H_2(t)$ time series based based on system
 \eqref{nleqsHaminter} in fig. \ref{H2eps3}.


\subsection{Comparison with another Hamiltonian system in $1:2:3$ resonance}
We will discuss our results for the inhomogeneous FPU chain with another
Hamiltonian system in $1:2:3$ resonance, and compare the instability types
of the $x_2$ normal mode.

For the inhomogeneous FPU lattice in $1:2:3$ resonance we found complex
instability ($C$) of the $x_2$ normal mode and no cases of $HH$
instability. Both cases, $HH$ and $C$ lead to a non-integrable normal form
but the dynamics is different. See \cite{OC}.

\begin{figure}[htp]
\begin{center}
\resizebox{16cm}{!}{
\includegraphics{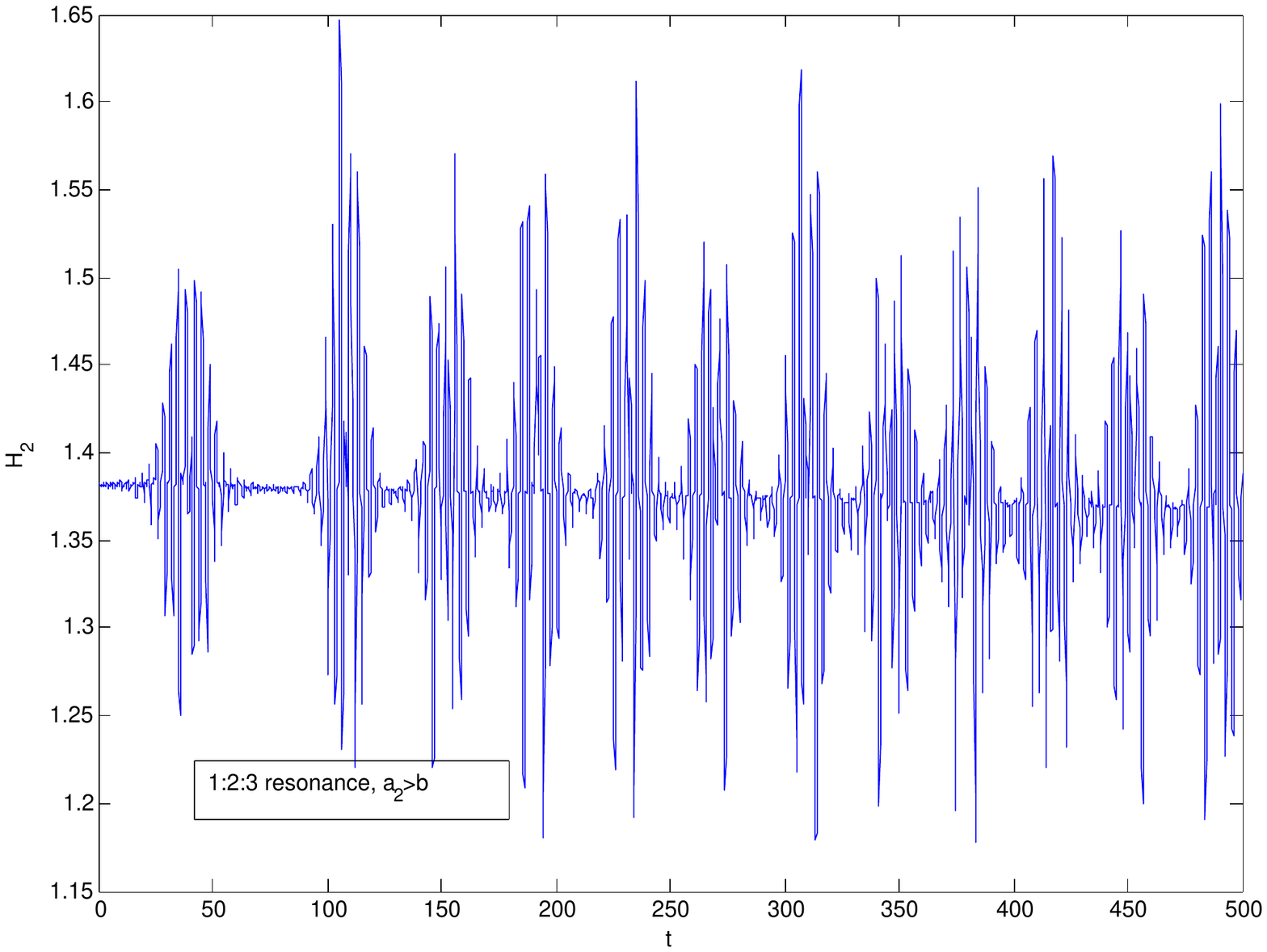}\,\, \includegraphics{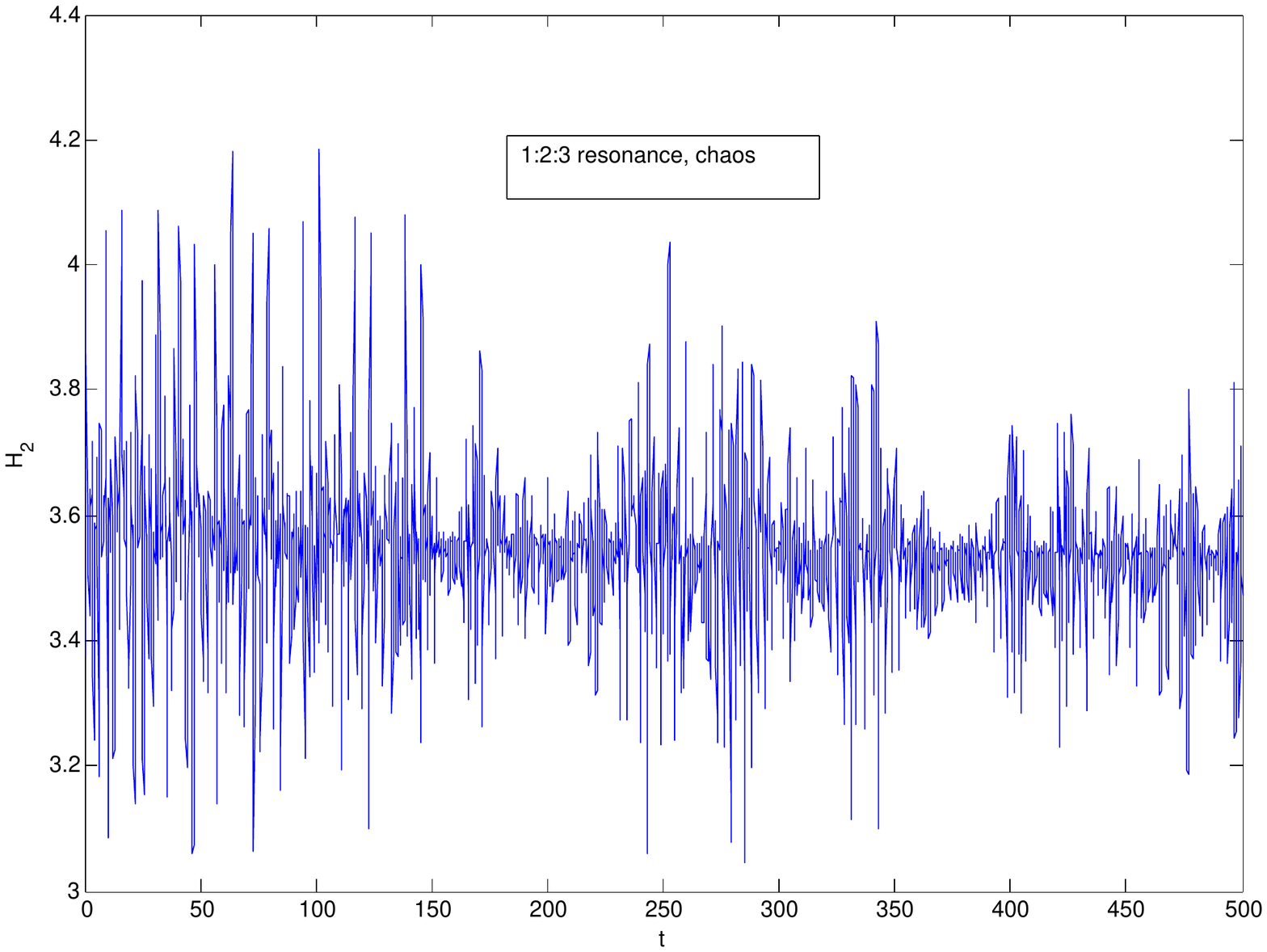}}
\end{center}
\caption[]{ Two $H_2(t)$ time series based on Hamiltonian (\ref{hamHHC}). On
the left $x_1(0)=0.1,\; x_2(0)=1,\; x_3(0)=0.5$, and on the right
$x_1(0)=2,\; x_2(0)=1\; x_3(0)=-.05$. For both time series we use
$\varepsilon = 0.5, a_2=3,a_3=1, b=1$, and
$\dot x_1(0)=\dot x_2(0)=\dot x_3(0)=0$.

For the $x_2$ normal mode we have instability $HH$ on the left and
instability $C$ on the right. In both cases the Hamiltonian flow is chaotic
but in the right picture the system has undergone Devaney-Shilnikov
 bifurcation.

 Horizontal scales: time in $[0,500]$, vertical scales: energy in
 $[1.15,1.65]$ (left) and in $[3,4.5]$ (right). \label{HHC}}
\end{figure}

To illustrate the different dynamics consider the Hamiltonian presented as
an example in \cite{FV15}:
\begin{equation} \label{hamHHC}
H(p, q)= \frac{1}{2}(p_1^2+q_1^2)+(p_2^2+q_2^2)+\frac{3}{2}(p_3^2+q_3^2)-
\varepsilon q_1^2(a_2q_2+a_3q_3)- \varepsilon b q_1q_2q_3.
\end{equation}
This system is in $1:2:3$ resonance but it is not derived from a FPU chain.
We present $H_2(t)$ for both cases in fig. \ref{HHC}. The dynamics is
chaotic but in the case left, the $q_2$ normal mode is unstable with real
eigenvalues (HH); transverse homoclinic intersections produce chaotic
motion. On the right the $q_2$ normal mode is complex unstable (C) which
produces the Hamiltonian Devaney-Shilnikov phenomenon. This involves a
homoclinic orbit surrounded by an infinite number of unstable periodic
solutions producing more violent chaotic motion as predicted in \cite{D76}.

\section{Conclusions}
\begin{itemize}
\item[] {\bf General}
\begin{itemize}\item For an inhomogeneous periodic FPU-chain with four
particles, most frequency ratios occur for a one-dimensional variety of
mass ratios. The frequency ratios $1:2:1$ and $1:1:3$ arise for a finite
number of mass ratios, the ratios $1:2:2$, $1:1:1$ and $1:3:3$ do not occur
at all in this FPU-chain. See table \ref{tab-res}.

\item For any number of particles $n\geq 3$ the set of mass distributions
   for a given frequency distribution has a relatively simple algebraic
   structure. For $n=4$ we describe algorithmically how to determine this set
   for a given frequency distribution. For $n\geq 4$ there are frequency
   distributions that do not correspond to any mass distribution.

\end{itemize}

\item[]{\bf The case of four particles in $1:2:3$ resonance}
\begin{itemize}

\item A special case of the resonance $1:2:3$ has the symmetry of two equal
masses and two quite different masses. Along the variety of mass ratios as
a limit case one of the masses tends to infinity.

\item The symmetric case of two equal masses differs dynamically from the
other cases. The transition corresponds to a Hamiltonian-Hopf bifurcation
with a Shilnikov-Devaney bifurcation producing chaotic dynamics. In a more
general context such behaviour of the $1:2:3$ resonance was described in
\cite{HV}.

\item The normalized system for the symmetric case of two equal masses is
integrable and has periodic solutions for each of the three eigenmodes (the
normal modes). Moreover, there are on the energy manifold two families of
periodic solutions connecting the second and the third eigenmode. This is a
degeneration in the sense described by Poincar\'e \cite{PMC}, vol. 1.

\item Under the transition away from the symmetric case, the eigenmodes
$x_1$ (associated with frequency $3$) and $x_2$ (associated with frequency
$2$) produce a periodic solution (normal mode) in the nonlinear system. The
periodic solution that was associated to the third eigenmode in the
symmetric case moves away along an edge of the action simplex. The two
continuous families of periodic solutions of the symmetric case break up
into four periodic solutions.

\item The inhomogeneous periodic FPU $\alpha$-chain with four particles is
characterized by a non-integrable normal form, except in the symmetric case
of two equal masses. The implication is that near stable equilibrium its
chaotic behaviour is not restricted to exponentially small sets as in the
case of two dof systems and as in the case of the classical FPU
$\alpha$-chain. In this sense the model of the classical FPU-chain is
misleading.

\end{itemize}
\end{itemize}

\appendix

\newcommand\een{\mathbf{e}}
\newcommand\z{\zeta}

\section{Further details for the spectrum induced by $H_2$} Here we give a
further discussion of results mentioned in section~\ref{br}.

\subsection{Fiber contained in a quadratic set}\label{sect-ell}For given
eigenvalues $\ld_1\geq \ld_2\geq \cdots\geq \ld_{n-1}>\ld_n=0$ of $A_n C_n$
(see subsection \ref{sect-mev})
we have the relations \eqref{pldrel} for the elements $(a_1,\ldots,a_n)$ in
the corresponding fiber. Here we will use the last two relations, given in
\eqref{pe-eq}.

\begin{lemma}
The polynomials $p_{n-1}$ and $p_{n-2}$ have the form indicated
in~\eqref{pe-eq}.
\end{lemma}
\prf If we replace the entries $-1$ at positions $(1,n)$ and $(n,1)$ in
$C_n$ by $0$ we obtain the Cartan matrix $\mathcal{C}_n$ for the root
system of type $A_n$. (See, eg., \cite[D\'ef. 3 in 1.5 of Chap. 6, and
Planche I]{Bo}.)
The determinant of $\mathcal{C}_n$ is known to be $n+1$.

If all $a_j$ are non-zero, the characteristic equation is equivalent to
$\det(C_n - \ld A^{-1})=0$. We determine first the factor of
$(-\ld)^{n-1}$. In the expansion of the determinant the term with $\ld$ at
all diagonal positions except at $(j,j)$ is equal to
\[ 2 \prod_{i\neq j} (-\ld a_i^{-1} )= 2 (-\ld)^{n-1} a_j/(a_1a_2\cdots
a_n)\,.\]
So the factor of $(-\ld)^{n-1}$ in $\det(AC_n - \ld I_n)$ is
$\sum_j 2 a_j  = p_{n-1}(a)$.

For the factor of $(-\ld)^{n-2}$ we have contributions of two types: Two
diagonal positions $j$ and $j+1$ (modulo $n$) lead to a contribution of the
form
$\det(\mathcal{C}_2)\allowbreak\, \prod_{i\neq j,j+1}(-\ld)/\allowbreak a_i$.
Two non-adjoining diagonal positions $j_1$, $j_2$ contribute
$2\cdot 2 \sum_{i\neq j_1,j_2} (-\ld)/a_i$. This leads to the description
of $p_{n-2}(a)$. \qed

By scaling we arrange that the vectors $(\ld_1,\ldots,\ld_{n-1},0)$ of
eigenvalues of $A_n C_n$ satisfy $\sum_{j=1}^{n-1}\ld_j=1$, and we put
$\eta= e_2\bigl( \ld_1,\ldots,\ld_{n-1}\}\bigr)$. Then the points of the
fiber of a given vector of eigenvalues are elements of the following set
$Q_\eta$:

\begin{prop}\label{prop-quadr}Let $n\geq 3$. For given $\eta>0$ denote by
$Q_\eta$ the set of points $a\in \RR^n$ satisfying
\begin{equation} p_{n-2}(a)= \eta\,, \quad
p_{n-1}(a)=\frac12\,.\end{equation}
Then
\begin{enumerate}
\item[a)] If $\eta<\frac12-\frac 3{4n}$, then $Q_\eta$ is a compact quadric
in the hyperplane $a_1+\cdots+a_n=\frac12$ in $\RR^n$ with a non-empty
intersection with $\RR_{>0}^n$.
\item[b)] If $\eta=\frac12-\frac 3{4n}$, then $Q_\eta$ consists of one point
in $\RR_{>0}^n$.
\item[c)] If $\eta>\frac12-\frac 3{4n}$, then $Q_\eta=\emptyset$.
\end{enumerate}
\end{prop}
\prf Let $P=C_n - 6I_n+4 E$, with $E $ the $n\times n$-matrix with all
elements equal to~$1$. Then, considering $a=(a_1,\ldots,a_n)$ as a row
vector, we have
\[ q_{n-2}(a) =\frac12 a P a^T\,.\]
To see this we check that $P$ is the matrix
\[ \begin{pmatrix} 0& 3& 4& \cdots& 4 & 3\\
3& 0 & 3 & \cdots & 4 & 4\\
4& 3 & 0 & \cdots & 4 & 4\\
\vdots & \vdots & \vdots & \ddots & \vdots & \vdots\\
4& 4& 4& \cdots & 0 & 3\\
3& 4&4& \cdots & 3& 0
\end{pmatrix}
\]
There are orthogonal matrices $U$ such that $U^T C_n U= \Ld$, where $\Ld$ is
the diagonal matrix with the eigenvalues $\ld_j$ of $C_n$ on the diagonal.
We put the eigenvalue $0$, with eigenvector $\een$ as the last one. Then
$\een=(0,\ldots,0,\sqrt n) U^T$.
This gives
\begin{align*}
p_{n-1}(a) &= 2 \sum_j a_j = 2 a \een^T = 2\sqrt n (a U)_n\,,\\
p_{n-2}(a) &= \frac12 a \bigl( C_n-6I_n )a^T + 2 a E a^T = \frac12 a U
\bigl( \Ld-6I_n\bigr)U^T a^t + \frac12 \bigl( p_{n-1}(a)\bigr)^2\,.
\end{align*}

The points in the hyperplane $p_{n-1}(a)=1$ can be described as
\[ a = \Bigl(x_1,x_2, \ldots, x_{n-1}, \frac1{2\sqrt n}\Bigr)
U^T \,.\]
We write $x=(x_1,\ldots,x_{n-1})^T$. We find the equation
\begin{equation}\label{ellipsoid}
\begin{aligned} \eta &= \frac12 (x,1/(2\sqrt n)) (\Ld -6 I)
(x,1/(2\sqrt n))^T +\frac12\\
&= - \sum_{j=1}^{n-1} \frac{6-\ld_j}2 \, x_j^2 + \frac12 -\frac 3{4n} \,.
\end{aligned}
\end{equation}
So the points run through a quadratic set in the hyperplane $p_{n-1}(a)=1$.
The eigenvectors of $C_n$ can be chosen as $(\z^k,\z^{2k},\cdots,\z^{nk})$
with $\z=e^{2\pi i/ n}$, which leads to eigenvalues
$2-2\cos2\pi k/n\in [0,4]$. So the $\ld_j-6$ are strictly negative. The
equation becomes
\begin{equation} \sum_{j=1}^{n-1} \frac{6-\ld_j}2 x_j^2 = \frac12-\frac
3{4n} - \eta\,.\end{equation}
In case b) in the proposition the single point $x=0$ corresponds to
$\frac1{2n}\een\in \RR_{>0}^n$. As $\eta$ decreases the quadric expands in
all directions, some of these stay inside $\RR_{>0}^n$. \qed

\begin{corol}\label{cor-im}If $n=3$ each choice of eigenvalues
$\lambda_1\geq \lambda_ 2 >\lambda_3=0$ of $A_3C_3$ occurs for some
positive diagonal matrix~$A_3$.

If $n\geq 4$, there are choices of eigenvalues
$\ld_1\geq \ld_2\geq \cdots\geq \ld_{n-1}>\ld_n=0$
for which there are no positive diagonal matrices $A_n$ such that $A_n C_n$
has these eigenvalues.
\end{corol}
\prf The choice $\ld_1=\cdots=\ld_{n-1}=\frac1{n-1}$ leads to
\[\eta=e_2\bigl( \{\ld_1,\ldots,\ld_{n-1} \}\bigr)=\binom{n-1}2 \bigm/
(n-1)^2 = \frac 12 - \frac1{2(n-1)}\,.\]
This is at most $\frac12-\frac 3{4n}$ if $n=3$. This establishes the second
assertion.

For $n=3$ we have $\ld_1+\ld_2=1$, hence
$\eta=\ld_1\ld_2\leq \frac 14 = \frac 12 - \frac 3{4\cdot 3}$.\qed

\subsubsection{Spherical coordinates. }\label{sect-spher} In the case $n=4$
we may take the orthogonal matrix in the proof of the proposition in the
form
\[ U= \begin{pmatrix}-\frac12&0&\frac{-1}{\sqrt 2}& \frac12\\
\frac12&\frac{-1}{\sqrt 2}&0& \frac12\\
-\frac12&0&\frac1{\sqrt 2}& \frac12\\
\frac12&\frac1{\sqrt 2}&0&\frac12
\end{pmatrix}\,,\]
corresponding to the eigenvalues $4,2,2,0$. This gives
\begin{equation}\begin{aligned} x_1&= \frac{-a_1+a_2-a_3+a_4}2\,,&\quad\
x_2&= \frac{a_4-a_2}{\sqrt 2}\,,\\
x_3&= \frac{a_3-a_1}{\sqrt 2} \end{aligned}\end{equation}
Points of the fiber give points on the ellipsoid
$x_1^2+2x_2^2+2x_3^2 = \frac 5{16}-\eta$. Then spherical coordinates $\psi$
and $\phi$ are determined by
\begin{equation}\label{pc1}
\begin{aligned}x_1&= \rho \sin \psi\,,\quad x_2= \frac\rho{\sqrt 2} \cos\psi
\cos \phi\,,
\quad x_3=\frac\rho{\sqrt2} \cos\psi \sin \phi\,,\\
\rho&= \sqrt{\frac5{16}-\eta}\,>\,0\,,\quad -\frac \pi 2\leq \psi
\leq\frac\pi 2\,,
\quad -\pi \leq \phi\leq \pi\,. \end{aligned}
\end{equation}
These are the spherical coordinates used in fig.~\ref{fig-ell123}.

\subsection{Conditions for the fibers to be non-empty}
\begin{figure}[bf, htp]
\begin{center}
\includegraphics[width=5cm]{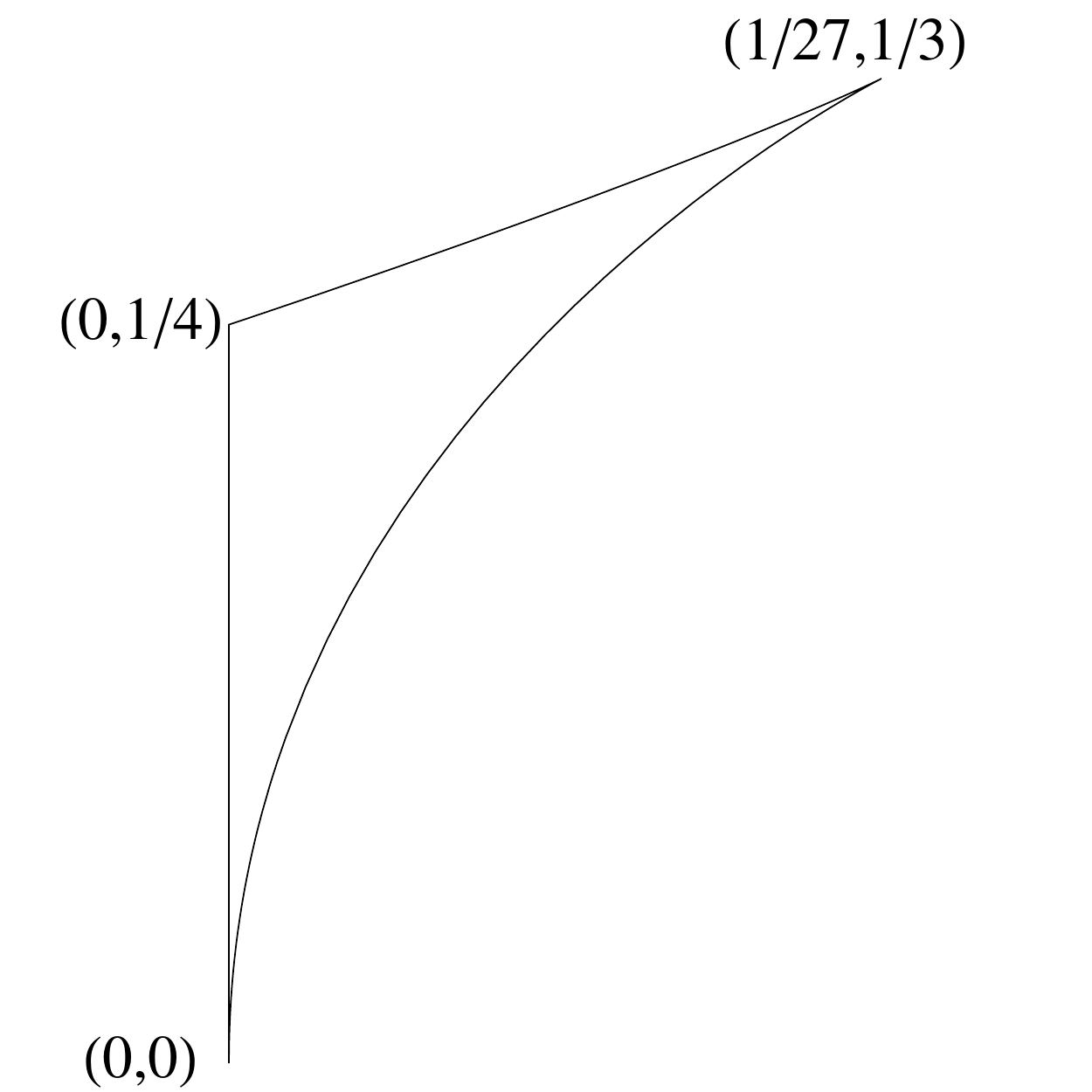}
\end{center}
 \caption{Region in the $\xi$-$\eta$-plane corresponding to choices of
 positive eigenvalues. (Horizontal axis: $\xi$; vertical axis:
 $\eta$.)}\label{fig-btim}
 \end{figure}

For $n=4$ the
equations \eqref{eq4} determine whether points of the fibers exist. In
particular, a (scaled) choice of eigenvalues determines $\xi,\eta>0$ which
determine the equations for the fiber. We first consider the values of
$(\xi,\eta)$ that can occur:
\begin{prop}\label{prop-ldim}Let $n=4$. The set of
$(\xi,\eta)=\bigl(e_3(\{\ld_1,\ld_2,\ld_3\}),e_2(\{\ld_1,\ld_2,\ld_3\})\bigr)$
where $(\ld_1,\ld_2,\ld_3)$ runs through the open triangle in $\RR^3_{>0}$
given by $\ld_1+\ld_2+\ld_3=1$, satisfy
\begin{equation} 0<\xi\leq\frac1{27}\,,\quad 0<\eta\leq \frac13\,,\quad
T(\xi,\eta)\leq 0\,, \end{equation}
where
\begin{equation}\label{Tdef} T(\xi,\eta) =
27\xi^2+4\eta^3-18\xi\eta-\eta^2+4\xi\,.\end{equation}
\end{prop}
\noindent\emph{Illustration }in fig.~\ref{fig-btim}.
 \prf We have to determine the image $X$ of the triangle
$T=\bigl\{(\ld_1,\ld_2,\ld_3)\in \RR_{>0}\;:\; \ld_1+\ld_2+\ld_3=1\bigl\}$
under the map
\[\Phi:(\ld_1,\ld_2,\ld_3) \mapsto(\xi,\eta) =
(\ld_1\ld_2\ld_3,\ld_1\ld_2+\ld_2\ld_3+\ld_3\ld_1)\,.\]

If a point $(\ld_1\ld_2,\ld_3)\in T$ is mapped to the boundary of the image
$X$, then the gradient of $\Phi$ has rank less than $2$ at that point. That
occurs if two of the coordinates are equal. By $S_{\!3}$-symmetry it
suffices to consider $\ld_2=\ld_3$. The image of the open segment
$\bigl\{(1-2y,y,y)\;:\; 0<y<\frac12\bigr\}$ consists of the points
\[ (\xi,\eta) = \bigl( y^2(1-2y),y(2-3y)\bigr)\,.\]
These are points of the curve $T(\xi,\eta)=0$. They run from $(0,0)$ to the
cusp at $\bigl(\frac1{27},\frac13\bigr)$ and then to
$\bigl(0,\frac14\bigr)$.

The boundary of $T$ consists of three segments, one of them
$\bigl\{ (0,x,1-x)\;:\: 0\leq x \leq 1\bigr\}$. The image is
$\bigl\{\bigl(0,x(1-\nobreak x)\bigr)\;:\;0\leq x\leq 1\bigr\}$, the
segment from $(\xi,\eta)=(0,0)$ to $\bigl(0,\frac14\bigr)$. By
$S_3$-invariance the two other boundary segments have the same image.

The image $X$ is the region enclosed by these boundary curves. \qed

The points $(\xi,\eta)$ for which the fiber is non-empty form a subset of
the region in Proposition~\ref{prop-ldim}. Corollary~\ref{cor-im} tells us
that the fiber is empty for $\bigl( \frac 1{27}, \frac 13\bigr)$. We give a
description of the set of $(\xi,\eta)$ corresponding to non-empty fibers. A
proof can be given along the same lines as that of
Proposition~\ref{prop-ldim}, but takes much more work. In the determination
of the fibers according to the computational scheme in the next subsection
it becomes clear anyhow whether the fiber is empty or not.

\begin{figure}[bf, htp]
\begin{center}
\includegraphics[width=7cm]{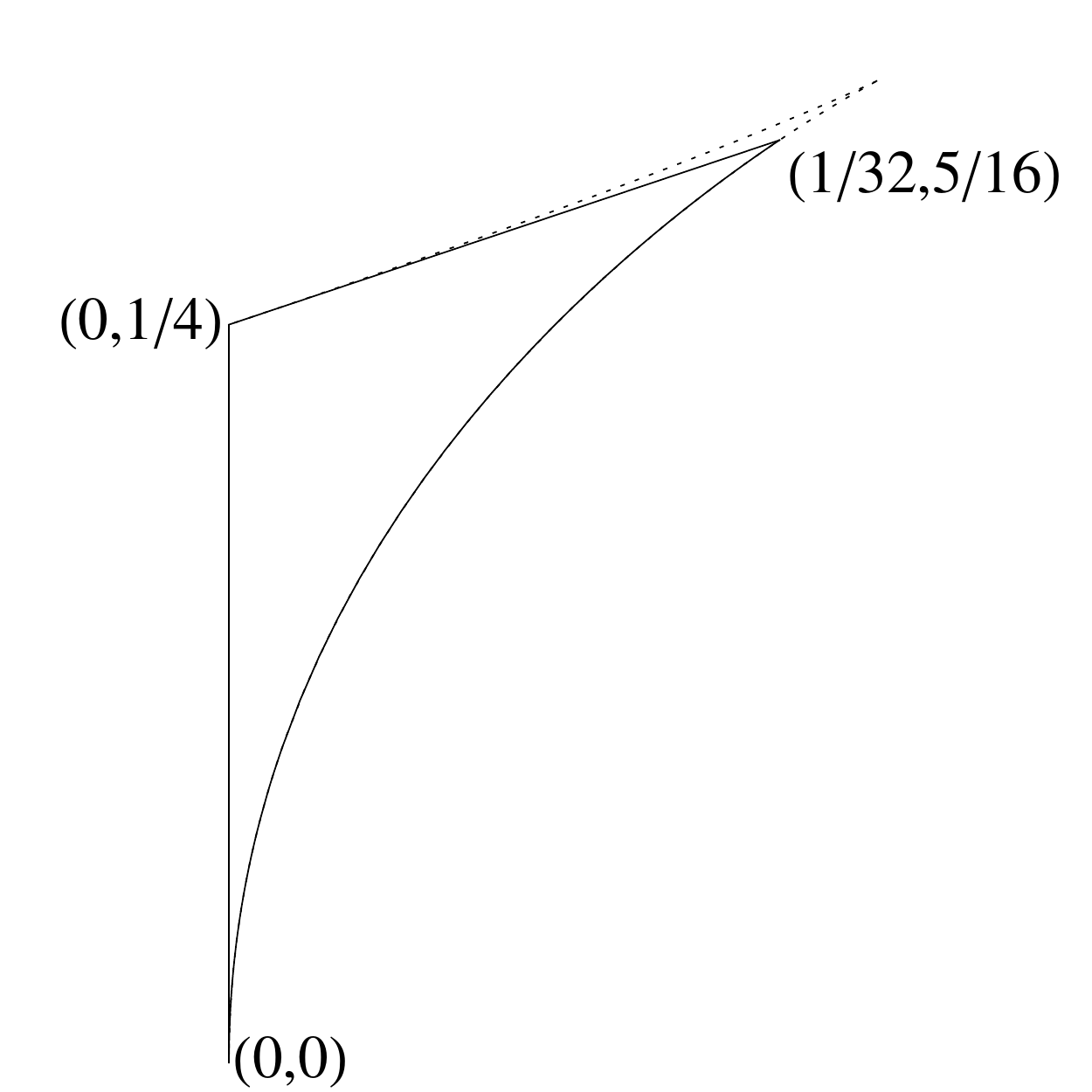}
\hspace{\fill}
\includegraphics[width=7cm]{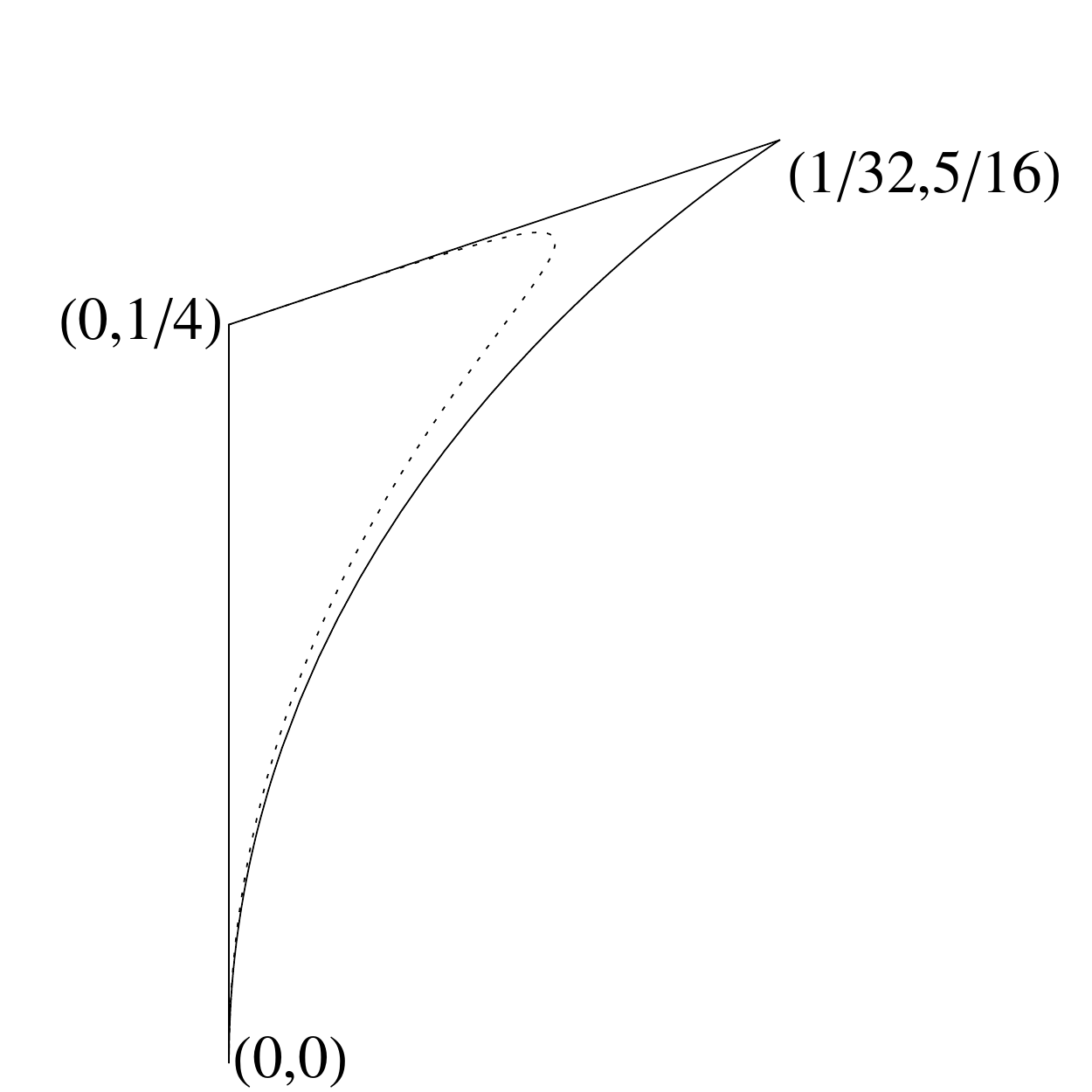}
\end{center}
\caption[]{ On the left the region in \eqref{imal} of points $(\xi,\eta)$
corresponding to non-empty fibers. The dotted line gives (part of) the
boundary of the larger region in Proposition~\ref{prop-ldim}, corresponding
to choices of positive eigenvalues. This shows that most of the possible
combinations $(\xi,\eta)$ correspond to a non-empty fiber.

On the right is again the region in \eqref{imal}, with the subregion in
\eqref{rncp} indicated by the dotted line. The points strictly to the right
of the dotted line correspond to compact fibers.} \label{fig-imal}
\end{figure}

\begin{prop}\label{prop-imal}The set of points $(\xi,\eta)$ corresponding to
a non-empty fiber is equal to
\begin{equation}\label{imal} \Bigl\{ (\xi,\eta)\in \RR_{>0}^2\;:\; 0<\xi\leq
\frac 1{32},\; 0<\eta \leq 2\xi+\frac14\,, \; T(\xi,\eta)\leq 0\Bigr\}\,,
\end{equation}
with $T$ as defined in \eqref{Tdef}.

The points $(\xi,\eta)$ for which the fiber is not compact constitute the
subset
\begin{equation}\label{rncp} \Bigl\{ (\xi,\eta) \in
\bigl(0,\frac1{32}\bigr)\times\bigl(0,\frac5{16}\bigr)\;:\;
8\xi^2+\eta^3-5\xi\eta-\frac14\eta^2+\frac 98 \xi\leq 0\Bigr\}\,.
\end{equation}
\end{prop}
\noindent\emph{Illustrations } in fig.~\ref{fig-imal}.

\subsection{Computation of fibers} The computation carried out in subsection
\ref{sect-comp123} for the resonance $(1:2:3)$ is guided by the use of the
action of the dihedral group $D_4$ on the solutions. We start with the
quantities $\xi,\eta,\eta_1,\eta_2$, which are invariant under the whole
group $D_4$.

In the next stage we consider the quantity
$4\left(-a_1+a_2-a_3+a_4\right)=\sqrt{1-16 \eta_2}$ which is invariant
under the subgroup $V_4\subset D_4$ generated by the permutations $[1,3]$
and $[2,4]$. This quantity is sent to its negative by $[1,2][3,4]$. The
group $V_4$ also leaves invariant $s_{13}, s_{24}, p_{13}, p_{24}$.
(If $\eta=4\xi+\frac 3{16}$ then we can take $\eta_2=\frac1{16}$. In that
situation $p_{13}$ and $p_{24}$ are not uniquely determined.)

In the next stage we determine $a_1$ and $a_3$, invariant under $[2,4]$ and
exchanged by $[1,3]$. Similarly $a_2$ and $a_4$ are invariant under $[1,3]$
and exchanged by $[2,4]$. The total solution $(a_1,a_2,a_3,a_4)$ is changed
by non-trivial elements of $D_4$, except in cases with additional symmetry.

\begin{table}[bf, htp] \fbox{\begin{minipage}{.95\textwidth}
\renewcommand\theenumi{\roman{enumi}}
\begin{enumerate}
\item For given positive eigenvalues put
\[
\eta=\frac{\ld_1\ld_2+\ld_2\ld_3+\ld_3\ld_1}{(\ld_1+\ld_2+\ld_3)^2}\,,\quad
\xi= \frac{\ld_1\ld_2\ld_3}{(\ld_1+\ld_2+\ld_3)^3}\,. \]
\item Write $\eta=4\eta_1+3\eta_2$. Determine the subinterval
$I_1\subset \bigl( 0,\frac 1{16}\bigr)$ such that $\eta_1>0$ for
$\eta_2\in I_1$.
\item Compute
\begin{align*} s_{13}&=\frac{1-\sqrt{1-16\eta_2}}4\,,&
\quad s_{24}&=\frac{1+\sqrt{1-16\eta_2}}4\,,\\
p_{13}&= \frac{\xi/4-s_{13}\,\eta_1}{\sqrt{1-16\eta_2}/2}\,,&\quad p_{24}&=
\eta_1-p_{13}\,.
\end{align*}
 Determine the subset $I_2\subset I_1$ such that $p_{13}>0$ and $p_{24}>0$
 for $\eta_2\in I_2$.
\item Compute $d_{13}=s_{13}^2-4p_{13}$ and $d_{24}=s_{24}-4p_{24}^2$.
Determine $I_3\subset I_2$ such that $d_{13}\geq 0$ and $d_{14}\geq 0$.
\item Compute $a_1,a_3 = \frac12\bigl(s_{13}\mp \sqrt{d_{13}}\bigr)$ and
$a_2,a_4= \frac12\bigl(s_{24}\mp \sqrt{d_{24}}\bigr)$. Determine the subset
$I_4\subset I_3$ such that $a_j>0$ for $j=1,\ldots,4$ for $\eta_2\in I_4$.
\item Apply all symmetries in the dihedral group $D_4$ to the points
$(a_1,\ldots,a_4)$.
\end{enumerate}\end{minipage}}
\caption{Instructions to compute fibers for the case $n=4$. In these
instructions we assume that $\eta \neq 4\xi+\frac 3{16}$. Otherwise we also
have to consider $p_{13}\in (0,\xi)$ and investigate whether this leads to
further solutions. } \label{tab-comp}
\end{table}

\begin{figure}[bf, htp]
\begin{center}
\includegraphics[width=7cm]{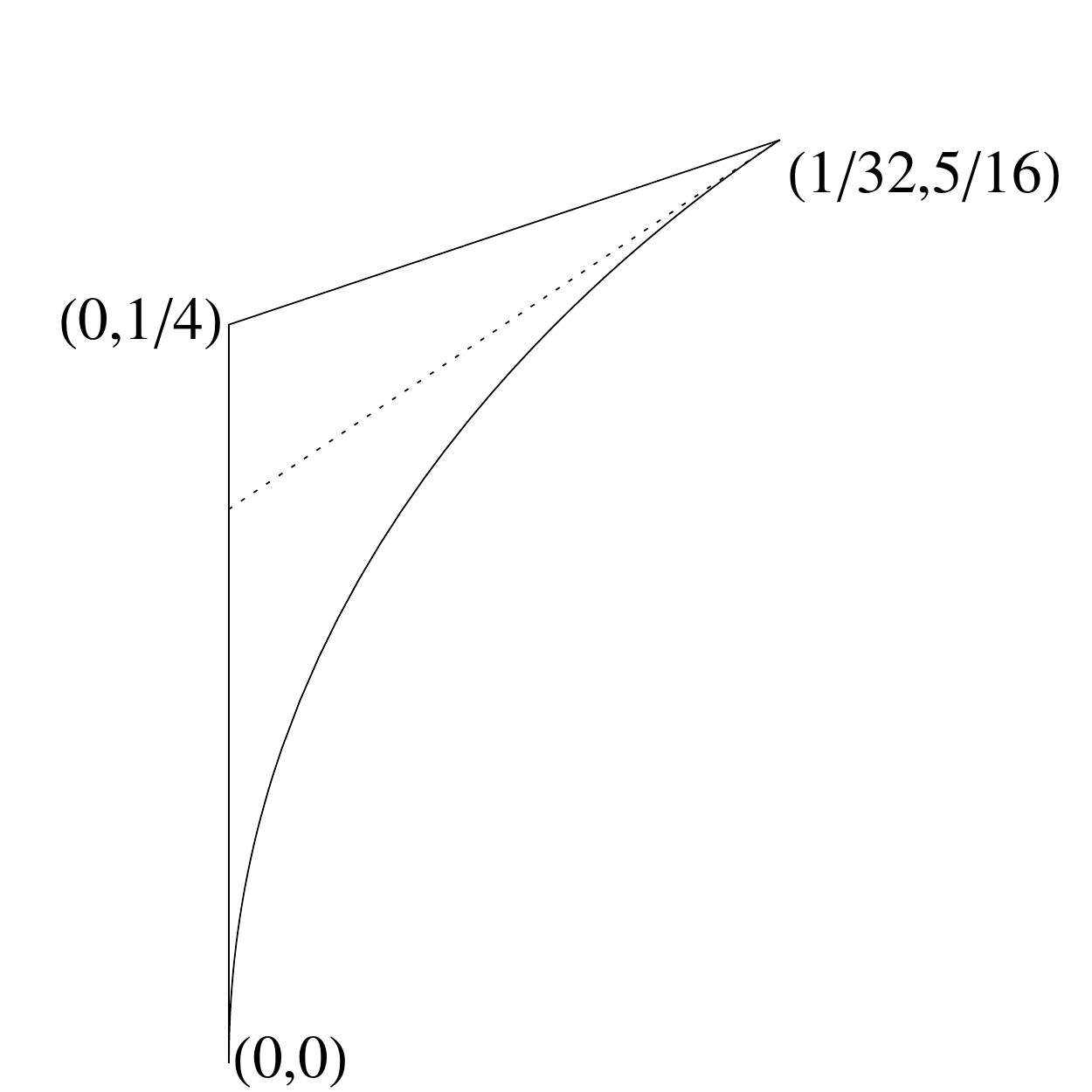}
\end{center}
\caption[]{The points on the dotted line $\eta=4\xi+\frac 3{16}$ between
$\bigl(0, \frac 3{16}\bigr)$ and $\bigl( \frac 1{32}, \frac 5{16}\bigr)$
correspond to fibers for which the computational scheme in
Table~\ref{tab-comp} is incomplete.

None of the resonances in Table~\ref{tab-res} correspond to points on this
exceptional line. } \label{fig-ec}
\end{figure}

In Table~\ref{tab-comp} the resulting computational scheme is described. It
works under the assumption that the point $(\xi,\eta)$ is not on the line
$\eta=4\xi+\frac 3{16}$, illustrated in fig.~\ref{fig-ec}. The parameter
$u$ was specially adapted to the resonance $(1:2:3)$. Here we use
$\eta_2\in \bigl(0,\frac1{16}\bigr)$ as the parameter.

We apply the computational scheme to the resonances $(1:2:2)$, $(1:1:2)$,
$(1:3:6)$ and $(2:3:4)$. Together with the resonance $(1:2:3)$ considered
in subsection \ref{sect-comp123} these are representative examples of the
cases in Table~\ref{tab-res}.
\begin{figure}[bf, htp]
\begin{center}
\includegraphics[width=8cm]{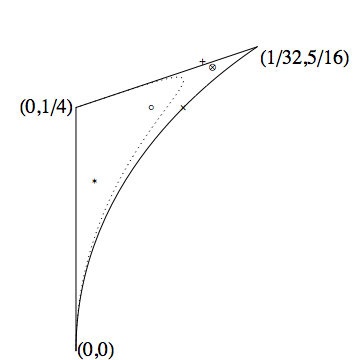}
\end{center}
\caption[]{Points corresponding to the resonances $(1:2:3)$ ($\circ$),
$(1:2:2)$ ($+$), $(1:1:2)$ ($\times$), $(1:3:6)$ ($\ast$), and $(2:3:4)$
($ \otimes$). } \label{fig-exsr}
\end{figure}

\subsubsection{Resonance $(1:2:2)$} To
$(\ld_1,\ld_2,\ld_3)=\bigl(\frac 49,\frac49,\frac19\bigr)$ corresponds
$(\xi,\eta) = \bigl( \frac{16}{729},\frac 8{27}\bigr)$. In
fig.~\ref{fig-exsr} it is hard to see whether it is in the region described
in~\eqref{imal}. A direct computation shows that $\eta>2\xi+\frac14$, so
the fiber is empty.

If we carry out the steps in the computational scheme, the set of values
that $\eta_2$ may have becomes empty when we check whether $d_{24}\geq 0$.

\subsubsection{Resonance $(1:1:2)$} \label{A112}
With $(\ld_1,\ld_2,\ld_3)=\bigl( \frac 23,\frac 16, \frac 16\bigr)$ we have
$(\xi,\eta)=\bigl( \frac 1{54}, \frac14 \bigr)$. The corresponding point
seems to be on the boundary of the region for a non-empty fiber. It turns
out that $T(\xi,\eta)$ is exactly~$0$.

\begin{figure}[bf, htp]
\begin{center}
\includegraphics[width=8cm]{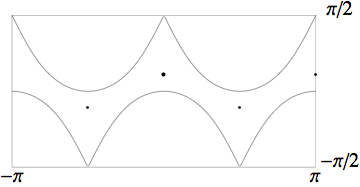}
\end{center}
\caption[]{The fiber for the resonance $(1:1:2)$ (subsection \ref{A112}) in
spherical coordinates as described in~subsection \ref{sect-spher}. The
thick point corresponds to the vector in \eqref{p112}, the other points are
its translates under elements of $D_4$. The curved line indicates the
boundary of the region with positive coordinates. } \label{fig-112}
\end{figure}

Following the computational scheme the expression for $d_{13}$ in terms of
$\eta_2$ turns out to be non-positive for
$\eta_1\in \bigl(0,\frac1{16}\bigr)$, with a zero only at
$\eta_2=\frac1{18}$. This leads to the solution
\begin{equation}\label{p112}
(a_1,\ldots,a_4) = \Bigl(\frac1{12},\frac{2-\sqrt 2}{12}, \frac1{12},
\frac{2+\sqrt 2}{12}\Bigr)\,.\end{equation}
It is invariant under the substitution $(13)$ in the dihedral group. See
fig.~\ref{fig-112}.

\subsubsection{Resonance $(1:3:6)$}For
$(\ld_1,\ld_2,\ld_3)=\bigl( \frac{18}{23},\frac9{46},\frac 1{46}\bigr)$ we
have $(\xi,\eta)=\bigl( \frac{81}{24334},\frac{369}{2116}\bigr)$. The
corresponding point in fig.~\ref{fig-exsr} is to the left of the dotted
line. This indicates that the fiber contains open curves.

\begin{figure}[bf, htp]
\begin{center}
\includegraphics[width=.95\textwidth]{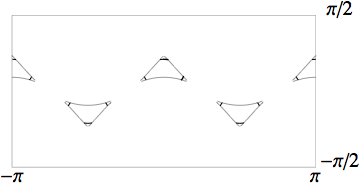}
\end{center}
\caption[]{The fiber for the resonance $(1:3:6)$ in spherical coordinates,
as described in~subsection \ref{sect-spher}.

The interior of the four triangles correspond to the region with positive
coordinates. The fiber consists of twelve open curves, three in each
triangle. The curves obtained with the computational scheme are draw
thicker than their translates under the dihedral group. } \label{fig-136}
\end{figure}

The computational scheme gives solutions for
\[ \eta_2 \in \bigl[ \frac{11}{1058}, h_2\bigr)\cup (h_3,h_4)\,,\]
with algebraic numbers $h_2\approx .112814$, $h_3\approx .0501346$,
$h_4\approx .0548411$. For $\eta_2=\frac{11}{10	 58}$ we find a point that
is invariant under $(13)\in D_4$. Fig.~\ref{fig-136} illustrates the fiber.

\subsubsection{Resonance $(2:3:4)$} For
$(\ld_1,\ld_2,\ld_3)=\bigl( \frac{16}{29},\frac 9{19},\frac4{29}\bigr)$ we
have $(\xi,\eta)=\bigl(\frac{576}{24389},\frac{244}{841}\bigr)$. The
corresponding point in fig.~\eqref{fig-exsr} is in the region where the
fiber is compact. With the relations in subsection \ref{sect-spher} one can
check that all $a_j$ are positive on the ellipsoid for
$\eta=\frac{244}{841}$.

\begin{figure}[htp]
\begin{center}
\includegraphics[width=9cm]{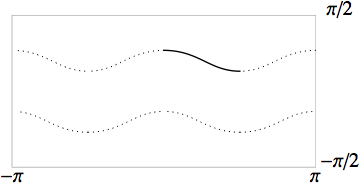}
\end{center}
\caption[]{The fiber for the resonance $(2:3:4)$ in spherical coordinates,
as described in subsection~\ref{sect-spher}.

The thick line corresponds to the solutions obtained by the computational
scheme. The dotted lines are formed by the translates under $D_4$ of the
computed part. In this case all points on the ellipsoid have positive
coordinates. } \label{fig-234}
\end{figure}

The computational schema gives a family of solutions depending on
$\eta_2 \in \bigl[ \frac{42}{841}, \frac{99}{1682}\bigr]$. The end points
give symmetric solutions: $a_1=a_3$ for $\eta_2= \frac{42}{841}$, and
$a_2=a_4$ for $\eta_2= \frac{99}{1682}$. In fig.~\ref{fig-234} we see that
the fiber consists of two closed curves.

\subsection{Transformation matrices for the resonance
$(1:2:3)$}\label{app-trm}
In subsection \ref{sect-comp123} we computed functions $u\mapsto a_j(u)$,
$1\leq j \leq 4$, on the interval $[0,u_1)$ as diagonal elements of a
diagonal matrix $A_4(u)$ such that $A_4(u)\, C_4$ has eigenvalues
$\frac 9{14}$, $\frac 4{14}$, $\frac 1{14}$, $0$. For the transformation to
eigenmodes of the Hamiltonian we need in subsection \ref{sect-hamtr} a
family of orthogonal matrices $u\mapsto U(u)$ such that $U(u)$ diagonalizes
$A_4(u)^{1/2}\, C_4\, A_4(u)^{1/2}$. For any value of $u$ such orthogonal
matrices can be found numerically. Here we want to describe explicitly the
dependence on $u\in [0,u_1)$. The version 9.0.1.0 of {\sc Mathematica} that
we used had difficulties with the symbolic computations. Hence we indicate
how we proceeded.

\begin{lemma}Let $A_4 $ be a positive diagonal matrix with diagonal elements
$a_1,\ldots,a_4$. Let $\ld $ be an eigenvalue of $A_4 C_4$ such that
$\ld\neq 2a_j$ for $1\leq j\leq 4$. Put
\[ \mu_j = \frac1{2-\ld/a_j}\,.\]
Then
\[ \bigl( \mu_1(\mu_2+\mu_4), \mu_2,\mu_3(\mu_2+\mu_4),\mu_4\bigr)\]
is an eigenvector of $A_4 C_4$ for the eigenvalue $\ld$.
\end{lemma}
\prf We have
\[ C_4-\ld A^{-1}=\begin{pmatrix}
\mu_1^{-1}&-1&0&-1\\
-1&\mu_2^{-1}&-1&0\\
0&-1&\mu_3^{-1}&-1\\
-1&0&-1&\mu_4^{-1}\end{pmatrix}\,.\]
We try to solve $(C-\ld A^{-1})v=0$ with $v=(p,x,q,y)$. The first and third
lines give $x+y=\mu_1^{-1}p=\mu_3^{-1}q$. Similarly, we get
$p+q=\mu_2^{-1}x=\mu_4^{-1}y$ Since $\ld$ is an eigenvalue of $A_4 C_4$
there are non-zero solutions, for which $x$ and $y$ both have to be
non-zero. So there is a solution with $x=\mu_2$. Then we obtain the vector
in the lemma.\qed

Now we take for $a_j$ the expressions in \eqref{a13} and \eqref{a24}. It is
clear that $2a_j(u)-\ld_i$ is not identically zero in $u$ for any of the
four eigenvalues $\ld_i$ and and any $j$. So we obtain vectors $v_i(u)$,
$1\leq i \leq 4$, that are eigenvectors of $A_4(u) \, C_4$ for the
eigenvalue $\ld_i$ for generic values of $u$.
\begin{table}[bf, htp] {\small
\begin{align*}
L_{1,1}&= \frac{\sqrt{u+6} \left(\sqrt{16-u}
   (20-3 (u-4) u)-18 \sqrt{2} \sqrt{5-u} \sqrt{6-u} \sqrt{u}\right)}{192
   \sqrt{35}
   (u-5)}\,,
   \displaybreak[0]\\
   L_{1,2}&=\frac{\sqrt{4-u} \left(\sqrt{2} \sqrt{6-u} (u (3 u-22)-20)+16
   \sqrt{5-u} \sqrt{-(u-16)
   u}\right)}{64 \sqrt{105} (u-5)}\,,
   \displaybreak[0]\\
   L_{1,3}&=\frac{\sqrt{10-u} \left(\sqrt{u}
   ((28-3 u) u-76)+2 \sqrt{2} \sqrt{5-u} \sqrt{6-u} \sqrt{16-u}\right)}{64
   \sqrt{21} (u-5)}\,,
   \displaybreak[0]\\
    L_{1,4}&=-\frac{\sqrt{1200-u (3 (u-22)
   u-484)} \sqrt{40-u (3 (u-8)
   u+64)}}{96 \sqrt{14} (u-5)}\,;
   \displaybreak[0]\\
   L_{2,1}&=\frac{\sqrt{16-u} \left(\sqrt{u+6}
   (-3 (u-16) u-160)+18 \sqrt{8-2 u} \sqrt{5-u} \sqrt{10-u}\right)}{192
   \sqrt{35} (u-5)}\,,
   \displaybreak[0]\\
   L_{2,2}&=\frac{\sqrt{6-u} \left(\sqrt{8-2 u} (u (3 u-38)+60)+16
   \sqrt{5-u} \sqrt{-(u-10)
   (u+6)}\right)}{64 \sqrt{105}
   (u-5)}\,,
   \displaybreak[0]\\
   L_{2,3}&=-\frac{\sqrt{u} \left(\sqrt{10-u}
   (u (3 u-32)+96)+2 \sqrt{8-2 u} \sqrt{-(u-5) (u+6)}\right)}{64 \sqrt{21}
   (u-5)}\,,
   \displaybreak[0]\\
   L_{2,4}&=-\frac{\sqrt{u (-3 (u-22)
   u-484)+1200} \sqrt{u (-3 (u-8)
   u-64)+40}}{96 \sqrt{14} (u-5)}\,;
   \displaybreak[0]\\
   L_{3,1}&=\frac{\sqrt{u+6} \left(\sqrt{16-u}
   (20-3 (u-4) u)+18 \sqrt{2} \sqrt{5-u} \sqrt{-(u-6)
   u}\right)}{192 \sqrt{35} (u-5)}\,,
   \displaybreak[0]\\
   L_{3,2}&=\frac{\sqrt{4-u} \left(\sqrt{2} \sqrt{6-u} (u (3 u-22)-20)-16
   \sqrt{5-u} \sqrt{-(u-16)
   u}\right)}{64 \sqrt{105} (u-5)}\,,
   \displaybreak[0]\\
   L_{3,3}&=-\frac{\sqrt{10-u} \left(\sqrt{u}
   (u (3 u-28)+76)+2 \sqrt{2} \sqrt{5-u} \sqrt{6-u} \sqrt{16-u}\right)}{64
   \sqrt{21} (u-5)}\,,
   \displaybreak[0]\\
   L_{3,4}&=-\frac{\sqrt{1200-u (3 (u-22)
   u+484)} \sqrt{40-u (3 (u-8)
   u+64)}}{96 \sqrt{14} (u-5)}\,;
   \displaybreak[0]\\
   L_{4,1}&=-\frac{\sqrt{16-u} \left(\sqrt{u+6} (3 (u-16)
   u+160)+18 \sqrt{8-2 u} \sqrt{5-u} \sqrt{10-u}\right)}{192 \sqrt{35}
   (u-5)}\,,
   \displaybreak[0]\\
   L_{4,2}&=\frac{\sqrt{6-u} \left(\sqrt{8-2 u} (u (3 u-38)+60)-16
   \sqrt{5-u} \sqrt{-(u-10)
   (u+6)}\right)}{64 \sqrt{105}
   (u-5)}\,,
   \displaybreak[0]\\
   L_{4,3}&=\frac{\sqrt{u} \left(\sqrt{10-u}
   ((32-3 u) u-96)+2 \sqrt{8-2 u} \sqrt{-(u-5) (u+6)}\right)}{64 \sqrt{21}
   (u-5)},,
   \displaybreak[0]\\
   L_{4,4}&=-\frac{\sqrt{1200-u (3 (u-22)
   u+484)} \sqrt{40-u (3 (u-8)
   u+64)}}{96 \sqrt{14} (u-5)}\,.
\end{align*}}
\caption{The transformation matrix $L(u)=A_4(u)^{1/2} U(u)$ that we use to
transform the FPU-chain with 4 particles to eigenmodes.} \label{L-expl}
\end{table}

These eigenvectors are the starting point of further computations with {\sc
Mathematica}. We give {\sc Mathematica} many additional substitution rules,
taking into account that $u\in [0,1)$ in the handling of square roots.

The vectors $w_i=A_4(u)^{-1/2}v_i$ are eigenvectors of
$A_2(u)^{1/2}C_4 A_4(u)^{1/2}$. Since the four eigenvalues are different,
the $w_i$ are orthogonal. We take $\tilde w_i=n_i^{-1}w_i$ with
$n_i = \sqrt{w_i\cdot w_i}$ to get an orthonormal basis. There is the
freedom to choose the sign. We multiply $\tilde w_1$ with $-1$, to get
consistency with our earlier computations.

The $\tilde w_i$ can be chosen as the columns of the orthogonal matrix
$U(u)$. Then the vectors
\[\tilde v_i = A_4(u)^{1/2} \tilde w_i= n_i^{-1}v_i(u)\]
are the columns of the transformation matrix $L(u)=A_4(u)^{1/2} U(u)$. In
Table~\ref{L-expl} we give our choice.

The construction of the $v_i$ allows the components to have singularities.
The or\-thonor\-mal\-iza\-tion removes any singularities, so the matrix
elements of $L(u)$ are continuous functions on $[0,u_1)$, given by
algebraic expressions. An explicit expression for the other transformation
matrix $K(u)=A_4(u)^{-1/2} U(u)= A_4(u)^{-1}L(u)$ follows easily.

A check of our computations (including our substitution rules) is carried
out, and gives
\[ K(u)^T A_4(u)K(u)=I_4\,,\qquad L(u)^T C_4 L(u) = \begin{pmatrix}\frac
9{14}&0&0&0\\
0&\frac 27& 0&0\\
0&0&\frac 1{14}&0\\
0&0&0&0\end{pmatrix}\,,\]
in accordance with equation~\eqref{H2trf}.

\end{document}